\documentclass[final,authoryear]{elsarticle}
\usepackage{savesym}
\usepackage{amsmath}
\savesymbol{iint}
\usepackage{txfonts}        % or other package that defines \iint
\restoresymbol{TXF}{iint}
\usepackage{framed,multirow}
\usepackage{doi}

\usepackage{amssymb}
\usepackage{latexsym}
\usepackage{geometry}
\usepackage{url}
\usepackage{xcolor}
\definecolor{newcolor}{rgb}{.8,.349,.1}

\usepackage[title]{appendix}%
\usepackage{booktabs}
\usepackage{bbding} % For checkmarks
\usepackage{pifont}
\usepackage{wasysym}

\usepackage{algorithm}
\usepackage{algorithmicx}
\usepackage{algpseudocode} % For pseudocode styling
\usepackage{slantsc}

\usepackage{diagbox}
\usepackage{xcolor,colortbl}
\usepackage{multirow}
\usepackage{tabu}
\usepackage{epstopdf}
\usepackage{longtable}
\usepackage{caption}
\usepackage{lscape}
\usepackage{graphicx}
\usepackage{latexsym}
\usepackage{wasysym}
\usepackage{rotating}
\usepackage{array}

\usepackage{ifthen}
\usepackage[normalem]{ulem} % for \sout
\usepackage{mathtools}
\newsavebox{\astrutbox}
\sbox{\astrutbox}{\rule[-5pt]{0pt}{20pt}}

\newcommand{\be}{\begin{equation}}
\newcommand{\ee}{\end{equation}}
\newcommand{\bee}{\begin{eqnarray}}
\newcommand{\eee}{\end{eqnarray}}

\newcommand {\delux}[2]{\frac{\partial{#1}}{\partial{#2}}}

\definecolor{Gray}{gray}{0.85}
\definecolor{LightCyan}{rgb}{0.88,1,1}
\definecolor{Red}{rgb}{1,.5,0}
\newcolumntype{g}{>{\columncolor{Gray}}c}
\newcolumntype{y}{>{\columncolor{LightCyan}}c}
\newcolumntype{o}{>{\columncolor{Red}}c}

\newcommand{\bc}{\hat{b}}
\newcommand{\bs}{\tilde{b}}
\newcommand{\Jc}{\hat{J}}
\newcommand{\Js}{\tilde{J}}
\newcommand{\Kc}{\hat{V}}
\newcommand{\Ks}{\tilde{V}}
\newcommand{\kc}{\hat{k}}
\newcommand{\ks}{\tilde{k}}
\newcommand{\uc}{\hat{u}}
\newcommand{\us}{\tilde{u}}

\newcommand{\wc}{\hat{w}}
\newcommand{\ws}{\tilde{w}}
\newcommand{\ec}{\hat{e}}
\newcommand{\es}{\tilde{e}}
\newcommand{\sigmac}{\hat{\sigma}}
\newcommand{\sigmas}{\tilde{\sigma}}
\newcommand{\omegac}{\hat{\omega}}
\newcommand{\omegas}{\tilde{\omega}}

\newcommand{\sigmacj}{\sigmac_j}
\newcommand{\sigmasj}{\sigmas_j}
\newcommand{\omegacj}{\omegac_j}
\newcommand{\omegasj}{\omegas_j}
\newcommand{\kcj}{\kc_j}
\newcommand{\ksj}{\ks_j}
\newcommand{\kci}{\kc_i}

\newcommand{\ucj}{\uc_j}
\newcommand{\usj}{\us_j}

\newcommand{\wcj}{\wc_j}
\newcommand{\wsj}{\ws_j}

\newcommand{\ecj}{\ec_j}
\newcommand{\esj}{\es_j}

\newcommand{\narsimha}[1]{{\color{blue}{#1}}}
\newcommand{\rev}[1]{{\color{black}{#1}}}
\definecolor{rev}{RGB}{0, 0, 0} % Black

\newcommand{\gfsfull}{{Generalized Fourier Series} }
\newcommand{\gfs}{GFS}

\newboolean{showdeleted}
\setboolean{showdeleted}{true} % Change to false to hide deleted text

\newcommand{\deleted}[1]{%
  \ifthenelse{\boolean{showdeleted}}{\textcolor{red}{\sout{#1}}}{}%
  }

\newcolumntype{H}{>{\setbox0=\hbox\bgroup}c<{\egroup}@{}}

\newtheorem{remark}{Remark}
\begin{document}
 
%\verso{Given-name Surname \textit{etal}}
% \verso{N. R. Rapaka, M. K. Riahi}

\begin{frontmatter}

\title{Generalized Fourier Series: An $N\log_{2}(N)$ extension for aperiodic functions that eliminates Gibbs oscillations}

%A novel adaptive spectral method for approximating non-periodic functions without the Gibbs oscillations
%A novel adaptive spectral method for non-periodic domains: annihilation of the Gibbs oscillations

% \tnotetext[label1]{}
\author[label_ku,label_entc]{Narsimha Reddy Rapaka\corref{cor1}}
\author[label_ku,label_entc]{Mohamed Kamel Riahi}
 
% \ead[url]{home page}
% \fntext[label2]{}
\cortext[cor1]{Corresponding author. Tel.: +971 2 3125416, 
\textit{Email address:} narsimha.rapaka@ku.ac.ae (Narsimha R. Rapaka)
%\ead{narsimha.rapaka@ku.ac.ae}
}
% \author[label_kaust]{Ravi \snm{Samtaney}\corref{cor3}}
%\ead{Ravi.Samtaney@kaust.edu.sa}
\address[label_ku]{Department of Mathematics, College of Computing and Mathematical Sciences, Khalifa University of Science and Technology, Abu Dhabi, UAE.}
\address[label_entc]{Emirates Nuclear Technology Center,  Khalifa University of Science and Technology, Abu Dhabi, UAE.}
%\fntext[fn1]{This is author footnote for second author.}  
% \cortext[cor1]{
% %Corresponding author. Tel.: +966 12 8087521. 
% %\textit{E-mail:} \texttt{narsimha.rapaka@kaust.edu.sa} (Narsimha Reddy Rapaka)
% Corresponding author. Tel.: +971 2 3125416. 
% \textit{E-mail:} \texttt{narsimha.rapaka@ku.ac.ae} (Narsimha Reddy Rapaka)
% %\ead{narsimha.rapaka@kaust.edu.sa}
% }
% \cortext[cor3]{Deceased coauthor.}

%\address[1]{Affiliation 1, Address, City and Postal Code, Country}
%\address[2]{Affiliation 2, Address, City and Postal Code, Country}
% \address[label_kaust]{Mechanical Engineering, Division of Physical Science and Engineering, King Abdullah University of Science and Technology, Thuwal, Saudi Arabia}

%\received{1 May 2018}
%\finalform{10 May 2018}
%\accepted{13 May 2018}
%\availableonline{15 May 2018}
%\communicated{R. Samtaney}

\begin{abstract}
%%%
This article introduces the Generalized Fourier Series (GFS), a novel spectral method that extends the classical Fourier series to non-periodic functions. GFS addresses key challenges such as the Gibbs phenomenon and poor convergence in non-periodic settings by decomposing functions into periodic and aperiodic components. The periodic part is represented using standard Fourier modes and efficiently computed via the Fast Fourier Transform (FFT). The aperiodic component employs adaptive, low-rank sinusoidal functions with non-harmonic modes, dynamically tuned to capture discontinuities and derivative jumps across domain boundaries.

Unlike conventional Fourier extension methods, GFS achieves high accuracy without requiring computational domain extensions, offering a compact and efficient representation of non-periodic functions. The adaptive low-rank approach ensures accuracy while minimizing computational overhead, typically involving additional complex modes for the aperiodic part. Furthermore, GFS demonstrates a high-resolution power, with degrees of freedom comparable to FFT in periodic domains, and maintains $N\log_{2}(N)$ computational complexity. The effectiveness of GFS is validated through numerical experiments, showcasing its ability to approximate functions and their derivatives in non-periodic domains accurately. With its robust framework and minimal computational cost, GFS holds significant potential for advancing applications in numerical PDEs, signal processing, machine learning, and computational physics by providing a robust and efficient tool for high-accuracy function approximations.

% Statement on performance advantage?
% discuss the inverse transform?
% show figures including the basis functions: highlihgt the spatially varying amplitude of the basis functions
% does the approximated function converge to the original (exact) function in the limit dx->0?
% how does the modes depend on the smoothness?
% show that as the jump conditions approach zero, the non-periodic modes also approach zero: note that the rank and norm of the jump matrix approaches zero for functions that are smooth at the domain boundaries.

%Present method involves approximation of the derivatives at the end points to compute jumps in the function derivatives which is used to find out the non-integer modes constituting the aperiodic part. It may be possible to find out the non-integer modes with a procedure that does not involve computation of the derivatives. Perhaps, the jumps can be expressed in terms of the nodal values exactly or formulated in terms of a function minimization!
%The key lies in successful decomposition of the original function into periodic and aperiodic parts.
%%%%
\end{abstract}
%\begin{keyword}
%% MSC codes here, in the form: \MSC code \sep code
%% or \MSC[2008] code \sep code (2000 is the default)
%\MSC 41A05\sep 41A10\sep 65D05\sep 65D17
%% Keywords
%\KWD Keyword1\sep Keyword2\sep Keyword3
%\end{keyword}

\end{frontmatter}
% \linenumbers 
%======================================================================
\section{Introduction}
%======================================================================

%%%%%%%%%%%%%%%%%
The accurate approximation of smooth, non-periodic functions within bounded domains is a classical yet persistently challenging problem in numerical analysis. This issue has profound implications for various applications, particularly in the numerical solution of partial differential equations (PDEs). Traditional approaches, such as finite difference, finite volume, finite element methods, and Fourier or Chebyshev spectral methods, have been extensively employed to address this challenge. While non-spectral methods offer broad applicability to problems with arbitrary boundary conditions, they often suffer from limited resolution, requiring a large number of degrees of freedom per wavelength. On the other hand, spectral methods provide superior resolution power but are inherently constrained to periodic domains or homogeneous boundary conditions.  This article aims to broaden the application of Fourier spectral methods to non-periodic functions with minimal computational overhead, focusing solely on the lack of periodicity while assuming sufficient smoothness in the interior.

\rev{The Fast Fourier Transform (FFT) is widely used for approximating smooth, periodic functions due to its efficiency and spectral convergence. However, when applied to non-periodic functions, FFT suffers from the Gibbs phenomenon, producing spurious oscillations near domain boundaries \cite{Hewitt:1979,GottliebS:1997,Jerri1998gibbs}. Over the years, several strategies have been developed to mitigate this limitation. These include filtering techniques to suppress oscillations, domain extension methods that embed non-periodic functions into larger periodic domains \cite{Iserles:2008,Huybrechs:2010,thesis_Adcock:2010,Adcock:2011,Adcock:2014,Geronimo:2020}, polynomial-based periodization transformations \cite{Krylov:1906,Eckhoff:1993,Roache:1978}, signal-processing approaches such as Prony’s method \cite{Prony:1795,Plonka:2018}, and more recently, rational approximation techniques such as the AAA algorithm \cite{Nakatsukasa:2018,Nakatsukasa:2020,HuybrechsT:2023,Discoll:2024}. Notably, Eckhoff's method \cite{Eckhoff:1995} account for multiple isolated discontinuties and its convergence has been analyzed in \cite{Nersessian:2006,Barkhudaryan:2007,Poghosyan:2013}. While these approaches offer partial remedies, they also present significant drawbacks. In particular, domain extension methods increase computational overhead by introducing additional degrees of freedom, involves solution of ill-conditioned systems \cite{AdcockHM:2014}, 
% \sout{may suffer from artificial boundary effects,} 
and often require problem-dependent tuning \cite{Adcock:2014JCP}; moreover, they are not practical in applications where the domain is fixed by data, memory, or geometry constraints.
% (e.g., PDE solvers with prescribed boundary conditions).
% While these approaches offer partial remedies, domain extension methods introduce additional degrees of freedom and computational overhead and may not be practical in applications where the domain is fixed. 
Polynomial-based modifications \cite{Roache:1978,Eckhoff:1995} often lack robustness and deliver limited accuracy for highly oscillatory functions, as demonstrated in this work. Similarly, Prony’s method is ill-conditioned for large system sizes and requires prior knowledge of the number of modes. The AAA algorithm, while powerful and broadly applicable, can incur high computational cost for large-scale problems and may overfit when applied to noisy data \cite{Nakatsukasa:2018}.}

\rev{To tackle these challenges, this work introduces a novel method called the \gfsfull (\gfs), which extends Fourier spectral methods to non-periodic functions with minimal computational overhead. The key idea is to decompose a function into periodic and aperiodic components: the periodic part is efficiently approximated using standard FFT. In contrast, the aperiodic part is represented with an adaptive set of $n=\mathcal{O}(1)$ complex sine and cosine modes.
This strategy avoids artificial domain extensions, preserves the efficiency and simplicity of the Fourier framework, and provides a robust mechanism for capturing both smooth and oscillatory non-periodic behavior. Numerical validations support these results.}

\rev{
Exising methods augment standard Fourier series with a fixed set of functions to approximate non-periodic functions, e.g., polynomials \cite{Roache:1978}, Birkhoff–Hermite polynomials \cite{HybrechsIN:2010}, Bernoulli polynomials \cite{Eckhoff:1995}.
In contrast, GFS adds a dynamic set of non-harmonic sinusoidal modes that adapt to the function's non-periodic nature, efficiently capturing its aperiodic component. We show that GFS is more robust, achieves superior numerically accuracy and convergence. 
% Moreover, by keeping the periodic component $C^{4n-1}$ continuous, the overall method achieves higher accuracy and superior convergence compared to existing domain extension methods~\cite{Iserles:2008,Huybrechs:2010,thesis_Adcock:2010,Adcock:2011, Adcock:2014,Geronimo:2020}.

% The GFS method, like classical Fourier-based methods, is inherently applicable to multi-dimensional problems on tensor-product domains, and we discuss potential extensions to more general geometries.
The GFS method, like Eckhoff’s \cite{Eckhoff:1995,Eckhoff:1998}, is multi-dimensional and flexible, allowing for future expansion to handle multiple discontinuities in complex geometries, such as cubic outer boxes with several immersed obstacles modeled as discontinuities.} 

This paper is organized as follows. Section \ref{sec:motivation} provides the mathematical motivation for the approach, detailing how jump conditions influence the aperiodic component. Section \ref{sec:Approach} outlines the formulation of the \gfsfull and describes the continuous spectrum of aperiodic signals. The adaptive basis construction is then presented, along with the computational algorithms used to compute the modes and expansion coefficients. 
% Section \ref{sec:convergence_aper} addresses the convergence analysis, while 
Section \ref{sec:comp_complexity} outlines the computational complexity of the proposed method. Finally, Section \ref{sec:Numerical_Evidences} demonstrates the method's efficacy through numerical tests on representative non-periodic functions, showcasing the superior resolution power and convergence rate of \gfs~compared to finite difference and the standard FFT methods.

By introducing an adaptive decomposition of the function space and leveraging the efficiency of FFT, the proposed \gfs~method offers a new perspective on handling non-periodic functions in bounded domains. The method retains linear computational complexity, avoids the need for domain extensions, and delivers high-resolution power akin to spectral methods in periodic domains. These features make the \gfs~a compelling alternative for problems involving non-periodic functions, with potential applications in numerical PDEs, signal processing, machine learning, and computational physics.

\rev{
\section{\rev{State-of-the-art methods}}
We compare GFS with the existing methods described briefly below:
\subsection{Eckhoff Method}
The Eckhoff method \cite{Eckhoff:1995} reconstructs a $2\pi$-periodic, piecewise smooth function $u(x)$ from a finite number of its Fourier coefficients $\hat{u}_k$. Assuming $M$ singularities at unknown locations $\gamma_j$ with jumps $A_{j}^{n}\equiv u^{(n)}(\gamma_j^+)-u^{(n)}(\gamma_j^-)$ in the $n$th derivative, the function is decomposed as:
\begin{align}
    u(x) = v(x) + \sum_{n=0}^{q-1} \sum_{j=1}^{M} A_{j}^{n} V_n(x; \gamma_j)\label{eqn:Eckhoff}
\end{align}
where $v(x)$ is smooth and $V_n(x; \gamma_j)$ are periodic functions derived from Bernoulli polynomials with controlled singularities at $\gamma_j$.
\begin{align*}
    V_n(x;\beta) = - \frac{(2\pi)^n}{(n+1)!} B_{n+1}\left(\frac{\xi }{2\pi}\right),    \quad \xi = \mod(x-\beta + 2\pi, 2\pi), 0<\xi<2\pi,
\end{align*}
where, $B_n(x), n=1,2,\hdots,$ are the Bernoulli polynomials \cite{Eckhoff:1995}. Further, $\delux{}{x}V_n(x;\beta) = V_{n-1}(x;\beta)$.
To locate the singularities, one constructs the nonlinear system:
\begin{align*}
    \sum_{n=0}^{q-1} \sum_{j=1}^{M} \frac{A_{j}^{n}}{(ik)^{n}} e^{-ik \gamma_j} = 2\pi i k (\hat{u}_k - \hat{v}_k) \approx 2\pi i k \hat{u}_k,
\end{align*}
for $k=N/2-Mq,N/2-Mq,\hdots,N/2-1$, from which the jump locations and magnitudes are estimated algebraically (see \cite{Eckhoff:1995} for the solution procedure). 
In this work, we consider endpoint singularities only, i.e., $x=-\pi$ (and, by periodicity of $V_n(x;\gamma_j)$, $x=\pi$) so that $M=1$ and $\gamma_1=-\pi$. 

We find that evaluating jumps numerically via Fourier coefficients as in Eckhoff method \cite{Eckhoff:1995} is not robust (also noted in \cite{Eckhoff:1998}) for the functions considered here. Moreover, the Eckhoff method yields only a first-order approximation for the highest derivative jump \cite{Eckhoff:1995}, which can substantially degrade accuracy for highly oscillatory signals. Eckhoff \cite{Eckhoff:1998} examined finite-difference method and asymptotic expansions near discontinuities, while \cite{Barkhudaryan:2007,Poghosyan:2013} studied various choices of Fourier-tail indices to accelerate convergence of jump computations. Nonetheless, their numerical tests used only smooth functions, leaving robustness for highly oscillatory cases unclear.  Here, we evaluate the jumps $A_{1}^{n}\coloneqq u^{(n)}(-\pi)-u^{(n)}(\pi)$ in Eq. \eqref{eqn:Eckhoff} analytically for comparison with GFS. In the GFS method, the end point jumps are denoted $J_n$, related by $J_n=-A_1^n$. 

\subsection{Roache Method}
Roache's method \cite{Roache:1978} handles non-periodic problems using FFTs via ``reduction to periodicity''. The idea is to decompose the target function $f(x)$ into:
\[
f(x) = g(x) + \tilde{f}(x), x\in[-\pi,\pi]
\]
where $g(x)$ is a polynomial that matches the jumps in the boundary derivatives of $f$, denoted with $J_m, m=0,1,\hdots,$, up to order $q-1$, ensuring $\tilde{f}(x)$ becomes continuous at the end points up to $q-1$ derivatives:
\begin{align*}
g(x) &= \sum_{k=0}^{q} a_k x^k, \quad \text{with} \quad \tilde{f}^{(m)}(\pi) - \tilde{f}^{(m)}(-\pi) = 0, \quad 0 \leq m < q,\\
    a_N &= \frac{1}{2\pi} \frac{J_{q-1}}{q!},\\
    a_k &= \frac{1}{2\pi} \left[ \frac{J_{k-1}}{k!} - \sum_{m=k+1}^q a_m \frac{m!}{k! (m-k+1)!} ((\pi)^{m-k+1}-(-\pi)^{m-k+1}) \right].
\end{align*}
Then, FFT is applied to $\tilde{f}(x)$, and $g(x)$ is differentiated analytically.

\subsection{Prony Method}
Prony's method (\cite{Prony:1795}, Chapter 10 of \cite{Plonka:2018}) reconstructs an exponential sum from sampled data $h(k)\equiv h(x_k), x_k=k\Delta x, k=0,1,\hdots,N-1$ with $N\geq 2M$:
\begin{align}
    h(k) = \sum_{j=1}^{M} c_j e^{\phi_j x_k} = \sum_{j=1}^{M} c_j z_j^k, \quad k = 0, \ldots, 2M - 1, \label{eqn:prony}
\end{align}
where, $z_j = e^{\phi_j \Delta x}$ are unknown, distinct complex parameters.

It constructs the Prony polynomial:
\[
p(z) = \prod_{j=1}^{M} (z - z_j) = z^M + p_{M-1}z^{M-1} + \dots + p_0
\]
Solving the linear system:
\[
\sum_{k=0}^{M} p_k h(k + m) = 0, \quad m = 0, \ldots, M-1, \quad p_M=1.
\]
yields $p_k$ and the $z_j$ (thereby $\phi_j$) via roots of the Prony polynomial (equivalent to eigenvalues of the companion matrix of $p(z)$), and then $c_j$ from a Vandermonde system given below (see Chapter 10 of \cite{Plonka:2018} for further details).
\begin{align}
\begin{bmatrix}
1 & 1 & \cdots & 1 \\
z_1 & z_2 & \cdots & z_M \\
% z_1^2 & z_2^2 & z_3^2 & \cdots & z_M^2 \\
\vdots & \vdots & \ddots & \vdots \\
z_1^{M-1} & z_2^{M-1} & \cdots & z_M^{M-1}
\end{bmatrix}
\begin{bmatrix}
    c_1\\
    c_2\\
    \vdots \\
    c_M
\end{bmatrix}
=
\begin{pmatrix}
    h_0\\
    h_1\\
    \vdots\\
    h_{M-1}
\end{pmatrix}
\label{eqn:Vandermonde_Prony}
\end{align}
}
{\color{rev}

% References from Reviewer 1: \cite{Barkhudaryan:2007}

\begin{table}[t]
\centering
\begingroup
\color{rev}
\resizebox{\textwidth}{!}{%
\begin{tabular}{@{}lcHlccccccccc@{}}
\toprule
\rotatebox{-90}{\textbf{Method}} 
  & \rotatebox{-90}{\textbf{Year}} 
  & \rotatebox{-90}{\textbf{Convergence}} 
  & \rotatebox{-90}{\textbf{Complexity}} 
  & \rotatebox{-90}{\textbf{Non-periodic}} 
  & \rotatebox{-90}{\textbf{Resolution (PPW)}} 
  & \rotatebox{-90}{\textbf{Jump Detection}} 
  & \rotatebox{-90}{\textbf{Total DOF}} 
  & \rotatebox{-90}{\textbf{Additional DOF}} 
  & \rotatebox{-90}{\textbf{Splitting}} 
  & \multicolumn{3}{c}{\rotatebox{-90}{\textbf{Aperiodic Recovery}}} \\
\cmidrule(l){11-13}
  &  &  &  &  &  &  &  &  &  
  & \rotatebox{-90}{\textbf{Polynomial}} 
  & \rotatebox{-90}{\textbf{Exponential}} 
  & \rotatebox{-90}{\textbf{Oscillatory}} \\
\midrule
FFT (standard)                    
  & 1965   
  & $\mathcal{O}(e^{-\alpha N})$    
  & $\mathcal{O}(N\log_2 N)$          
  & \ding{53}     
  & 2       
  & \ding{53}      
  & $N$    
  & \ding{53}     
  & \ding{53}       
  & \ding{53} & \ding{53} & \ding{53} \\

GFS  (present)      
  & 2025   
  & $\mathcal{O}(N^{-(q+1)})$      
  & $\mathcal{O}(N\,(q + \log_2 N) + 3(q/4)^3)$  
  & \checkmark    
  & $\sim 3$ 
  & \checkmark     
  & $N + q$ 
  & $q$      
  & dynamic  
  & \checkmark & \checkmark &  \checkmark \\

Eckhoff’s Spectral Reconstruction  \cite{Eckhoff:1995}
  & 1995   
  & $\mathcal{O}(N^{-1})$           
  & $\mathcal{O}(N(q+\log_2 N) + q^3)$          
  & \checkmark    
  & --       
  & \checkmark     
  & $N+q$    
  & $q$ 
  & fixed    
  & \checkmark & \checkmark & \ding{53} \\

Roache’s Jump-Fitting     \cite{Roache:1978}        
  & 2000   
  & $\mathcal{O}(N^{-1})$           
  & $\mathcal{O}(N(q + \log_2 N) + q^3)$   
  & \checkmark    
  & --       
  & \checkmark     
  & $N+q$    
  & $q$ 
  & fixed    
  & \checkmark & \checkmark & \ding{53} \\

Prony’s Method \cite{Prony:1795}
  & 1795   
  & exponential (exact arith.)      
  & $\mathcal{O}(NM + M^3)$              
  & \checkmark    
  & -- 
  & \checkmark     
  & $N+M$ 
  & $M$ 
  & \ding{53}  
  & \ding{53} & \checkmark & \ding{53} \\

Finite Differences (order $6$)    
  & 1945   
  & $\mathcal{O}(N^{-r})$           
  & $\mathcal{O}(N)$                
  & \checkmark    
  & $\mathcal{O}(10^1-10^2)^*$       
  & \ding{53}      
  & $N$    
  & \ding{53}     
  & \ding{53}      
  & \checkmark & \checkmark & \checkmark \\
\bottomrule
\end{tabular}
}
\endgroup
\caption{\rev{Comparison of approximation methods for non-periodic (aperiodic) functions. DOF refers to degrees of freedom, $N$ is the mesh size, $q$ is the number of jumps (in GFS $q=4n$ where $n \sim \mathcal{O}(1)$ is the number of aperiodic modes), $M$ is the number of exponential modes in Prony method with $N\ge 2M$. In Eckhoff's method only one discontinuous location at the domain bounaries is assumed for the complexity estimates. PPW is the points required per wavlength for a target error $\epsilon$. *For a standard sixth-order central finite difference (FD) scheme (one-sided near boundary nodes): PPW $\sim \mathcal{O}(10)$ for $\epsilon=10^{-3}$, PPW $\sim \mathcal{O}(10^2)$ for $\epsilon=10^{-10}$; for GFS: PPW $\sim 3$ for $\epsilon=10^{-10}$ (based on data from Table \ref{tab:complex_exp} and Fig. \ref{fig:convergence_multimode}). 
% Standard reporting consider periodic boundary conditions and when boundary conditions are considered, their accuracy and PPW gets deteriorated further.
}
}
\label{tab:method-comparison}
\end{table}
}
\rev{Table \ref{tab:method-comparison} contrasts GFS with representative numerical methods in computational complexity, degrees of freedom (DOF), scope, and other key features.}
%======================================================================
\section{Generalized Fourier Series (GFS):}\label{sec:Approach}
%======================================================================
%%%%%%%%%%%%%%%%%%%%%%%%%%%%%%%%%%%%%%%%%
\subsection{\rev{Motivation and non-harmonic modes}}\label{sec:motivation}
\begin{figure}[!ht]
 \centering
\begin{minipage}{0.32\linewidth}\centering
 \includegraphics[width=\linewidth,height=2.5in]{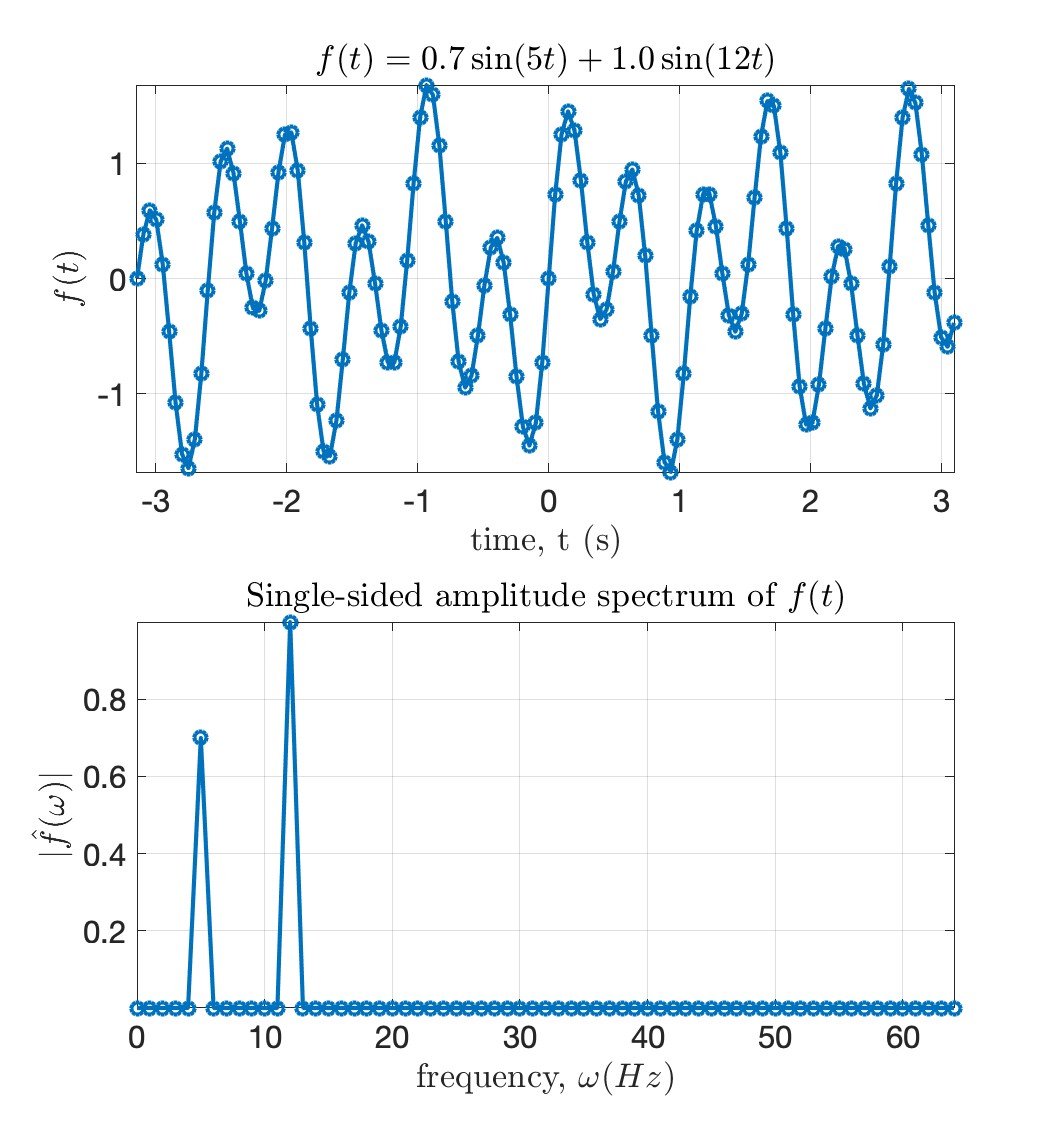} 
 \\(a) \end{minipage}~~
\begin{minipage}{0.32\linewidth}\centering
\includegraphics[width=\linewidth,height=1.25in,trim={0.in 8.in 0.in 0in},clip=true]{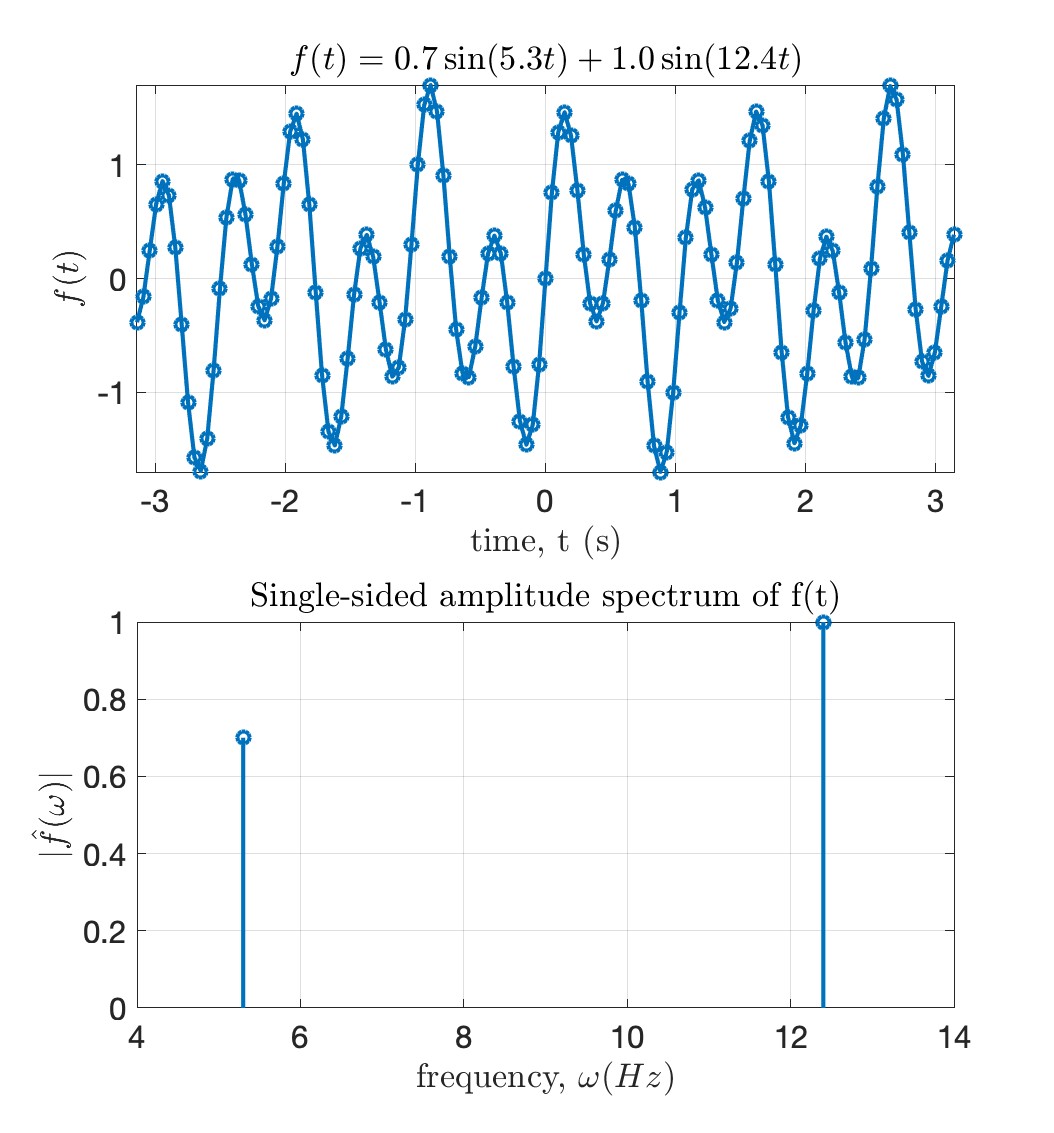}  
\includegraphics[width=\linewidth,height=1.25in,trim={0.in 0.in 0.in 8in},clip=true]{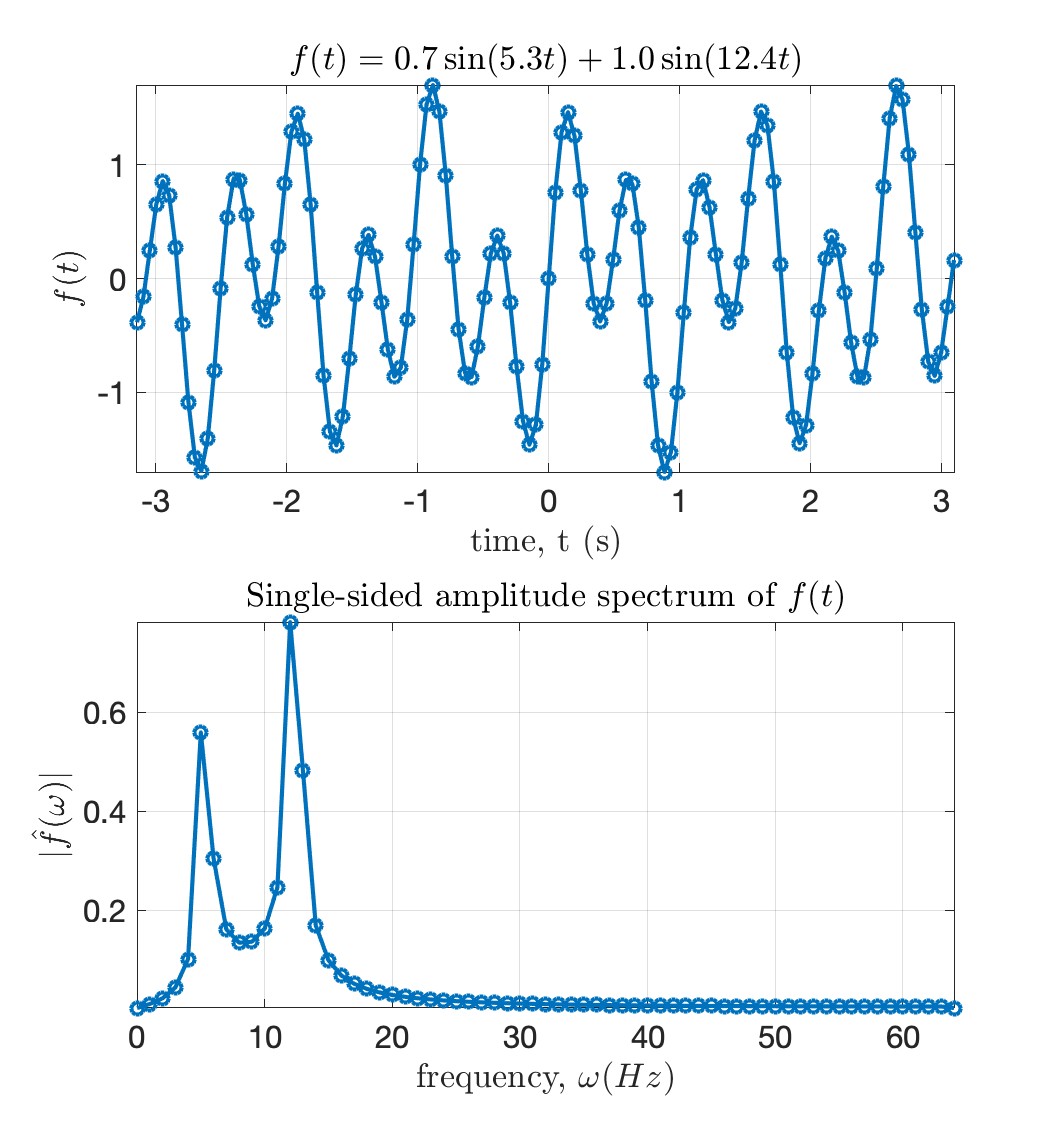}  
\\(b)
 \end{minipage}
\begin{minipage}{0.32\linewidth}\centering
\includegraphics[width=\linewidth,height=2.5in]{./Figures/f_ASM.jpg}  
\\(c)
 \end{minipage}
 \caption{Spectral leakage phenomena: Fast Fourier Transform (FFT) of a function $f(t)=0.7 \sin(2 \pi k_1 t) + \sin(2 \pi k_2 t), 0\le t\le 1$; (a) $f(t)$ is periodic with $k_1=5, k_2=12$, (b) $f(t)$ is non-periodic with $k_1=5.3, k_2=12.4$. Here, the sample frequency, $N=128$. (c) $f(t)$ is non-periodic with $k_1=5.3, k_2=12.4$ and is represented with generalized Fourier series with two non-periodic modes $k_1=5.3, k_2=12.4$.}
 \label{fig:spectral_leakage}
\end{figure}

% Spectral leakage occurs when a signal is truncated in a finite interval that is different from the interval associated with the fundamental mode. The resulting Fourier spectrum spreads energy from one frequency into nearby modes instead of concentrating it at a single mode. 

Spectral leakage occurs when a signal's frequency components are not perfectly aligned with the discrete frequencies of the Fourier transform, such as in the Discrete Fourier Transform (DFT). 
This misalignment results in energy spreading into adjacent frequency bins, creating ``leakage." 
% This misalignment causes energy to spread from one frequency into nearby Fourier modes, rather than being concentrated at a single frequency. 
It is most prominent with non-periodic signals in finite domains.

Fig. \ref{fig:spectral_leakage} illustrates the spectral leakage phenomena with an example. 
Consider a function $f(t)=0.7 \sin(2 \pi k_1 t) + \sin(2 \pi k_2 t), -\pi\le t\le \pi$.
Fig. \ref{fig:spectral_leakage}(a) shows a periodic signal (top) comprised of discrete harmonics $k_1=5, k_2=12$ of the fundamental frequency corresponding to the domain. The corresponding spectrum (bottom) exhibits discrete spikes at $k_1, k_2$. This is expected because the Fourier spectrum of a periodic signal comprises a discrete set of modes corresponding to harmonics (integer multiples) of a fundamental mode $2\pi/L$ associated with the domain length $L$. 
Fig. \ref{fig:spectral_leakage}(b) shows a non-periodic signal (top) comprised of modes $k_1=5.3, k_2=12.4$. The corresponding spectrum (bottom) is spread over nearby modes, yielding a side-lobe structure instead of discrete spikes at $k_1=5.3, k_2=12.4$. This is because the modes $k_1=5.3, k_2=12.4$ are not harmonics of the fundamental frequency ($k=1$). This phenomenon is known as spectral leakage. Fig. \ref{fig:spectral_leakage}(c) shows the non-periodic signal (top) comprised of modes $k_1=5.3, k_2=12.4$ represented by the present method (\gfs) whose spectrum (bottom) is sharp and composed merely of two discrete non-integer modes, same as those in the original function. The \gfs~algorithm described in later sections automatically detects these mode numbers.

\rev{
Spectral leakage primarily arises because the Fourier spectrum is sampled over integer multiples of a fundamental frequency (harmonic modes). By incorporating non-harmonic modes, we can recover a sharp spectrum, but identifying these modes without undue computational cost is challenging.
}

\rev{
The classical Prony method \cite{Prony:1795}\cite{Plonka:2018} can detect non-harmonic modes in noiseless, sampled data, but it requires prior knowledge of the number of modes $M$. Since $M$ is arbitrary, the resulting system (e.g. the Vandermonde system in Eq. \eqref{eqn:Vandermonde_Prony}) can become highly ill-conditioned when $M\gg 1$ and its computational cost scales as $\mathcal{O}(N^3)$ with mesh size $N$. 
}

\rev{
To reduce computational overhead and improve conditioning, we propose GFS, which decomposes the signal into periodic and non-periodic components. This approach, however, relies on knowledge of jump conditions, i.e., prior information about jumps or a numerical method to evaluate them.
}
\subsection{The continuous spectrum for an aperiodic signal}
\rev{In the standard Fourier framework, aperiodic signals exhibit a continuous spectrum due to their non-repeating nature, whereas periodic signals consist of discrete harmonics (integer multiples of a fundamental frequency).
As a result, the Fourier representation for aperiodic signals shifts from a sum of discrete frequencies to an integral over a continuous frequency domain. 

Our framework handles both periodic and aperiodic signals using non-harmonic modes that capture complex features, including real frequencies for sinusoidal behavior and imaginary frequencies for growth or decay. By employing non-harmonic modes, we can represent the signal with a finite set of real-valued frequencies, eliminating spectral leakage, as illustrated in Fig.\ref{fig:spectral_leakage}.

We decompose a function $u(x)$ on $[-\pi, \pi]$ as $u(x) = u_{p}(x) + u_a(x)$, where $u_{p}(x)$ is periodic and $u_a(x)$ is aperiodic.
The same analysis applies to any domain $x\in[a,b]$; we use the coordinate transformation $x^*=2\pi (x-x_0)/(b-a)$, with $x_0=(a+b)/2$, and for simplicity, take $a=-\pi, b=\pi$.
We split \( u_a(x) \) into symmetric (even) \( u_c(x) \) and antisymmetric (odd) \( u_s(x) \) components, using non-harmonic mode functions \(\cos(\kcj x)\) and \(\sin(\ksj x)\), where \(\kcj, \ksj \in \mathbb{C}\). The symmetric part's cosine modes are denoted with a hat (`$\widehat\quad$'), while the anti-symmetric part's sine modes are denoted with a tilde (`$\sim$').}
\begin{eqnarray}
u_c(x) &=& \sum_{j=1}^{n_c} \ucj \cos(\kcj x), \quad\quad\quad\quad
u_s(x) = \sum_{j=1}^{n_s} \usj \sin(\ksj x), \label{def:uc_us}\\
u_a(x) &=& u_c(x) + u_s(x) = \underbrace{\sum_{j=1}^{n_c} \ucj \cos(\kcj x)}_{\text{symmetric}} + \underbrace{\sum_{j=1}^{n_s}\usj \sin(\ksj x)}_{\text{anti-symmetric}}.
\label{def:ua}
\end{eqnarray}

Let's assume a set of complex modes $k_j=\sigma_j+i\omega_j$, where $\sigma_j,\omega_j \in \mathbb{R}_{+}$ and $i=\sqrt{-1}$. The Euler formula clearly illustrates the decay/growth behavior of the mode $k_j$.
\begin{align*}
e^{i k_{j} x} &= \left(\cos(\sigma_{j} x) + i\sin(\sigma_{j} x) \right)e^{-\omega_{j} x}, \\
e^{-i k_{j} x} &= \left(\cos(\sigma_{j} x) - i\sin(\sigma_{j} x) \right)e^{\omega_{j} x}.
\end{align*}
For $\kc_{j}=\sigmac_{j}+i\omegac_{j}$, we set the following
decay and growth representation for $\cos$ ``symmetric'' series, with complex coefficients $(\ucj)$:
\[
u_{c}(x)= \dfrac{1}{2}\sum_{j} \ucj \left( e^{i\kcj x} + e^{-i\kcj x} \right) = \sum_{j} \ucj \cos(\kcj x),
\]
similarly, for $\ks_{j}=\sigmas_{j}+i\omegas_{j}$, a decay and growth representation for $\sin$ ``anti-symmetric'' series, with complex coefficients $(\usj)$:
\[
u_{s}(x)= \dfrac{1}{2i}\sum_{j} \usj \left( e^{i\ksj x} - e^{-i\ksj x} \right) = \sum_{j} \usj \sin(\ksj x)
\]

\rev{The above certifies our use of $\cos(\kc_{j} x)$ and $\sin(\ks_{j} x)$, where the complex-valued modes $\kc_{j}, \ks_{j}$ introduce exponential decay or growth alongside oscillations, adapting to the function being approximated.

% Using a classical Fourier representation, the continuous spectrum is essential for transient or localized phenomena, such as decaying pulses, where energy spreads smoothly across frequencies rather than concentrating at specific harmonics. In contrast, our method employs complex-valued frequencies, enabling the effective capture of both periodic and aperiodic phenomena without requiring a smooth continuous spectrum. It precisely identifies the frequencies suitable for periodic signals and their aperiodic counterparts, including growth- or decay-driven signals.
}

\rev{Although represented by complex modes and amplitudes, $u_c(x)$ and $u_s(x)$ in Eq. \eqref{def:uc_us} are always real (see Table \ref{tab:usuc_kskc} and Sec. \ref{sec:solution_modes} for details) and need not have the same number of modes, $n_c$ and $n_s$, respectively. For convenience, we set $n_c=n_s=n$ without loss of generality; if $n_c\neq n_s$, we take $n=\max(n_c,n_s)$ for the complexity analysis of the numerical algorithm described in Sec. \ref{sec:algorithm}.}

\subsection{Jump's approach for aperiodic signal}
Let us denote the $m^{th}$ derivative of $u(x)$ as $u^{(m)}$ and assume that $u(x)$ is smooth in the interior. We have $u^{(m)} = u_{p}^{(m)} + u_a^{(m)}$. 
%Let $[u]\equiv u(\pi)-u(-\pi)$ be the jump in $u$ at the end points $x_i=\pm \pi$ and $J_m\equiv [u^{(m)}]$. 
Let us 
% denote the jump in $u$ across the end points as $[u] \equiv u(\pi) - u({-\pi})$ and
denote the jump in the $m^{th}$ derivative of $u$ as $J_m\equiv [u^{(m)}] $. By the definition of periodicity, we have $[u_{p}^{(m)}]=0$ which yields $[u^{(m)}] = [u_a^{(m)}]$ for $m=0,1,2,3,...,\infty$.
Note that the jump conditions may be known in certain applications; otherwise, they are approximated numerically as described in Sec. \ref{sec:fd}.

The even derivatives of $u_c(x)$ and odd derivatives of $u_s(x)$ are symmetric and, therefore, have zero jumps at the endpoints. So, we have $J_{2m}\equiv [u^{(2m)}] =[u_a^{(2m)}] = [u_s^{(2m)}]$ and $J_{2m+1}\equiv [u^{(2m+1)}] =[u_a^{(2m+1)}] = [u_c^{(2m+1)}]$ for $m=0,1,2,...,\infty$.
\begin{align}
%u_a(x) &\equiv \sum_{j=1}^n \ucj \cos(\kcj x) + \sum_{j=1}^n\usj \sin(\ksj x),\\
%u_a^{(2m)} &= \sum_j (-1)^{(m)} \kcj^{(2m)}\ucj \cos(\kcj x) +\sum_j (-1)^{(m)} \ksj^{2m} \usj \sin(\ksj x),\quad\quad\quad\quad\quad~{\text{ (even derivatives)}}\\
%u_a^{(2m+1)} &= \sum_j (-1)^{(m+1)} \kcj^{(2m+1)}\ucj \sin(\kcj x) + \sum_j (-1)^{(m)} \ksj^{(2m+1)}\usj \cos(\ksj x),\quad\quad~~{\text{ (odd derivatives)}}\\
%J_0 = [u_a] &= \sum_{j=1}^n 2\usj \sin(\ksj \pi),\quad \quad J_1 = \left[ u'_a \right] = \sum_{j=1}^n -2\kcj\ucj \sin(\kcj \pi), \\
J_{2m} = \left[u_s^{(2m)}\right] &=  \sum_{j=1}^n (-1)^m 2\ksj^{2m} \usj \sin(\ksj \pi), \quad\quad ~~~~ m=0,1,2,...,\infty, 
% ~{\text{ (jump in even derivatives of $u(x)$)}} 
\label{eqn:J2m}
\\ J_{2m+1} = \left[u_c^{(2m+1)}\right] &= \sum_{j=1}^n (-1)^{(m+1)} 2\kcj^{(2m+1)}\ucj \sin(\kcj \pi), ~~ m=0,1,2,...,\infty. 
% ~\text{ (jump in odd derivatives of $u(x)$)} 
\label{eqn:J2mp1}
%\\ m&=&0,1,2,...,\infty. \nonumber
%\\ m&=&0,1,2,...,2n-1. \nonumber
\end{align}
The jumps computation can lead to rank deficiency; thus, we employ a low-rank approximation with $n\sim \mathcal{O}(1)$, demonstrating its effectiveness in practical applications in Sec. \ref{sec:Numerical_Evidences}. This approach improves computational efficiency and ensures optimal approximation when jumps are involved.  

\rev{
Our GFS approach can be summarized  into  three steps: 
\begin{itemize}
\item First, we estimate the aperiodic part from derivative jumps at the domain boundaries. \item Second, we recover the periodic part by subtracting the aperiodic part from the original signal. 
\item Third, we approximate the derivative of the original signal as the sum of the derivatives of its periodic and aperiodic parts.
\end{itemize}
In the sequel, we detail the first step—identifying the aperiodic part via non-harmonic modes—while the second step is trivial (though necessary), and the third step is readily obtained through classical Fourier analysis.
}

\subsection{\rev{Dynamic non-harmonic modes $\ks_j, \kc_j$ and amplitudes $\us_j, \uc_j$}}
\subsubsection{\rev{Formulation of the system of equations}}
With $2n$ unknowns ($\ksj, \usj$) in Eq. \eqref{eqn:J2m}, we derive a system of $2n$ equations by evaluating Eq. \eqref{eqn:J2m} for $m=0,1,2,...,2n-1$. Similarly, another system of $2n$ equations for the unknowns ($\kcj, \ucj$) is obtained from Eq. \eqref{eqn:J2mp1} for $m=0,1,2,...,2n-1$. Both systems are presented in matrix form below.

%\begin{minipage}{0.48\linewidth}
\begin{gather}
 2
  \begin{pmatrix}
  1		&	1	& \hdots 	&	1 \\
  -\ks_1^2 & -\ks_2^2 & \hdots 	& -\ks_n^2 \\
  \ks_1^4 & \ks_2^4 & \hdots 	& \ks_n^4 \\
\vdots & \vdots & \vdots 	& \vdots \\
    (i\ks_1)^{4n-2} &   (i\ks_2)^{4n-2} & \hdots 	&   (i\ks_n)^{4n-2} 
 \end{pmatrix}
\begin{pmatrix} 
\us_1 \sin(\ks_1\pi) \\ 
\us_2 \sin(\ks_2\pi) \\
\vdots \\
\us_n \sin(\ks_n\pi) 
 \end{pmatrix}
 = 
  \begin{bmatrix}
   J_0 \\ J_2\\J_4
\\   \vdots \\J_{4n-2}
   \end{bmatrix}, \nonumber
\end{gather}
%\end{minipage}~~
%\begin{minipage}{0.48\linewidth}
\begin{gather}
 -2
 \begin{pmatrix}
  1		&	1	& \hdots 	&	1 \\
  -\kc_1^2 & -\kc_2^2 & \hdots 	& -\kc_n^2 \\
  \kc_1^4 & \kc_2^4 & \hdots 	& \kc_n^4 \\
\vdots & \vdots & \vdots 	& \vdots \\
  (i\kc_1)^{4n-2} &   (i\kc_2)^{4n-2} & \hdots 	&   (i\kc_n)^{4n-2}
 \end{pmatrix}
 \begin{pmatrix} 
\uc_1 \kc_1\sin(\kc_1\pi) \\ 
\uc_2 \kc_2\sin(\kc_2\pi) \\
\vdots \\
\uc_n \kc_n\sin(\kc_n\pi) 
 \end{pmatrix}
 = 
  \begin{bmatrix}
   J_1 \\ J_3\\J_5
\\ \vdots \\J_{4n-1}
   \end{bmatrix}. \nonumber
\end{gather}
%\end{minipage}
%\vspace{0.1in}
\\
The above equations can be rearranged to obtain a set of equations involving only the unknowns $\kcj$ and $\ksj$ as below,

\begin{gather}
 \underbrace{2\begin{pmatrix} 
\us_1 \sin(\ks_1\pi) \\ 
\us_2 \sin(\ks_2\pi) \\
\vdots \\
\us_n \sin(\ks_n\pi) 
 \end{pmatrix}}_{\ws}
 = 
 {\underbrace{\begin{pmatrix}
  1		&	1	& \hdots 	&	1 \\
  -\ks_1^2 & -\ks_2^2 & \hdots 	& -\ks_n^2 \\
\vdots & \vdots & \vdots 	& \vdots \\
 (i\ks_1)^{2n-2} &  (i\ks_2)^{2n-2} & \hdots 	&  (i\ks_n)^{2n-2} 
 \end{pmatrix}}_{\equiv \Ks}}^{-1}
  \underbrace{\begin{bmatrix}
   J_0 \\ J_2\\\vdots\\J_{2n-2}
   \end{bmatrix}}_{\bs}
 = 
 \begin{pmatrix}
  (i\ks_1)^{2n} & (i\ks_2)^{2n} & \hdots 	& (i\ks_n)^{2n} \\
  (i\ks_1)^{2n+2} & (i\ks_2)^{2n+2} & \hdots 	& (i\ks_n)^{2n+2} \\
\vdots & \vdots & \vdots 	& \vdots \\
  (i\ks_1)^{4n-2} & (i\ks_2)^{4n-2} & \hdots 	& (i\ks_n)^{4n-2} 
 \end{pmatrix}^{-1}
  \begin{bmatrix}
   J_{2n} \\ J_{2n+2}\\\vdots\\J_{4n-2}
   \end{bmatrix}. 
   \label{eqn:ujkj_sin}
\end{gather}
\begin{gather}
 \underbrace{-2\begin{pmatrix} 
\uc_1 \kc_1\sin(\kc_1\pi) \\ 
\uc_2 \kc_2\sin(\kc_2\pi) \\
\vdots \\
\uc_n \kc_n\sin(\kc_n\pi) 
 \end{pmatrix}}_{\wc}
 = 
  {\underbrace{\begin{pmatrix}
  1		&	1	& \hdots 	&	1 \\
  -\kc_1^2 & -\kc_2^2 & \hdots 	& -\kc_n^2 \\
\vdots & \vdots & \vdots 	& \vdots \\
  (i\kc_1)^{2n-2} & (i\kc_2)^{2n-2} & \hdots 	& (i\kc_n)^{2n-2} 
 \end{pmatrix}}_{\equiv \Kc}}^{-1}
  \underbrace{\begin{bmatrix}
   J_1 \\ J_3\\\vdots\\J_{2n-1}
   \end{bmatrix}}_{\bc}
 = 
 \begin{pmatrix}
  (i\kc_1)^{2n} & (i\kc_2)^{2n} & \hdots 	& (i\kc_n)^{2n} \\
  (i\kc_1)^{2n+2} & (i\kc_2)^{2n+2} & \hdots 	& (i\kc_n)^{2n+2} \\
\vdots & \vdots & \vdots 	& \vdots \\
  (i\kc_1)^{4n-2} & (i\kc_2)^{4n-2} & \hdots 	& (i\kc_n)^{4n-2} 
 \end{pmatrix}^{-1}
  \begin{bmatrix}
   J_{2n+1} \\ J_{2n+3}\\\vdots\\J_{4n-1}
   \end{bmatrix}. 
   \label{eqn:ujkj_cos}
\end{gather}
The expressions satisfy $[u^{(m)}]=[u_a^{(m)}]$, leading to $[u_{p}^{(m)}]=0$ for $m=0,1,2,\hdots,4n-1$. This means that the periodic part is smooth (continuous at the endpoints) up to $(4n-1)$ derivatives, so $u_{p}(x) \in C^{4n-1}$. Each mode $\ksj, j=1,2,\hdots,n,$ (respectively $\kcj$) must be unique (or simply separate) for the inverse of the transposed Vandermonde matrix $\Ks$ ( respectively $\Kc$) to exist. The unknowns of the even part ($\kc_j$, $\uc_j$) depend solely on the jumps relative to odd derivatives of $u(x)$, while those of the odd part ($\ks_j$, $\us_j$) rely only on the jumps relative to the even derivatives. Additionally, the second and third components in Eqs. \eqref{eqn:ujkj_sin}-\eqref{eqn:ujkj_cos} involving $\ksj$ and $\kcj$ are identical, with the even jumps $J_0, J_2,\hdots, J_{4n-2}$ in Eq. \eqref{eqn:ujkj_sin} replaced by the odd jumps $J_1, J_3,\hdots, J_{4n-1}$ in Eq. \eqref{eqn:ujkj_cos}. 
Thus, the solution methods for $\ksj$ and $\kcj$ are analogous.

We solve the above non-linear systems of equations to find the unknown modes ($\kcj, \ksj$) and their corresponding amplitudes ($\ucj, \usj$) in terms of the jumps $J_m$. We first consider the cases for $n=1,2,3$ in Appendix \ref{sec:proof} and then generalize the method for arbitrary $n$ below. This generalization lacks a rigorous proof due to the complexity of the algebra with larger $n$, but we validate it in Sec. \ref{sec:Numerical_Evidences} with multiple test cases.

\subsubsection{\rev{Solution method}}\label{sec:solution_modes}%: $u_{p}$ is $C^{(4n-1)}$--continuous}
First, let us define the elementary symmetric polynomials ($\es_j$) in $\ks_1^2, \ks_2^2, \hdots, \ks_n^2$ as,
\begin{eqnarray*}
\sum_{1\leq j\leq n} \ksj^2 = \es_1, 
\quad \sum_{1< j<k\leq n} \ksj^2 \ks_k^2 = \es_2,
\quad \hdots, 
%\vdots \quad~~\\%\sum_{1< j<k<l\leq n} \ksj^2 \ks_k^2 k_{il}^2 = \es_3
\quad \sum_{1< j_1<j_2<\cdots<j_m\leq n} \ks_{j_1}^2 \ks_{j_2}^2 \cdots \ks_{j_m}^2 = \es_m,
\end{eqnarray*}
and the elementary symmetric polynomials ($\ec_j$) in $\kc_1^2, \kc_2^2, \kc_3^2$ as,
\begin{eqnarray*}
\sum_{1\leq j\leq n} \kcj^2 = \ec_1, 
\quad \sum_{1< j<k\leq n} \kcj^2 \kc_k^2 = \ec_2,
\quad \hdots, 
\quad \sum_{1< j_1<j_2<\cdots<j_m\leq n} \kc_{j_1}^2 \kc_{j_2}^2 \cdots \kc_{j_m}^2 = \ec_m.
\end{eqnarray*}

Analogous to Eqs. \eqref{eqn:es1ec1}, \eqref{eqn:es2}-\eqref{eqn:ec2}, and \eqref{eqn:es_3}-\eqref{eqn:ec_3} for $n=1,2$ and $3$, respectively, a set of equations for the elementary symmetric polynomials $\es_j, \ec_j$ for arbitrary number ($n$) of modes are given below which satisfy the second equality in Eqs. \eqref{eqn:ujkj_sin}-\eqref{eqn:ujkj_cos}, respectively,

\begin{gather}
  \begin{bmatrix}
   \es_n \\ \vdots \\ \es_2 \\ \es_1
   \end{bmatrix}
 =
 -
 {\underbrace{\begin{pmatrix} 
 J_0 & J_2 & \hdots &  J_{2(n-1)} \\ 
 J_2 & J_4 & \hdots & J_{2(n)} \\
 \vdots &  \vdots & \ddots &  \vdots \\
 J_{2(n-1)} & J_{2n} & \hdots &  J_{2(2n-2)} 
 \end{pmatrix}}_{\equiv \Js}}^{\dagger}
  \begin{bmatrix}
   J_{2n} \\ J_{2(n+1)}\\\vdots \\J_{2(2n-1)}
   \end{bmatrix},
   \label{eqn:es_n}
\end{gather}
% ===================
\begin{gather}
  \begin{bmatrix}
   \ec_n \\\vdots \\ \ec_2 \\ \ec_1
   \end{bmatrix}
 =
 -
 {\underbrace{\begin{pmatrix} 
 J_1 & J_3 &	\hdots & J_{2n-1} \\ 
 J_3 & J_5 &	\hdots & J_{2n+1} \\
 \vdots &  \vdots & \ddots &  \vdots \\
 J_{2n-1} & J_{2n+1} &	\hdots & J_{4n-3} 
 \end{pmatrix}}_{\equiv \Jc}}^{\dagger}
  \begin{bmatrix}
    J_{2n+1} \\ J_{2n+3}\\\vdots\\ J_{4n-1}
   \end{bmatrix}. 
   \label{eqn:ec_n}
\end{gather}
A straightforward expansion shows that $\ks_1^2, \ks_2^2, ...\ks_n^2$ and $\kc_1^2, \kc_2^2, ...\kc_n^2$ are the roots ($\lambda$) of the polynomials below involving $\es_j, \ec_j$, respectively,
\begin{align}
\prod_{j=1}^n (\lambda-\ksj^2) &= \lambda^n - \es_1 \lambda^{(n-1)} +\hdots+ (-1)^k \es_k \lambda^{(n-k)} +\hdots + (-1)^n \es_n = 0 , 
\label{eqn:ks_characteristic}
\\\prod_{j=1}^n (\lambda-\kcj^2) &= \lambda^n - \ec_1 \lambda^{(n-1)} +\hdots+ (-1)^k \ec_k \lambda^{(n-k)} +\hdots + (-1)^n \ec_n = 0.
\label{eqn:kc_characteristic}
\end{align}
Note that $\kc_j = \pm \sqrt{\lambda_j}$ and $\ks_j = \pm \sqrt{\lambda_j}$ and the $\pm$ sign is irrelevant here because $\uc_j \cos(\kc_j x)$ and $\us_j \sin(\ks_j x)$ (therefore, $u_c(x)$ and $u_s(x)$) are even functions of $\kc_j$ and $\ks_j$, respectively. Here, we consider the complex square root defined as $\sqrt{z}\equiv \sqrt{|z|}e^{i\theta/2}$ where $z=|z|e^{i\theta}, -\pi\le \theta\le \pi$. Once the modes $\ks_j, \kc_j$ are obtained by solving for the roots of the above polynomials, the amplitudes $\us_j, \uc_j$ can be obtained from the first equality in Eqs. \eqref{eqn:ujkj_sin} and \eqref{eqn:ujkj_cos}, respectively. Note that the matrices $\Js, \Jc$ in Eqs. \eqref{eqn:es_n}, and \eqref{eqn:ec_n} may be rank-deficit and hence non-invertable. 
For example, 
%if $u(x)$ is nearly periodic and contain only two aperiodic modes then rank($\Js$)=rank($\Jc$)=2 and we 
a purely periodic function $u(x)$ has all the jumps ($J_m, \forall j$) equal to zero and therefore rank($\Js$)=rank($\Jc$)=0. In practice, one must use the Moor-Penrose generalized inverse, hence the notation of $\Js^\dagger$.  

% In practice, if determinant of $\Js,\Jc=0$ then we set $n_s=\text{rank}(\Js), n_c=\text{rank}(\Jc)$, respectively.

\begin{remark}\label{remark:modes}
Since $\kcj^2$ and $\ksj^2$ are roots of polynomials with
real coefficients, they must be real numbers or form complex conjugate pairs, according to the complex conjugate root theorem. So, the modes $\kcj$ and $\ksj$ must be real (if $\kcj^2, \ksj^2 \in \mathbb{R}_+$) or purely imaginary (if $\kcj^2, \ksj^2 \in \mathbb{R}_-$, $\operatorname{Re}(\kcj)=\operatorname{Re}( \ksj) = 0$) or form complex conjugate pairs (if $\kcj^2, \ksj^2 \in \mathbb{C}, \kcj^2 =\kci^{*2}, \ksj^2 =\ksj^{*2}$  for some $(i,j)$).

% Let $\{k_j\}_{j=1}^n$ denote a set of modes that can be separated into three disjoint subsets:
% The entire set of modes can thus be written as $\{k_j\}_{j=1}^n = \mathcal{R} \cup \mathcal{I} \cup \mathcal{C}$. 
   % \[
   % k_j \in \mathcal{R}; \quad
   % k_j \in \mathcal{I}; \quad
   % k_j \in \mathcal{C}, \text{with } k_q = {k_j^*}.
   % \]
Let us define three disjoint sets $\mathcal{R}, \mathcal{I}$, and $\mathcal{C}$ comprised solely of real numbers, purely imaginary numbers, and complex conjugate pairs, respectively, as follows,
\begin{align*}
   \mathcal{R} &= \{k_j \in \mathbb{R} \mid \operatorname{Im}(k_j) = 0\}, \hspace{1.38in} \text{(real)}\\
   \mathcal{I} &= \{k_j \in i\mathbb{R} \mid \operatorname{Re}(k_j) = 0, \operatorname{Im}(k_j) \neq 0\},\hspace{0.61in} \text{(purely imaginary)}\\
   \mathcal{C} &= \{(k_i, k_j) \in \mathbb{C} \times \mathbb{C} \mid k_j = {k_i^*}, \operatorname{Im}(k_j) \neq 0\}, \hspace{0.32in} \text{(complex conjugate pairs)}
\end{align*}
If $k_j\in \mathcal{I}$, then $k_j^2 \in \mathcal{R}$.
Let us separate $\ksj^2$ (and $\kcj^2$) into the following subsets, 
$$
    \{\ksj^2\}_{j=1}^n = \{(\ks_p^2, \ks_q^2)\}_{\mathcal{C}}\cup \{\ks_r^2\}_{\mathcal{R}},
$$
where the integers $p,q,r\in \mathbb{I}; 1\leq p\leq \tilde{n}_p; q=\tilde{n}_p+p$, and $ 2\tilde{n}_p+1\leq r\leq n$.  And
$$
    \{\kcj^2\}_{j=1}^n = \{(\kc_p^2, \kc_q^2)\}_{\mathcal{C}}\cup \{\kc_r^2\}_{\mathcal{R}}, 
$$
where the integers $p,q,r\in \mathbb{I}; 1\leq p\leq \hat{n}_p; q=\hat{n}_p+p$, and $ 2\hat{n}_p+1\leq r\leq n$. 

% \begin{align*}
%     \{\ksj\}_{j=1}^n &= \{\ks_p\}_{p=1}^{\tilde{n}_p} \cup \{\ks_q\}_{q=1}^{\tilde{n}_q} \cup \{\ks_r\}_{r=1}^{\tilde{n}_r} \cup \{\ks_s\}_{s=1}^{\tilde{n}_s}, \tilde{n}_p+\tilde{n}_q+\tilde{n}_r+\tilde{n}_s=n, \tilde{n}_p=\tilde{n}_q,  (\ks_p, \ks_q) \in \mathcal{C}, \ks_p=\ks_q^{*}, \ks_r \in \mathcal{R}, \ks_s \in \mathcal{I},\\
%     \{\kcj\}_{j=1}^n &= \{\kc_p\}_{p=1}^{\hat{n}_p} \cup \{\kc_q\}_{q=1}^{\hat{n}_q} \cup \{\kc_r\}_{r=1}^{\hat{n}_r} \cup \{\kc_s\}_{s=1}^{\hat{n}_s}, \hat{n}_p+\hat{n}_q+\hat{n}_r+\hat{n}_s=n, \hat{n}_p=\hat{n}_q,  (\kc_p, \kc_q) \in \mathcal{C}, \kc_p=\kc_q^{*}, \kc_r \in \mathcal{R}, \kc_s \in \mathcal{I},\\
%     \{\ksj^2\}_{j=1}^n &= \{\ks_p^2\}_{p=1}^{\tilde{n}_p} \cup \{\ks_q^2\}_{q=1}^{\tilde{n}_q} \cup \{\ks_r^2\}_{r=1}^{\tilde{n}_r+\tilde{n}_s}, \tilde{n}_p+\tilde{n}_q+\tilde{n}_r+\tilde{n}_s=n,  \tilde{n}_p=\tilde{n}_q, (\ks_p^2, \ks_q^2) \in \mathcal{C}, \ks_p^2=\ks_q^{*2}, \ks_r^2 \in \mathcal{R},\\
%     \{\kcj^2\}_{j=1}^n &= \{\kc_p^2\}_{p=1}^{\hat{n}_p} \cup \{\kc_q^2\}_{q=1}^{\hat{n}_q} \cup \{\kc_r^2\}_{r=1}^{\hat{n}_r+\hat{n}_s}, \hat{n}_p+\hat{n}_q+\hat{n}_r+\hat{n}_s=n,  \hat{n}_p=\hat{n}_q, (\kc_p^2, \kc_q^2) \in \mathcal{C}, \kc_p^2=\kc_q^{*2}, \kc_r^2 \in \mathcal{R}.
% \end{align*}

Let us rearrange Eq. \eqref{eqn:ujkj_sin} and Eq. \eqref{eqn:ujkj_cos}, respectively as follows,

    \[\sum_{j=1}^n \Ks_{ij} \wsj = \bs_i, \text{ and }\quad
    \sum_{j=1}^n \Kc_{ij} \wcj = \bc_i, \]
    where 
\begin{align}    
    \wsj &= 2 \usj \sin(\ksj \pi),\label{eqn:wsj}\\
    \wcj &= -2 \ucj \kcj \sin(\kcj \pi).\label{eqn:wcj}
\end{align}
%where, $\wsj = 2 \usj \sin(\ksj \pi)$ and $\wcj = -2 \ucj \kcj \sin(\kcj \pi)$.
Since the $j^{th}$ columns of the transposed Vandermonde matrices $\Ks$ and $\Kc$ depend solely on $\ksj^2$ and $\kcj^2$, i.e., $\Ks_{ij}=\Ks_{ij}(\ksj^2)$ and $\Kc_{ij}=\Kc_{ij}(\kcj^2)$, we can rearrange above equations as,
% Let us denote the $j^{th}$ columns of $\Ks$ and $\Kc$ with $\vsj(\ksj^2)$ and $\vcj(\kcj^2)$, respectively. Since $\vsj(\ksj^2)$ and $\vcj(\kcj^2)$ are functions of mode $\ksj^2$ and $\kcj^2$, we can write,
% the modes $\ksj$ and $\kcj$ into a set of real modes and a set of complex conjugate pairs.
% \begin{align*}
%     \sum_j \vsj \wsj &= \bs_i,\\
%     \sum_j \vcj \wcj &= \bc_i
% \end{align*}

\begin{align*}
    \sum_{\substack{p=1\\\ks_p^2\in\mathcal{C}}}^{\tilde{n}_p} \Ks_{ip}(\ks_p^2) \ws_p + \sum_{\substack{q=\tilde{n}_p+1\\\ks_q^2\in\mathcal{C}}}^{2\tilde{n}_p} \Ks_{iq}(\ks_q^2) \ws_q + \sum_{\substack{r=2\tilde{n}_p+1\\\ks_r^2\in\mathcal{C}}}^{\tilde{n}} \Ks_{ir}(\ks_r^2) \ws_r &= \bs_i,\\
    \sum_{p=1}^{\hat{n}_p} \Kc_{ip}(\kc_p^2) \wc_p + \sum_{q=\hat{n}_p+1}^{2\hat{n}_p} \Kc_{iq}(\kc_q^2) \wc_q + \sum_{\substack{r=2\hat{n}_p+1\\\ks_r^2\in\mathcal{R}}}^{\hat{n}} \Kc_{ir(\kc_r^2)} \wc_r &= \bc_i.
\end{align*}
Since $\ks_q^2=\ks_p^{*2}, \forall (\ks_p^2, \ks_q^2)\in\mathcal{C}$ with $q=p+\tilde{n}_p$, and  $\kc_q^2=\kc_p^{*2}, \forall (\kc_p^2, \kc_q^2)\in\mathcal{C}$ with $q=p+\hat{n}_p$, it follows that $\Ks_{iq}(\ks_q^2) =  \Ks_{ip}^*(\ks_p^{*2})$, and $\Kc_{iq}(\kc_q^2) =  \Kc_{ip}^*(\kc_p^{*2})$. Therefore, we have
\begin{align*}
    \sum_{\substack{p=1\\\ks_p^2\in\mathcal{C}}}^{\tilde{n}_p} \left(\underbrace{\Ks_{ip}(\ks_p^2)}_{\text{complex}} \ws_p + \underbrace{\Ks_{ip}^*(\ks_p^{*2})}_{\text{conjugate}} \ws_{p+\tilde{n}_p}\right) + \sum_{\substack{r=2\tilde{n}_p+1\\\ks_r^2\in\mathcal{R}}}^{\tilde{n}} \underbrace{\Ks_{ir}(\ks_r^2)}_{\text{real}} \ws_r &= \underbrace{\bs_i}_{\text{real}},\\
    \sum_{\substack{p=1\\\kc_p^2\in\mathcal{C}}}^{\hat{n}_p} \left(\underbrace{\Kc_{ip}(\kc_p^2)}_{\text{complex}} \wc_p + \underbrace{\Kc_{ip}^*(\kc_p^{*2})}_{\text{conjugate}} \wc_{p+\hat{n}_p}\right) + \sum_{\substack{r=2\hat{n}_p+1\\\ks_r^2\in\mathcal{R}}}^{\hat{n}} \underbrace{\Kc_{ir(\kc_r^2)}}_{\text{real}} \wc_r &= \underbrace{\bc_i}_{\text{real}}.
\end{align*}
Since the right-hand sides of the above equations are real, the left-hand sides must also be real. Thus, $\Ks_{ir}(\ks_r^2)$ and $\Kc_{ir}(\kc_r^2)$ being real implies that $\ws_r$ and $\wc_r$ are real as well. From Eq. \eqref{eqn:wsj}, it follows that $\usj\in \mathcal{R}$ if $\wsj\in \mathcal{R}$ and $\ksj\in\mathcal{R}$, and  $\usj\in \mathcal{I}$  if $\wsj\in \mathcal{R}$ and $\ksj\in\mathcal{I}$. Similarly, from Eq. \eqref{eqn:wcj}, $\ucj\in \mathcal{R}$ if $\wcj\in \mathcal{R}$, regardless of whether  $\ksj\in\mathcal{R}$ or $\ksj\in\mathcal{I}$.

For the sum of the first two terms to be real, the complex conjugates $\Ks_{ip}(\ks_p^2)$ and $\Ks_{ip}^*(\ks_p^{*2})$ require that $\ws_p$ and $\ws_{p+\tilde{n}_p}$ are also complex conjugates. Consequently, from Eq. \eqref{eqn:wsj}, $\us_p$ and $\us_{p+\tilde{n}_p}$ are also complex conjugates. Similarly, $\wc_p$ and $\wc_{p+\hat{n}_p}$ are conjugates, which means $\uc_p$ and $\uc_{p+\hat{n}_p}$ are conjugates according to Eq. \eqref{eqn:wcj}. Therefore, we have:
\begin{align*}
    \ws_{p+\tilde{n}_p} &= \ws_{p}^* \Longrightarrow \us_{p+\tilde{n}_p} = \us_p^*,\\
    \wc_{p+\hat{n}_p} &= \wc_{p}^* \Longrightarrow \uc_{p+\hat{n}_p} = \uc_p^*.
\end{align*}

Table \ref{tab:usuc_kskc} shows possible combinations of mode numbers $\ksj, \kcj$ and mode amplitudes $\usj, \ucj$. Specifically, $\ucj$ can be either real ($\mathcal{R}$) or a complex conjugate pair ($\mathcal{C}$), while $\usj$ can be real ($\mathcal{R}$), imaginary ($\mathcal{I}$), or complex conjugate pair ($\mathcal{C}$).
%%%%%%%%%%%%%%%%%%%%%%%%%
\begin{table}[]
\centering
%\resizebox{\textwidth}{!}
\begin{minipage}{0.45\textwidth}
\centering
{%
\begin{tabular}{cccc}
\cmidrule(lr){1-4}
\multirow{2}{*}{} & \multicolumn{3}{c}{$\ksj^2$}   \\ \cmidrule(lr){2-4}
  & \multicolumn{1}{c}{$\mathbb{R}_+$} & \multicolumn{1}{c}{$\mathbb{R}_-$}      & $\mathcal{C}$ \\ 
  \cmidrule(lr){1-4}
$\ksj$  & \multicolumn{1}{c}{$\mathcal{R}$} & \multicolumn{1}{c}{$\mathcal{I}$} & $\mathcal{C}$ \\ 
$\wsj$            & \multicolumn{1}{c}{$\mathcal{R}$} & \multicolumn{1}{c}{$\mathcal{R}$}      & $\mathcal{C}$ \\
$\usj$            & \multicolumn{1}{c}{$\mathcal{R}$} & \multicolumn{1}{c}{$\mathcal{I}$} & $\mathcal{C}$  
\\\cmidrule(lr){1-4}
\end{tabular}%
}
\end{minipage}
% \hspace{1in}
\begin{minipage}{0.45\textwidth}
\centering
{%
\begin{tabular}{cccc}
\cmidrule(lr){1-4}
\multirow{2}{*}{} & \multicolumn{3}{c}{$\kcj^2$} \\ \cmidrule(lr){2-4}
   & \multicolumn{1}{c}{$\mathbb{R}_+$} & \multicolumn{1}{c}{$\mathbb{R}_-$}      & $\mathcal{C}$ \\ \cmidrule(lr){1-4}
$\kcj$ & \multicolumn{1}{c}{$\mathcal{R}$} & \multicolumn{1}{c}{$\mathcal{I}$} & $\mathcal{C}$ \\ 
$\wcj$            & \multicolumn{1}{c}{$\mathcal{R}$} & \multicolumn{1}{c}{$\mathcal{R}$}      & $\mathcal{C}$ \\
$\ucj$            & \multicolumn{1}{c}{$\mathcal{R}$} & \multicolumn{1}{c}{$\mathcal{R}$}      & $\mathcal{C}$ 
\\\cmidrule(lr){1-4}
\end{tabular}%
}
\end{minipage}
\caption{Mode amplitudes $\usj$ and $\ucj$ correspond to mode numbers $\ksj$ and $\kcj$ as described in Eqs. \eqref{eqn:ujkj_sin} and \eqref{eqn:ujkj_cos}; $\wsj$ and $\wcj$ are defined in \eqref{eqn:wsj}, and \eqref{eqn:wcj}. The disjoint sets $\mathcal{R}, \mathcal{I}$, and $\mathcal{C}$ consist of real numbers, purely imaginary numbers, and complex conjugate pairs, respectively.}
\label{tab:usuc_kskc}
\end{table}
\end{remark}
%%%%%%%%%%%%%%%%%%%%%%%%%
\subsection{Algorithm}\label{sec:algorithm}
In this section, we present our GFS algorithm designed to compute derivatives of smooth functions in non-periodic domains. The method relies on decomposing the function into periodic and aperiodic components. The periodic part is managed using FFT, while the aperiodic part is represented by a set of adaptively constructed non-harmonic modes, effectively capturing derivative jumps at domain boundaries.

This algorithm dynamically constructs the aperiodic basis to align with the function’s smoothness and discontinuities. The jump conditions required for this basis are either analytically derived or numerically computed. Subtracting the adaptively modeled aperiodic part from the original function isolates the periodic component, which is then differentiated using FFT. The periodic and aperiodic derivatives are recombined to produce the final result.

A key strength of this method lies in its ability to achieve high accuracy without requiring domain extensions, a common limitation of Fourier extension techniques. The algorithm demonstrates computational efficiency, with complexity scaling as $\mathcal{O}(N(n+\log_{2}(N)))$, where $N$ stands for the number of grid points and $n$ stands for the number of non-harmonic modes (see Sec. \ref{sec:comp_complexity}). Numerical experiments highlight its superior resolution and convergence rates compared to finite-difference methods, with significantly lower degrees of freedom required per wavelength.

The steps involved in computing the derivative $u'(x_i), i=1,2,\hdots,N$ are summarized in Algorithm \ref{algo:gfs}.
\begin{algorithm}
\caption{GFS algorithm for decomposing a function into periodic and non-periodic parts and computing derivatives}
\begin{algorithmic}[1]
\State Choose the number of symmetric modes $n_c$ and anti-symmetric modes $n_s$ in Eq. \eqref{def:ua}, typically setting $n_c=n_s=n\sim \mathcal{O}(1)$.
\State Evaluate the jumps in the derivatives of $u(x)$ across the domain boundaries $J_m$, either analytically (if known) or numerically.
% \If{rank$(\Js) < n_s$ in Eq. \eqref{eqn:es_n} or  rank$(\Jc) < n_c$ in Eq. \eqref{eqn:ec_n}}
%     \State Set $n_s = \text{rank}(\Js), n_c=\text{rank}(\Jc)$ accordingly.
%     \State Update $\Js, \Jc$ to match the new $n_s$ and $n_c$. %Repeat step 2 before continuing to step 4.
% \EndIf
\State Compute the elementary symmetric polynomials $\esj$ and $\ecj$ from Eqs. \eqref{eqn:es_n}-\eqref{eqn:ec_n}.
\State Solve for the ($n$) roots of Eq. \eqref{eqn:ks_characteristic} and Eq. \eqref{eqn:kc_characteristic} to obtain $\ksj$ and $\kcj$, for $j=1,2,\hdots,n$.
\State Compute $\usj$ and $\ucj$ using the first equality in Eqs. \eqref{eqn:ujkj_sin}-\eqref{eqn:ujkj_cos}.
\State Calculate $u_a(x_i) = \sum_{j=1}^n \left(\ucj \cos(\kcj x_i) + \usj \sin(\ksj x_i)\right)$, $u'_a(x_i) = \sum_{j=1}^n \left(-\kcj\ucj \sin(\kcj x_i) + \ksj\usj \cos(\ksj x_i)\right)$.
\State Determine $u_{p}(x_i) = u(x_i)-u_a(x_i)$ and obtain $u'_{p}(x)$ using the Fast Fourier Transform (FFT).
\State Finally, compute the derivative $u'(x_i) = u'_{p}(x_i) + u'_a(x_i)$. Higher derivatives, if required, can be computed analogously.
\end{algorithmic}
\label{algo:gfs}
\end{algorithm}
\subsection{Numerical computation of jump conditions}\label{sec:fd}
Now, we briefly describe the classical one-sided finite difference schemes used to compute the derivatives, \(u^{(m)}(x)\), at the domain boundaries, which are then used to evaluate the jumps $J_m$ in $u^{(m)}(x)$ across the domain boundaries. The computation of higher-order derivatives with high accuracy in a stable manner is a significant challenge in numerical analysis and scientific computing. Numerical differentiation, especially for higher-order derivatives, can suffer from stability issues, rounding errors, and amplification of small numerical errors.

Let us express the $d^{th}$ derivative of a function $u(x)$ at a discrete node $j$ ($u_j\equiv u(x_j)$) in terms of the neighbor nodes $u_{j+m}$ and use Taylor series expansion w.r.t. $u_j$ to have,
\begin{equation}
 \Delta x^d~ u^{(d)}_j = \sum_{m=0}^M a_m u_{j+m} = \sum_{m=0}^M a_m \sum_{n=0}^\infty \frac{u_j^{(n)}}{n!} (m\Delta x)^n = \sum_{n=0}^\infty \left(\sum_{m=0}^M a_m  \frac{m^{(n)}}{n!} \right) \Delta x^n u_j^{(n)}.
 \label{eqn:Taylor_gen}
\end{equation}
Here, $M+1=d+r$ stands for the width of the computational stencil used in finite difference approximation. The weights $a_m$ ($M+1$ unknowns) can be computed by solving the system of equations below obtained by comparing equal order terms on each side with a formal order of accuracy of $\Delta x^{M+1-d}$ in Eq. \eqref{eqn:Taylor_gen}, 
\begin{eqnarray}
 \sum_{m=0}^M a_m \frac{m^n}{n!} = \delta_{n,d} = 
 \begin{cases} 
 1, \hspace{0.25in} n=d,\\
 0, \hspace{0.25in} n\neq d,
 \end{cases}
\quad \quad n=0,1,2,3,...,M.
\end{eqnarray}

Above expressions involving a forward differencing stencil can be used on the left-side boundary, and for the right-side boundary, a one-sided backward differencing stencil must be used, and the corresponding system of equations can be obtained by replacing $m$ with $-m$ in the above equations.

\section{Computational cost}\label{sec:comp_complexity}
The computational complexity of each step of the algorithm, described in Sec. \ref{sec:algorithm}, is shown in Table \ref{tab:complexity}.
Numerical evaluation of $d^{th}$ derivative via finite difference approximation with $\mathcal{O}(\Delta x)^r$ accuracy requires $\mathcal{O}(r+d)$ operations where $(r+d)$ is the stencil width. The total cost for evaluating $J_m, m=0,1,2,\hdots,4n-1$ is $\mathcal{O}(\sum_{d=1}^{4n-1} (r+d)) = \mathcal{O}((4n-1) (r+2n)) \approx \mathcal{O}(4nr+8n^2)$. Evaluation of the elementary symmetric polynomials involve inversion of a $n\times n$ matrix which requires $\mathcal{O}(\frac{2}{3}n^3+\frac{1}{2}n^2+\frac{n}{6})$ operations through LU decomposition and a matrix vector product of $\mathcal{O}(2n^2+n)$ operations. Finding roots of a $n^{th}$ degree polynomial requires another $\mathcal{O}(n^3)$ operations.
Computation of the expansion coefficients $\us_j, \uc_j$ involve inversion of a $n\times n$ matrix and require $\mathcal{O}(n^3+n)$ operations. \rev{Evaluation of $u_a(x_j), u'_a(x_j), j=1,2,\hdots,N$ requires $\mathcal{O}(2Nn)$ operations each}. Evaluation of $u'_{p}(x_j)$ through FFT requires $\mathcal{O}(N\log_2 N)$ operations. Finally, $u'(x_j)=u'_{p}(x_j) + u'_a(x_j)$ requires $\mathcal{O}(N)$ operations. Overall, the total cost of the computational scales as $\mathcal{O}(n^3+9n^2+N(\rev{4}n+\log_2 N))$ operations. 

We demonstrate later that $n=\mathcal{O}(1)\ll N$ is sufficient for practical purposes so that the overall cost remains tractable with an effective complexity of $\mathcal{O}(N(\rev{4}n+\log_2 N))$. We can conclude that for $n=\mathcal{O}(1)\ll N$, the computational penalties due to the evaluation of the jumps $J_m$, the modes $\ks, \kc$ and the expansion coefficients $\us, \uc$ are negligible and the major contribution comes from the evaluation of $u_a(x_j), u'_a(x_j)$ at each node ($x_j$), similar to the reduction-to-periodicity technique of \cite{Roache:1978}.
\begin{table}[tbp]
\centering%\scriptsize
\begin{tabular}{|c|c|}
\hline
 step  & cost \\
\hline 
 Numerical evaluation of $J_m$ & $\mathcal{O}(4nr+8n^2)$\\
 Computation of $\es_j, \ec_j$ & $\mathcal{O}(n^3+n^2)$\\
 Polynomial roots $\ks_j, \kc_j$ & $\mathcal{O}(n^3)$ \\
 Computation of $\us_j, \uc_j$ & $\mathcal{O}(n^3+n)$\\
 Computation of $u_a(x_j), u'_a(x_j)$ & $\mathcal{O}(\rev{4}Nn)$\\
 Computation of $u'_{p}(x_j)$ via FFT & $\mathcal{O}(N\log_2 N)$\\
 Computation of $u'(x_j) = u'_p(x_j)+u'_a(x_j)$& $\mathcal{O}(N)$\\
 \hline
 % Total cost & $\mathcal{O}(3n^3+9n^2+n(4r+1)+N(n+1)+N \log_2 N)$\\
 Total cost & $\mathcal{O}(3n^3+9n^2+4nr+\rev{4}Nn+N \log_2 N)$\\
 Effective cost ($n=\mathcal{O}(1)\ll N, r=\mathcal{O}(1)$) & $\mathcal{O}(N (\rev{4}n+\log_2 N) \rev{+ 3 n^3})$\\
 \hline
\end{tabular}
\caption{Estimate of computational cost for \rev{evaluating derivative of $u(x)$} with GFS. Here, $M+1=r+d$ is the stencil width used in finite difference (FD) approximation of the jumps $J_m$, $r$ is the accuracy of a FD scheme for $d^{th}$ derivative. For the highest derivative $d=4n-1$.}
\label{tab:complexity}
\end{table}

%\cite{Roache:1978}: The reduction-to-periodicity technique adds to the operation count of the usual pseudospectral FFT method. For $N = 1/\Delta x >1$, the penalties due to the boundary evaluations of D, and the evaluations of the coefficients a, are negligible, and the major contribution comes from the evaluation of the reducing polynomial g and its derivatives at each node point. 
% ===================

%======================================================================
\section{Numerical Evidences}\label{sec:Numerical_Evidences}
%======================================================================
In this section, we compare the accuracy and convergence of our GFS method for computing the first derivative with \rev{Eckhoff \cite{Eckhoff:1995}, Roache \cite{Roache:1978}, Prony \cite{Prony:1795},} finite difference (FD) and \rev{standard} FFT methods. \rev{For the Eckhoff method, discontinuities are assumed at the boundary and analytical derivativees are used for computing all the jumps ($J_m, m=0,1,\hdots,q-1$) for both Eckhoff and Roache methods.} For the FD method, standard central schemes are applied at interior nodes, while one-sided stencils are used near boundary nodes with a computational stencil width $r+1$, leading to $\mathcal{O}(\Delta x^r)$ accuracy, where $\Delta x$ is the grid spacing (see Sec. \ref{sec:fd} for details). In the GFS method, the jumps $J_m$ for $m=0,1,\hdots,4n-1$ are evaluated both analytically and numerically using one-side FD schemes with a stencil width $4n-1+r$, giving the highest derivative's jump at the end points, $J_{4n-1}$, a formal accuracy of $\mathcal{O}(\Delta x^r)$. However, numerical evaluation of jumps introduces additional error to the GFS method.  
Computing higher derivatives of highly oscillatory functions is challenging due to the inherent ill-posedness known as Runge's phenomenon, with increasing rounding errors reducing overall accuracy as the magnitude of the higher derivatives increases.
%Computing higher derivatives of highly oscillatory functions is challenging due to the inherent ill-posedness known as Runge's phenomenon, and round-off errors become significant as the higher derivative's magnitude increases, limiting  overall accuracy. 
We present the numerical error of the GFS method, comparing results with jumps $J_m$ evaluated both analytically and numerically to highlight the impact of numerical evaluation.
The error ($e$) for each method is computed against the analytic values and quantified using  the $L^p$-norm of $e$, defined as follows:
\begin{align*}
 ||e||_p = \left(\Delta x \sum_{i=0}^N |e(x_i)|^p  \right)^{1/p}   
\end{align*}

% \subsection{Test function: $u(x)=e^{a(x+\pi)}\sin[b(x+\pi)]$}
\subsection{Modulated sine function}
\begin{figure}[htp]
 \centering
 \begin{minipage}{0.45\linewidth}\centering
 \includegraphics[width=2.25in,trim={0.in 0.in 0.in 0in},clip=true]{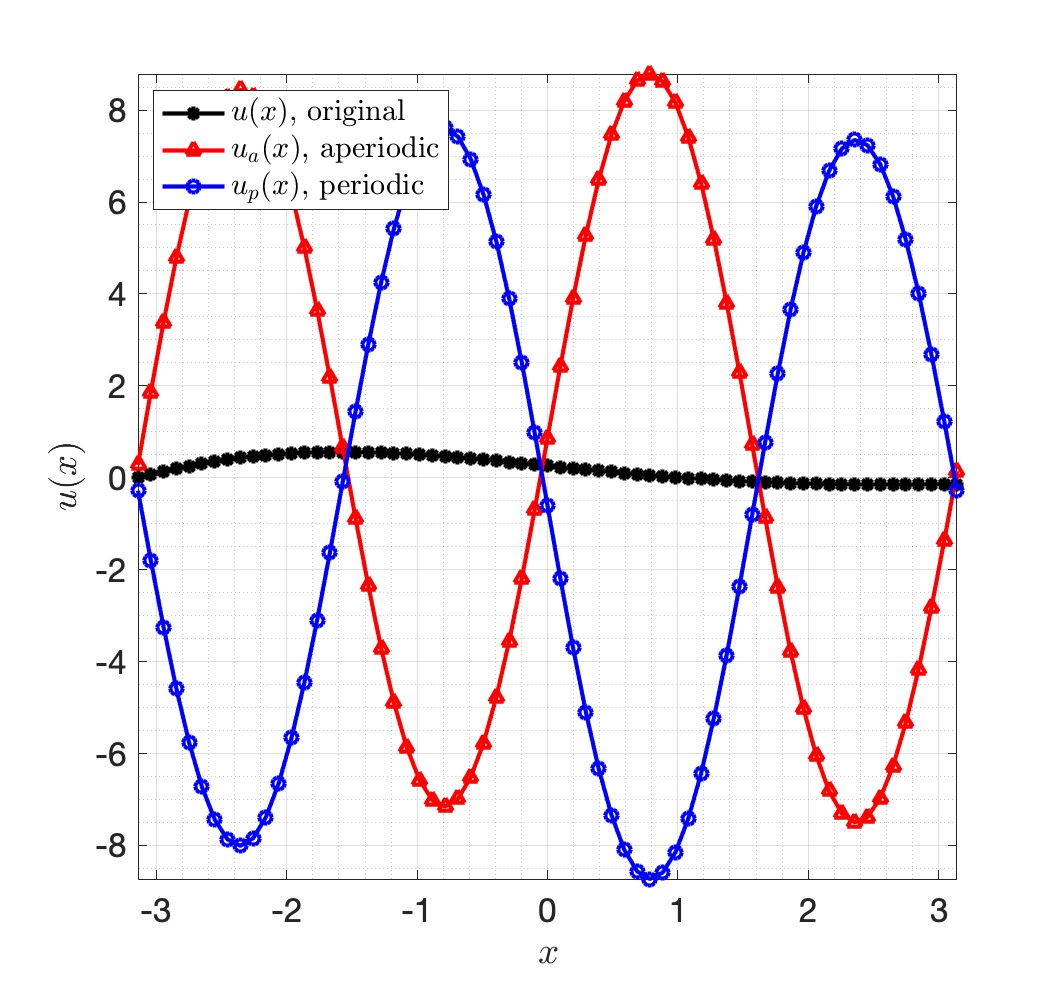} 
 \\(a) \end{minipage}~~
\begin{minipage}{0.5\linewidth}\centering
\includegraphics[width=2.25in,trim={0.in 0.in 0.in 0in},clip=true]{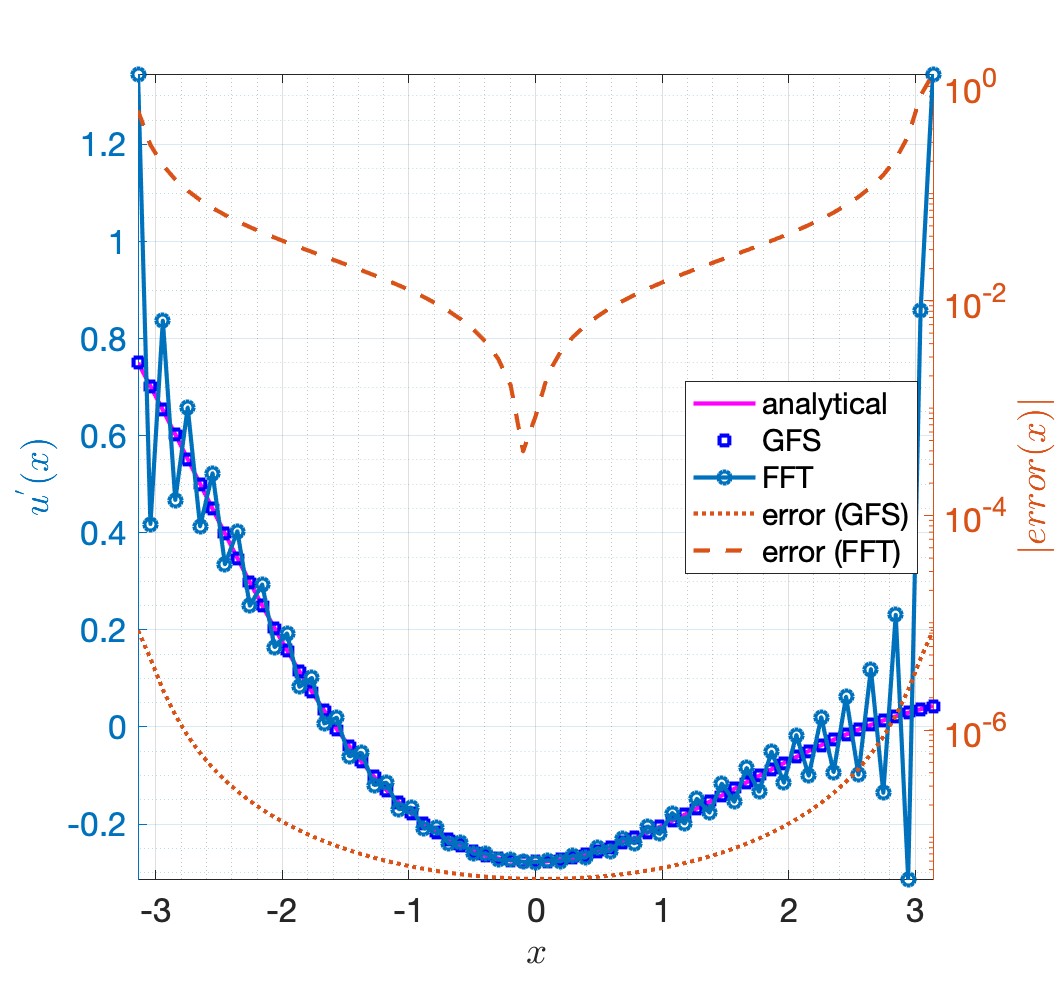}  
\\(b)
 \end{minipage}\\
 \begin{minipage}{0.45\linewidth}\centering
 \includegraphics[width=2.25in,trim={0.in 0.in 0.in 0in},clip=true]{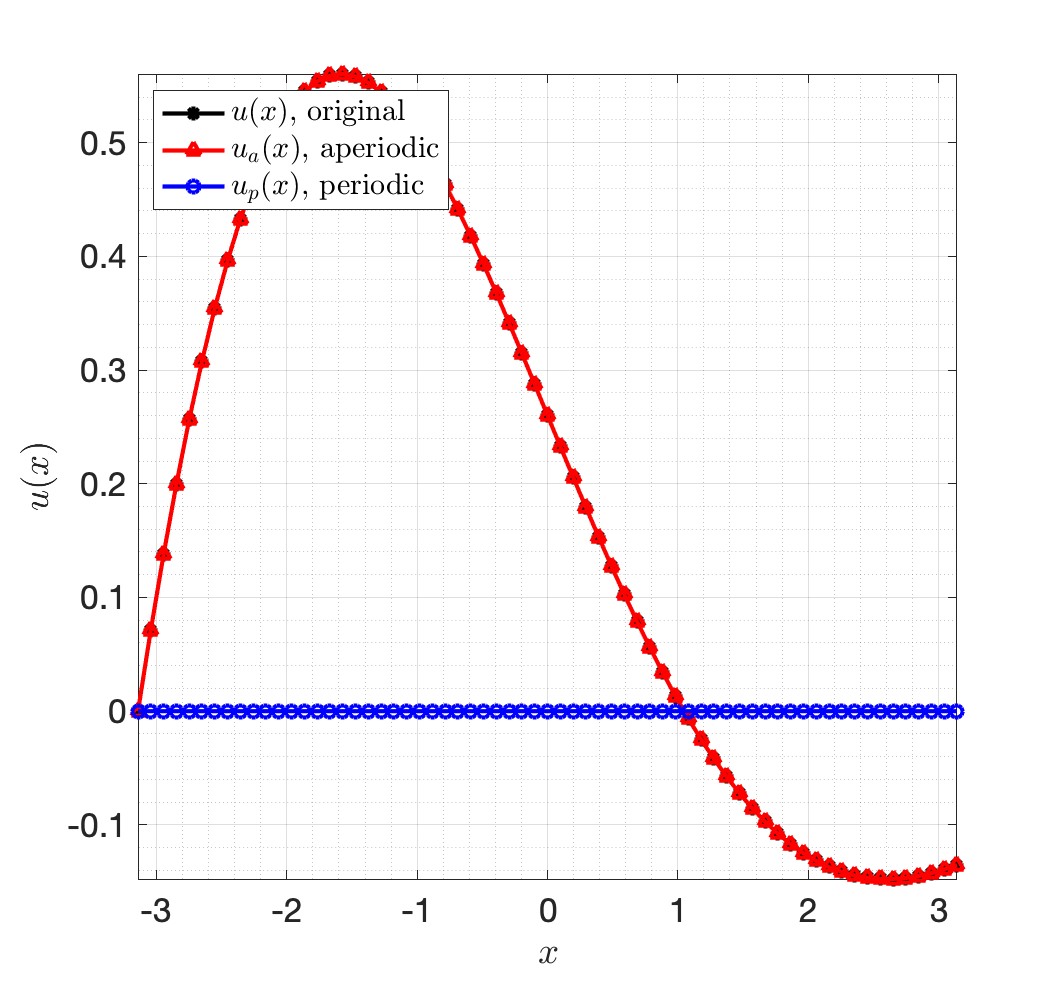} 
 \\(c) \end{minipage}~~
\begin{minipage}{0.5\linewidth}\centering
\includegraphics[width=2.25in,trim={0.in 0.in 0.in 0in},clip=true]{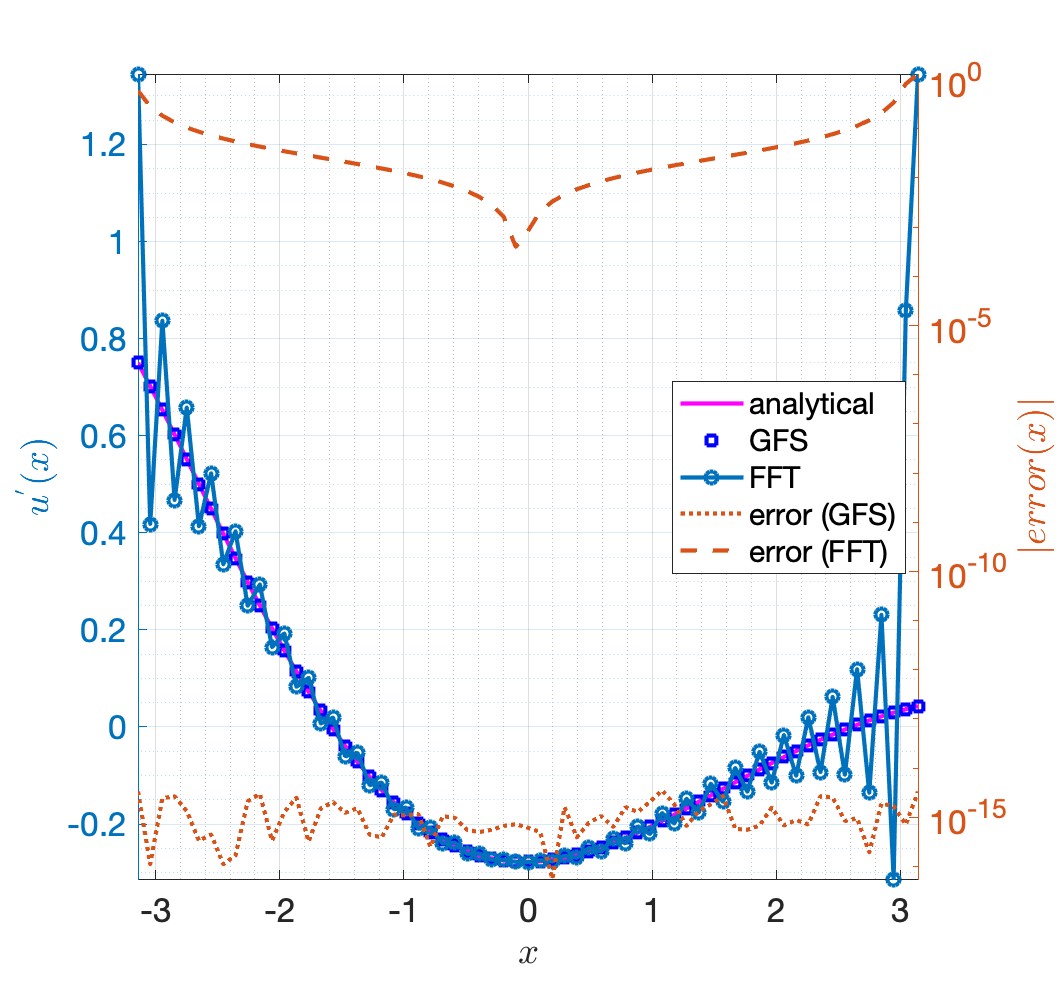}  
\\(d)
 \end{minipage}
 \caption{A modulated sine function, defined in Eq. \eqref{eqn:modulated_sine}, is approximated using GFS with $N=64$ mesh points. The decomposition into periodic and aperiodic parts is shown in  (a, c), while (b,d) compares the numerical error in the first derivative  for GFS and FFT  against analytical values for $n=1$ and 2 aperiodic modes, respectively. The FFT shows Gibbs oscillations, resulting in $\mathcal{O}(1)$ error, whereas GFS 
avoids these oscillations, yielding a significantly lower numerical error by several orders of magnitude.}
 \label{fig:eax}
\end{figure}

Consider the modulated sine function 
\begin{align}
    u(x)=e^{a(x+\pi)}\sin[b(x+\pi)], x\in [-\pi,\pi],\label{eqn:modulated_sine}
\end{align}
where, $a=-1/\pi, b=3/4$. 
%This is equivalent to the static test function $u(x)=e^{ax}\sin(2\pi b x), x\in [0,1]$ with $a=-2, b=3/4$ considered by \cite{Roache:1978}. 
%Here, the function $u(x)$ has two degrees of freedom (exponents $a,b$) in the space spanned by complex exponential basis. 
Figs. \ref{fig:eax}(a,c) illustrate the function $u(x)$ and its decomposition into periodic ($u_{p}(x))$ and aperiodic ($u_a(x)$) components for $n=1$ and $2$, respectively. 
In Fig. \ref{fig:eax}(a), for $n=1$, the magnitude of both parts exceed that of the original function, and their mutual cancellation reproduces the original function, indicating a non-normal basis. Conversely, Fig. \ref{fig:eax}(c) demonstrates that $u_a(x)$ completely resolves $u(x)$ for $n=2$, resulting in $u_{p}(x)=0$.

Figs. \ref{fig:eax}(b,d) compare the GFS and FFT methods to the analytical values, demonstrating that GFS avoids Gibbs oscillations, unlike FFT. Table \ref{tab:eax} quantifies the $L^2, L^\infty$ norms of the approximation error for the function ($e$) and its derivative ($e^\prime$). As the number of aperiodic modes $n$ in GFS increases, the error in the first derivative $e^\prime$  converges rapidly, reaching machine precision by $n=2$ with analytically evaluated jump conditions $J_m$. This faster convergence is attributed to the two degrees of freedom (exponents $a,b$) of the function $u(x)$ in the space of complex exponential functions. Even when $J_m$ is evaluated numerically, $e^\prime$ decays quickly, with accuracy reliant on the finite difference (FD) schemes. Enhanced accuracy in the FD scheme for evaluating $J_m$ decreases numerical error, enabling the overall error to approach that of the analytically evaluated case.
%However, when $J_m$ are evaluated numerically using one-sided finite difference schemes of $\mathcal{O}(\Delta x^r)$ accuracy (FD in Table \ref{tab:eax}), the rate of convergence with $n$ decreases and the overall accuracy is reduced. 

\begin{table}[htp]
\centering
\begin{tabular}{|c|cc|cc|c| }
\hline
%\diagbox{$n$}{$k_n$} 
$n$ & $\max |e|$ & $||e||_2$ & $\max |e^\prime|$ & $||e^\prime||_2$ & Jump ($J_m$)\\
\hline
  1 & $8.02\times 10^{-15}$ & $6.19\times 10^{-15}$ & $8.49\times 10^{-6}$ & $4.49\times 10^{-6}$ & analytical\\
  1 & $7.23\times 10^{-15}$ & $5.55\times 10^{-15}$ & $1.19\times 10^{-5}$ & $5.14\times 10^{-6}$ & FD2\\
  1 & $7.74\times 10^{-15}$ & $5.30\times 10^{-15}$ & $8.49\times 10^{-6}$ & $4.49\times 10^{-6}$ & FD4\\
  1 & $1.06\times 10^{-14}$ & $5.29\times 10^{-15}$ & $8.49\times 10^{-6}$ & $4.49\times 10^{-6}$ & FD6\\
  1 & $4.69\times 10^{-15}$ & $4.43\times 10^{-15}$ & $8.49\times 10^{-6}$ & $4.49\times 10^{-6}$ & FD8\\
  \hline
  2 & $5.55\times 10^{-17}$ & $1.74\times 10^{-17}$ & $3.52\times 10^{-15}$ & $3.80\times 10^{-15}$ & analytical\\
  2 & $1.11\times 10^{-16}$ & $3.48\times 10^{-17}$ & $1.00\times 10^{-10}$ & $3.17\times 10^{-11}$ & FD2\\
  2 & $5.55\times 10^{-17}$ & $1.74\times 10^{-17}$ & $4.80\times 10^{-14}$ & $1.85\times 10^{-14}$ & FD4\\
  2 & $1.11\times 10^{-16}$ & $3.48\times 10^{-17}$ & $2.22\times 10^{-13}$ & $8.27\times 10^{-14}$ & FD6\\
  2 & $1.36\times 10^{-20}$ & $4.25\times 10^{-21}$ & $6.73\times 10^{-13}$ & $2.97\times 10^{-13}$ & FD8\\
  \hline
\end{tabular}
\caption{A modulated sine function defined in Eq. \eqref{eqn:modulated_sine} is approximated using the GFS for $N=64$. The errors for the function ($e$) and its first derivative ($e^\prime$) are analyzed as the number of sine or cosine modes ($n$) increases. Jumps $J_m$ for $m=0,1,\hdots,4n-1$ are calculated analytically or numerically using one-sided finite difference schemes with stencil width $4n-1+r$ and $\mathcal{O}(\Delta x^r)$ accuracy for $J_{4n-1}$, noted as ``FDr'' in the $J_m$ column.}
\label{tab:eax}
\end{table}

\subsection{Gaussian}
\begin{figure}[htp]
 \centering
 \begin{minipage}{0.45\linewidth}\centering
 \includegraphics[width=2.25in,trim={0.in 0.in 0.in 0in},clip=true]{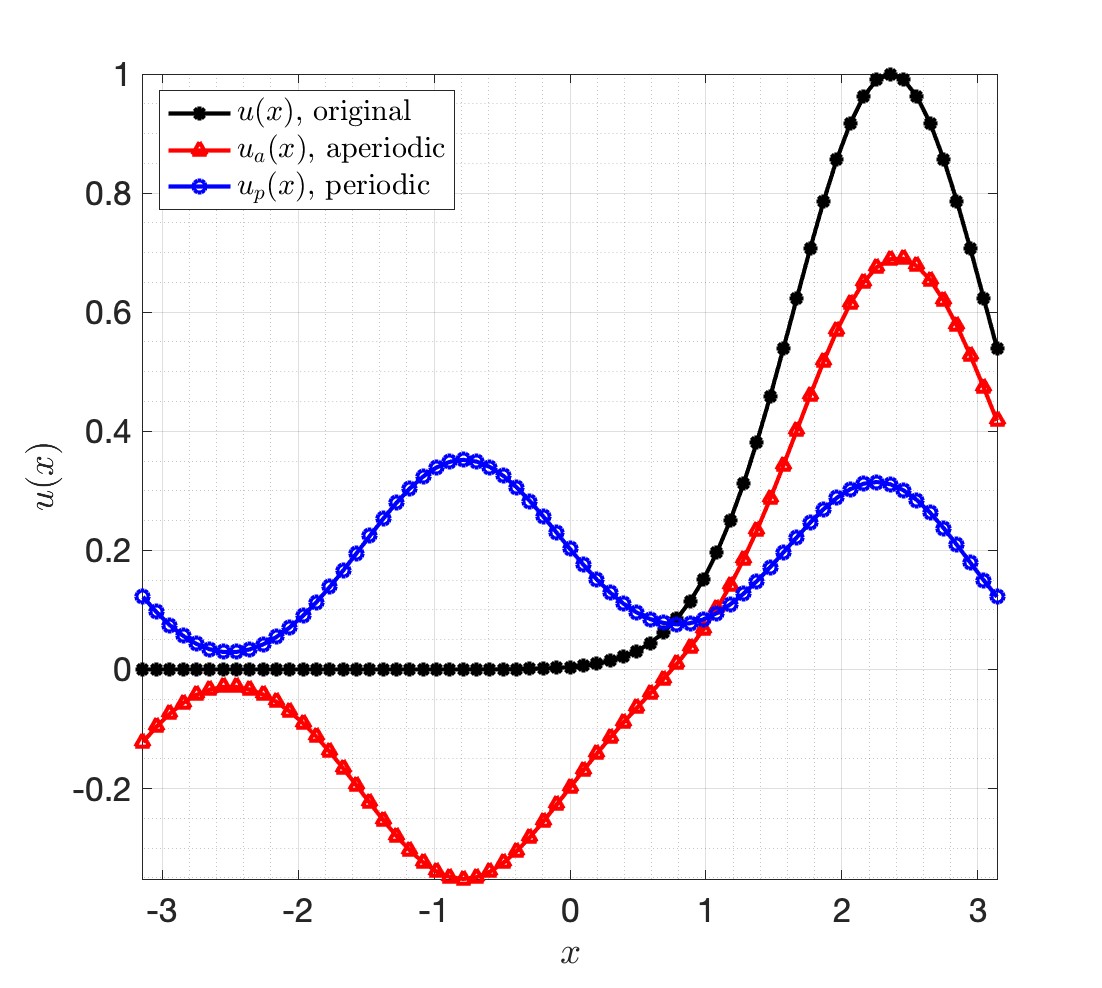} 
 \\(a) \end{minipage}~~
\begin{minipage}{0.5\linewidth}\centering
 \includegraphics[width=2.25in,trim={0.in 0.in 0.in 0in},clip=true]{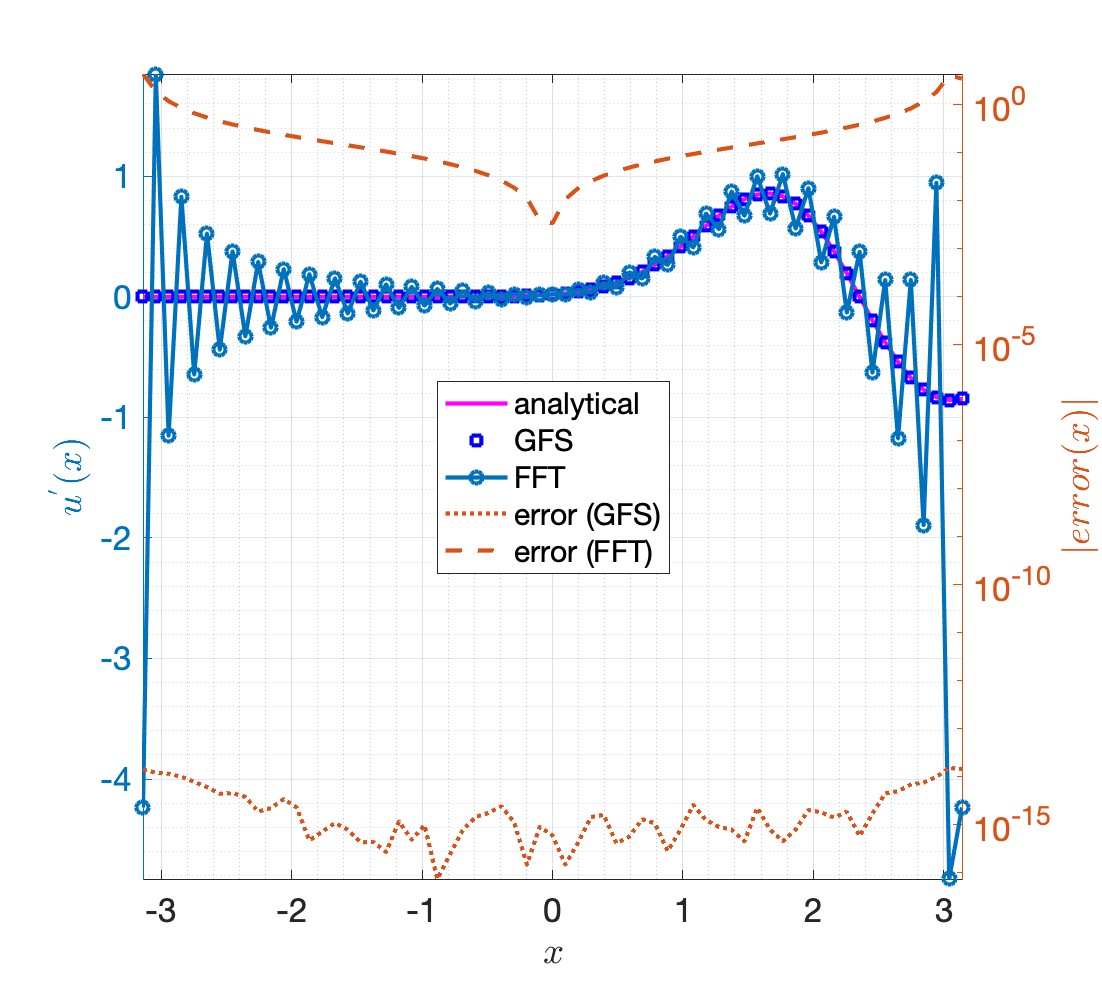} 
\\(b)
 \end{minipage}
 \caption{A Gaussian, defined in Eq. \eqref{def:gauss_func}, is approximated using GFS with $N=64$ grid points, and $n=3$: (a) shows its decomposition into periodic and aperiodic components, while (b) compares the first derivative approximation using GFS and FFT against analytical values. GFS aligns well with the analytical results without exhibiting Gibbs oscillations, unlike FFT, which shows these oscillations.}
 \label{fig:gauss}
\end{figure}

\begin{figure}[htp]
 \centering
 % \begin{minipage}{\linewidth}\centering
 % \includegraphics[width=3.5in,trim={0.in 0.in 0.in 0in},clip=true]{./Figures/convergence_gauss_2norm.jpg}
 \includegraphics[height=2.25in]{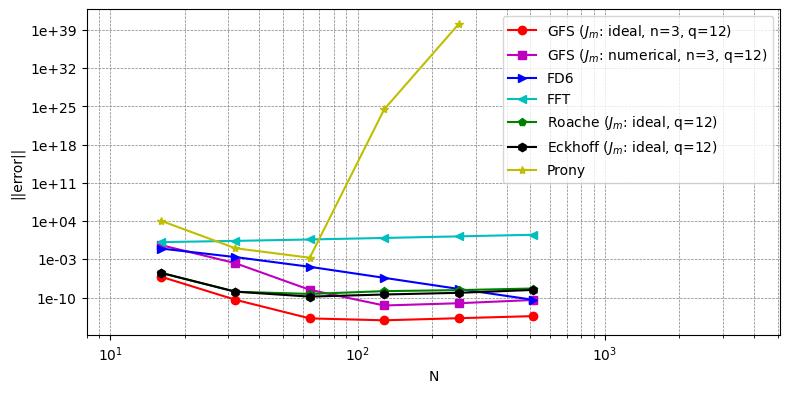}
 % \end{minipage}
  \caption{\rev{For the Gaussian function defined in Eq. \eqref{def:gauss_func}, convergence of the numerical error ($L^\infty-$ norm) in first derivative is compared for GFS, a sixth-order finite difference (FD6), FFT, Roache, Eckhoff, and Prony methods. For Roache and Eckhoff methods the jumps $J_m$ are evaluated analytically, and for GFS method $J_m$ are evaluated both analytically and numerically. Here, ``$J_m:$ ideal'' refers to the case where the jumps $J_m, m=0,1,\hdots,q-1, q=4n,$ are evaluated analytically. GFS is robust and exhibits superior convergence with mesh size $N$ for the same number of jumps ($q$), whether $J_m$ are evaluated analytically or numerically. The Prony method ($M=N/2$) is numerically unstable for large mesh sizes ($N>64$).}}
  \label{fig:convergence_gauss}
\end{figure}

Consider the non-periodic Gaussian function defined on the domain $[-\pi,\pi]$,
\begin{equation}
u(x)=\exp(-[(x-x_0)/w]^2), \quad  -\pi\le x \le \pi.\label{def:gauss_func}
\end{equation}
where, $x_0=3\pi/4, w=1$. Fig. \ref{fig:gauss}(a) displays $u(x)$, along with the periodic part ($u_{p}(x)$) and the aperiodic part ($u_a(x)$) derived from GFS for $n=3$ and $N=64$. Fig. \ref{fig:gauss}(b) compares $u'(x)$ from  GFS, FFT, analytical values, showing that GFS avoids Gibbs oscillations and closely matches the analytic solution, while FFT exhibits Gibbs oscillations. The $L^\infty, L^2$ norms of the error in $u'(x)$ (denoted as $||e^\prime||_p$) are presented in Table \ref{tab:gauss} for $n=3$, with increasing grid resolution $N$, and are compared to a sixth-order finite difference scheme ($r=6$) and FFT. Among the methods, FFT performs worst with $\mathcal{O}(1)$ error due to the Gibbs oscillations. \rev{Fig. \ref{fig:convergence_gauss} illustrates the superior convergence of GFS in comparison to FD, FFT, Roache, Eckhoff, and Prony methods. With analytically evaluated $J_m$, GFS converges rapidly with increasing $N$, in contrast to all the other methods, and also yields significantly lower error. Although numerical evaluation of $J_m$ introduces additional error in GFS, its convergence remains far superior to that of all the other methods. The Prony approximation computed with $M=N/2$ in Eq. \eqref{eqn:prony} becomes ill-conditioned for large $N$ and the computation blows up for $N>64$.}
%This shows that the numerical evaluation of $J_m$ limits the overall accuracy and convergence of the GFS. 

\begin{table}[htp]
\centering
\resizebox{\linewidth}{!}{
\begin{tabular}{|c||c|cc|cc||cc||cc| }
\hline
 Method  &  \multicolumn{5}{c||}{GFS}	& \multicolumn{2}{c||}{FD: $\mathcal{O}(\Delta x^r)$}& \multicolumn{2}{c|}{FFT}\\
 \hline
 Remark & &  \multicolumn{2}{c|}{$J_m:$ analytical} & \multicolumn{2}{c||}{$J_m:$ numerical} & \multicolumn{2}{c||}{$r=6$} & \multirow{2}{*}{} & \multirow{2}{*}{} \\
\hline
$N$ &$n$ & $||e^\prime||_\infty$ & $||e^\prime||_2$ & $||e^\prime||_\infty$ & $||e^\prime||_2$ & $||e^\prime||_\infty$ & $||e^\prime||_2$ &$||e^\prime||_\infty$ & $||e^\prime||_2$ \\
\hline
16 & 3 & 6.48e-07    &	 1.36e-06 & 3.98e-01    &	 2.56e-01 & 9.37e-02    &	 6.04e-02 & 1.39e+00    &	 1.31e+00 \\
32 & 3 & 3.86e-11    &	 4.76e-11 & 1.92e-04 & 8.52e-05 & 2.52e-03 &  1.14e-03 & 2.33e+00 &  1.74e+00 \\
64 & 3 & 1.50e-14    &	 1.20e-14 & 2.55e-09 & 7.99e-10 & 4.18e-05 &  1.35e-05 & 4.24e+00 &  2.40e+00 \\
128 & 3 & 6.66e-15    &	 5.48e-15 & 3.45e-12 & 1.08e-12 & 4.00e-07 &  9.78e-08 & 8.04e+00 &  3.36e+00 \\
256 & 3 & 1.62e-14    &	 1.07e-14 & 8.91e-12 & 1.43e-12 & 3.73e-09 &  8.45e-10 & 1.57e+01 &  4.73e+00 \\
512 & 3 & 3.90e-14    &	 2.15e-14 & 3.52e-11    &	 3.93e-12    &	 3.82e-11    &	 1.03e-11 & 3.09e+01 &  6.67e+00 \\
% 1024 & 3 & 8.18e-14    &	 4.83e-14 & 2.07e-08 & 2.73e-09 & 1.10e-12 &  1.73e-13 & 6.14e+01 &  9.42e+00 \\
\hline
\end{tabular}
}
\caption{A Gaussian defined in Eq. \eqref{def:gauss_func} is approximated using GFS for various grid sizes $N$. The first derivative approximation error ($e^\prime$) is compared among GFS, finited difference (FD), and FFT methods. Jumps $J_m, m=0,1,\hdots,4n-1,$ are calculated analytically or numerically using one-sided FD schemes with stencil width $M+1=r+4n-1$ such that $J_{4n-1}$ has $\mathcal{O}(\Delta x^6)$ accuracy (see Sec. \ref{sec:fd}).}
\label{tab:gauss}
\end{table}
%%%%%%%%%%%%%%%%%%%%%%%%%
%%%%%%%%%%%%%%%%%%%%%%%%%

\subsection{Log function}
\begin{figure}[htp]
 \centering
 \begin{minipage}{0.45\linewidth}\centering
 \includegraphics[width=2.5in,trim={0.in 0.in 0.in 0in},clip=true]{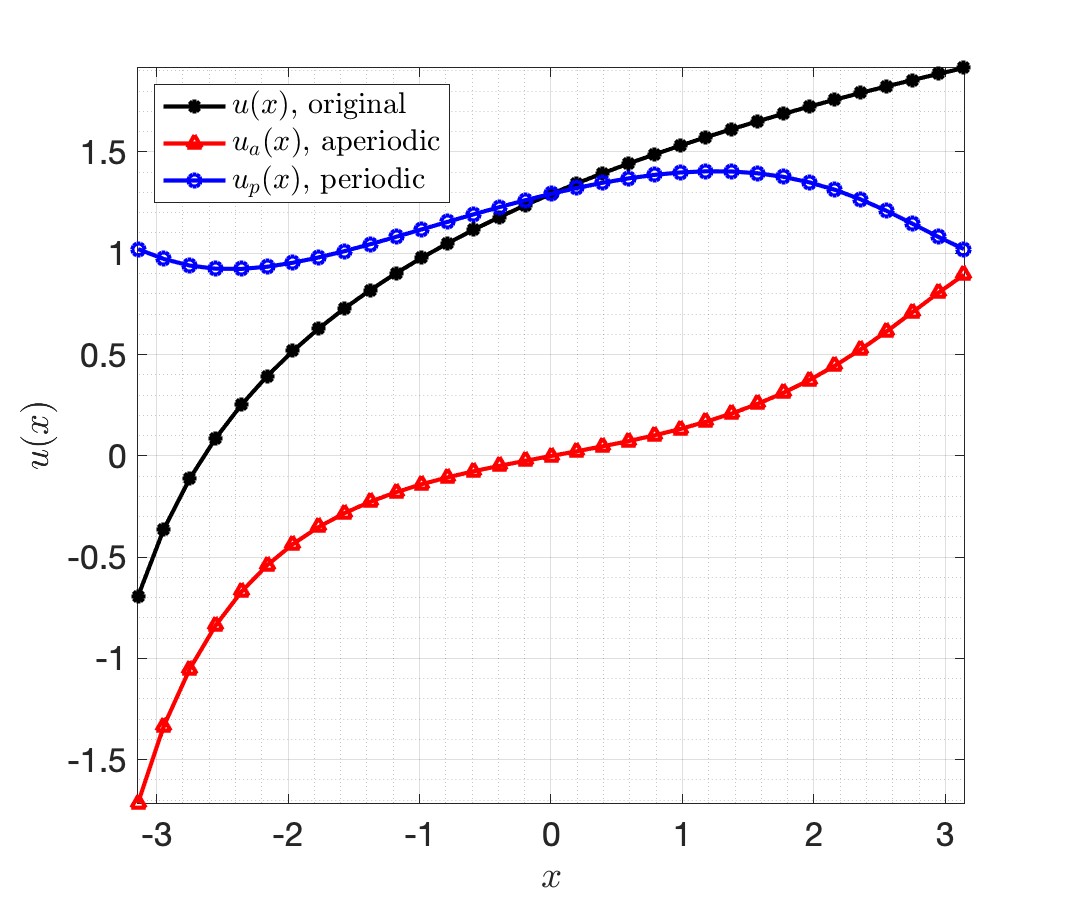} 
 \\(a) \end{minipage}~~
\begin{minipage}{0.5\linewidth}\centering
\includegraphics[width=2.5in,trim={0.in 0.in 0.in 0in},clip=true]{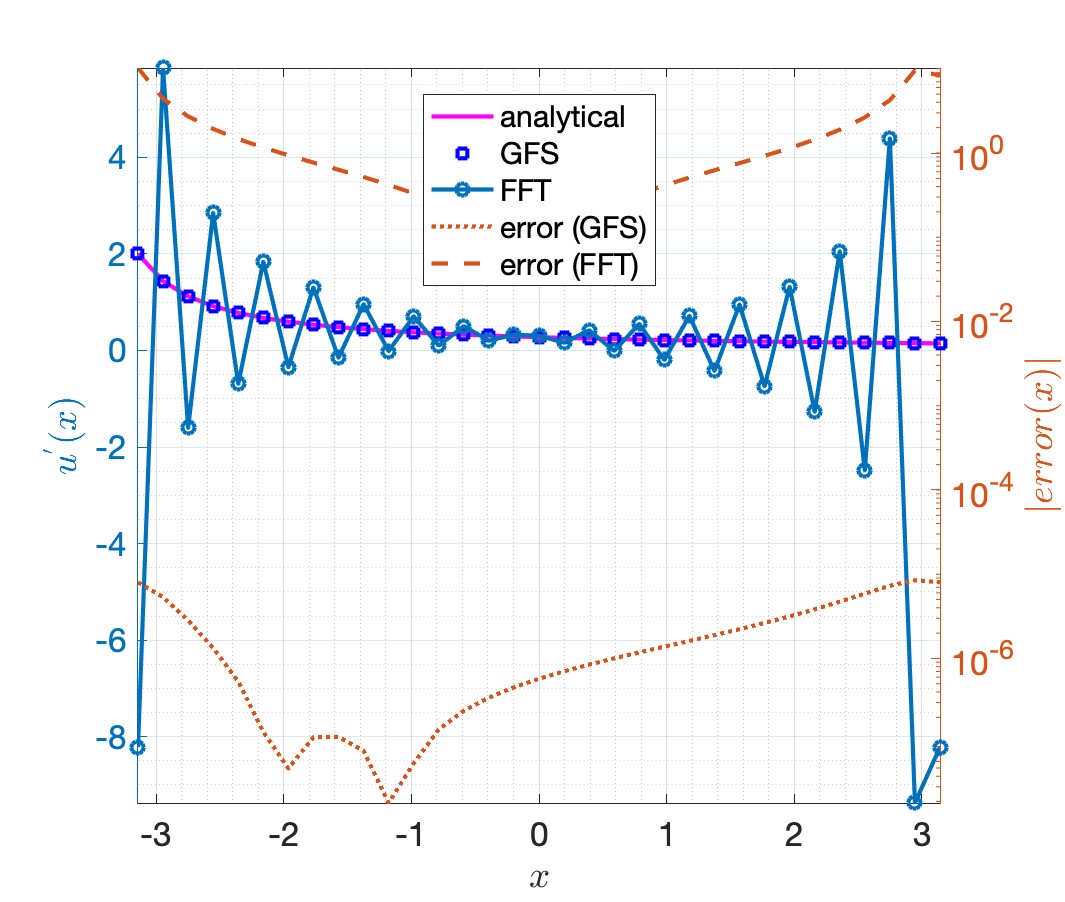}  
\\(b)
 \end{minipage}
 \caption{A function $u(x)=\log(x+\pi+1/2), -\pi\le x \le \pi$, is approximated with GFS using $N=32$ mesh points and $n=3$: Decomposition of $u(x)$ into periodic and aperiodic parts is shown in (a), approximation of the first derivative with GFS and FFT is compared with analytical values in (b); GFS does not exhibit Gibbs oscillations and aligns well with the analytic values while FFT shows Gibbs oscillations.}
 \label{fig:log}
\end{figure}

\begin{figure}[ht]
 \centering
 % \begin{minipage}{\linewidth}\centering
 % \includegraphics[width=3.5in,trim={0.in 0.in 0.in 0in},clip=true]{./Figures/convergence_log_2norm.jpg}
 \includegraphics[height=2.5in]{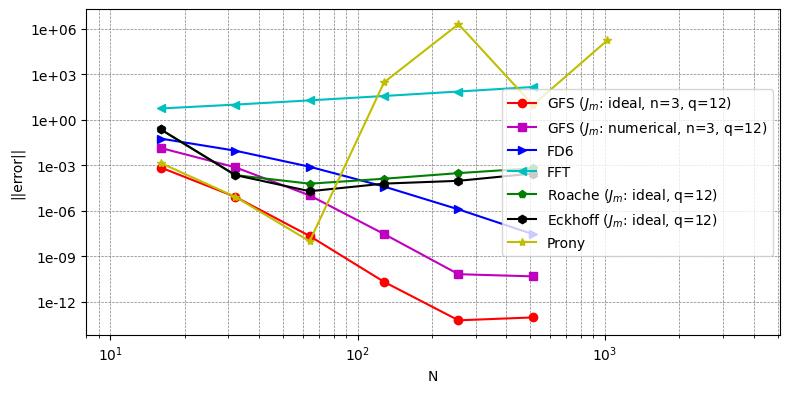}
 % \end{minipage}
  \caption{\rev{For the log function defined in Eq. \eqref{def:log_func}, convergence of the numerical error ($L^\infty-$ norm) in first derivative is compared for GFS, a sixth-order finite difference (FD6), FFT, Roache's, Eckhoff's, and Prony's methods. For the Roache and Eckhoff methods, the jumps $J_m$ are evaluated analytically, and for the GFS method, $J_m$ are evaluated both analytically and numerically. Here, ``$J_m:$ ideal'' refers to the case where the jumps $J_m, m=0,1,\hdots,q-1, q=4n,$ are evaluated analytically.  GFS is robust and shows superior convergence with mesh size $N$ for the same number of jumps ($q$), whether $J_m$ are evaluated analytically or numerically. The Prony method ($M=N/2$) is numerically unstable for large mesh sizes ($N>64$).}
  }
  \label{fig:convergence_log}
\end{figure}

Consider the non-periodic logarithmic function in the domain $[-\pi,\pi]$,
\begin{equation}
u(x)=\log\left(x+\pi+\frac{1}{2}\right), \quad  -\pi\le x \le \pi.\label{def:log_func}
\end{equation}
Fig. \ref{fig:log}(a) displays the function $u(x)$ along with its periodic part $u_{p}(x)$ and aperiodic part $u_a(x)$ obtained using GFS for $n=3$ and $N=32$. Fig. \ref{fig:log}(b) compares $u'(x)$ calculated via  GFS, FFT, and analytical methods. GFS effectively avoids Gibbs oscillations, closely aligning with the analytical method, while FFT exhibits these oscillations. The $L^\infty, L^2$ norms of the error in $u'(x)$ (denoted as $||e^\prime||_p$) are detailed in Table \ref{tab:log} for $n=3$, highlighting the impact of increased grid resolution $N$ and comparisons with a sixth-order finite difference scheme  and FFT. Among the three methods, FFT performs the worst with $\mathcal{O}(1)$ error due to Gibbs oscillations. \rev{Fig. \ref{fig:convergence_log} shows GFS method converging faster and with lower error than FD, FFT, Roache, Eckhoff, and Prony methods. When $J_m$ is evaluated analytically, GFS rapidly improves with increasing $N$, outpacing all the other methods. Even with numerical $J_m$, GFS remains markedly superior. The Prony method with $M=N/2$ in Eq. \eqref{eqn:prony} becomes ill-conditioned for large $N$ and fails for $N>64$.}
 
\begin{table}[htp]
\centering
\resizebox{\linewidth}{!}{
\begin{tabular}{|c||c|cc|cc||cc||cc| }
\hline
 \multirow{2}{*}{}  &  \multicolumn{5}{c||}{GFS}	& \multicolumn{2}{c||}{FD: $\mathcal{O}(\Delta x^r)$}& \multicolumn{2}{c|}{FFT}\\
 \hline
 & &  \multicolumn{2}{c|}{$J_m:$ analytical} & \multicolumn{2}{c||}{$J_m:$ numerical} & \multicolumn{2}{c||}{$r=6$} & \multirow{2}{*}{} & \multirow{2}{*}{} \\
\hline
$N$ &$n$ & $||e^\prime||_\infty$ & $||e^\prime||_2$ & $||e^\prime||_\infty$ & $||e^\prime||_2$ & $||e^\prime||_\infty$ & $||e^\prime||_2$ &$||e^\prime||_\infty$ & $||e^\prime||_2$ \\
\hline
16 & 3 & 7.22e-04    &  9.53e-04 & 1.40e-02    &  8.82e-03 & 5.84e-02    &  3.68e-02 & 5.68e+00    &  5.85e+00 \\
32 & 3 & 8.49e-06    &  8.84e-06 & 7.53e-04    &  3.34e-04 & 9.36e-03    &  4.18e-03 & 1.02e+01    &  8.12e+00 \\
64 & 3 & 2.21e-08    &  1.77e-08 & 1.05e-05    &  3.29e-06 & 8.17e-04    &  2.59e-04 & 1.94e+01    &  1.14e+01 \\
128 & 3 & 2.09e-11    &  1.21e-11 & 2.98e-08    &  6.62e-09 & 4.03e-05    &  9.06e-06 & 3.78e+01    &  1.61e+01 \\
256 & 3 & 6.09e-14    &  5.73e-14 & 6.69e-11    &  1.09e-11 & 1.26e-06    &  2.01e-07 & 7.46e+01    &  2.28e+01 \\
512 & 3 & 9.44e-14    &  7.54e-14 & 4.82e-11    &  5.99e-12 & 2.93e-08    &  3.30e-09 & 1.48e+02    &  3.22e+01 \\
% 1024 & 3 & 3.69e-13    &	 1.88e-13 & 4.89e-10    &	 4.23e-11 &  5.63e-10    &	 4.51e-11 & 2.95e+02    &	 4.55e+01 \\
% 2048 & 3 & 6.50e-13    &	 3.67e-13 & --    &	 -- &  1.10e-11    &	 6.45e-13 & 5.90e+02    &	 6.43e+01 \\
\hline
\end{tabular}
}
\caption{The log function defined in Eq. \eqref{def:log_func} is approximated using GFS with various grid sizes $N$. The first derivative approximation error ($e^\prime$) is compared among GFS, finited difference (FD), and FFT methods. Jumps $J_m, m=0,1,\hdots,4n-1,$ are calculated analytically or numerically using one-sided FD schemes with stencil width $M+1=r+4n-1$, ensuring that $J_{4n-1}$ achieves $\mathcal{O}(\Delta x^6)$ accuracy.}
\label{tab:log}
\end{table}
%%%%%%%%%%%%%%%%%%%%%%%%%
%%%%%%%%%%%%%%%%%%%%%%%%%

\subsection{multi-mode aperiodic function}
Consider the function composed of sinusoidal components with $N_k$ non-integer modes ($k_j$), which are non-periodic in the domain $[-\pi,\pi]$,
\begin{equation}
u(x)=\sum_{j=0}^{N_k-1} \sin(k_jx) + \cos(k_jx), \quad  -\pi\le x \le \pi, \quad  k_j = j+\Delta_j, \quad  \Delta_j\equiv\frac{1}{N_k} + \frac{j}{N_k}\left(\frac{N_k-2}{N_k-1}\right) \in \left[\frac{1}{N_k}, 1-\frac{1}{N_k} \right].\label{def:multimode}
\end{equation}
This function is highly oscillatory and enables us to assess the resolution power ($R$) of the GFS, defined as the degrees of freedom needed to resolve one wavelength to machine precision. This can be expressed as $R\equiv\frac{N}{\max k_j}\approx\frac{N}{N_k}$.

\begin{figure}[ht]
 \centering
 \begin{minipage}{0.32\linewidth}\centering
 \vspace{.2in}
 \includegraphics[width=\linewidth,height=1.7in]{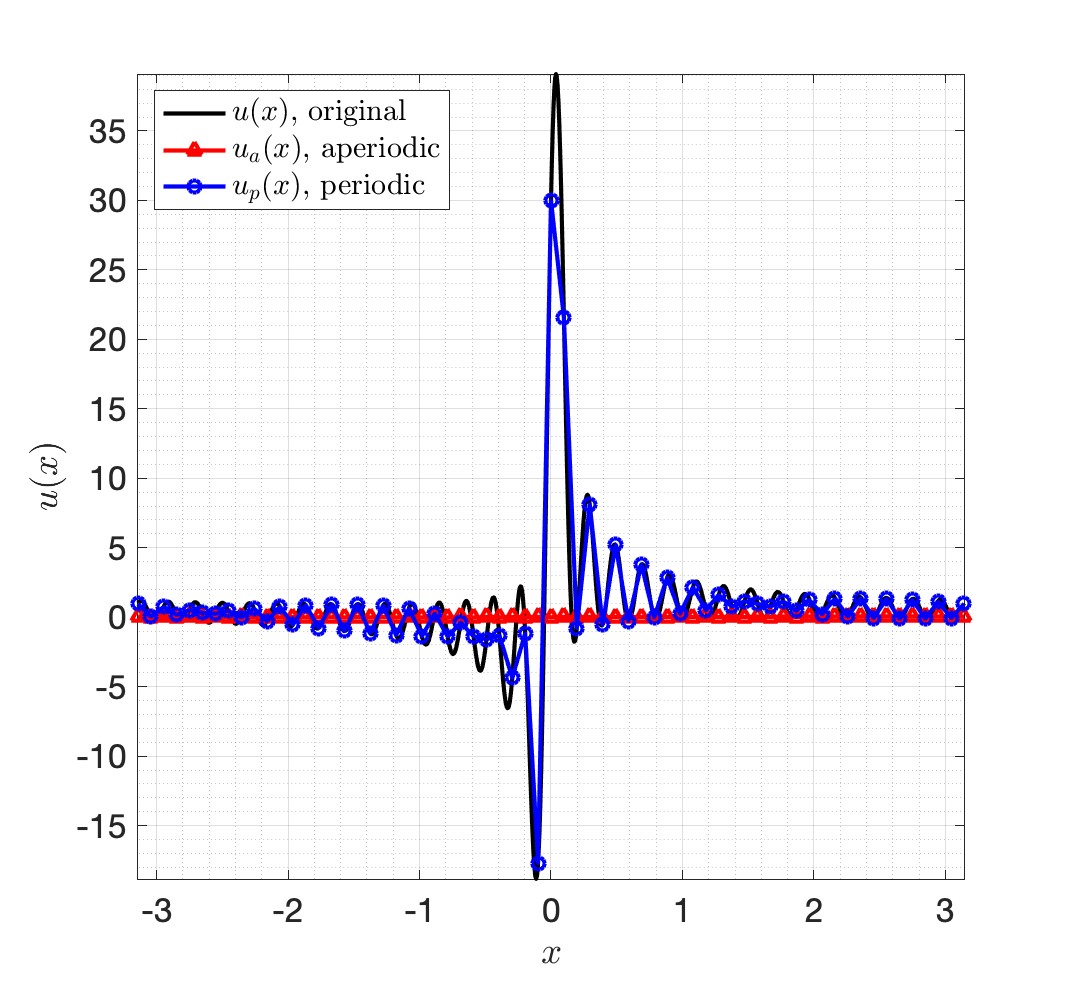} 
\\(a) 
 \end{minipage}~~
 \begin{minipage}{0.32\linewidth}\centering
 \vspace{.2in}
 \includegraphics[width=\linewidth,height=1.7in]{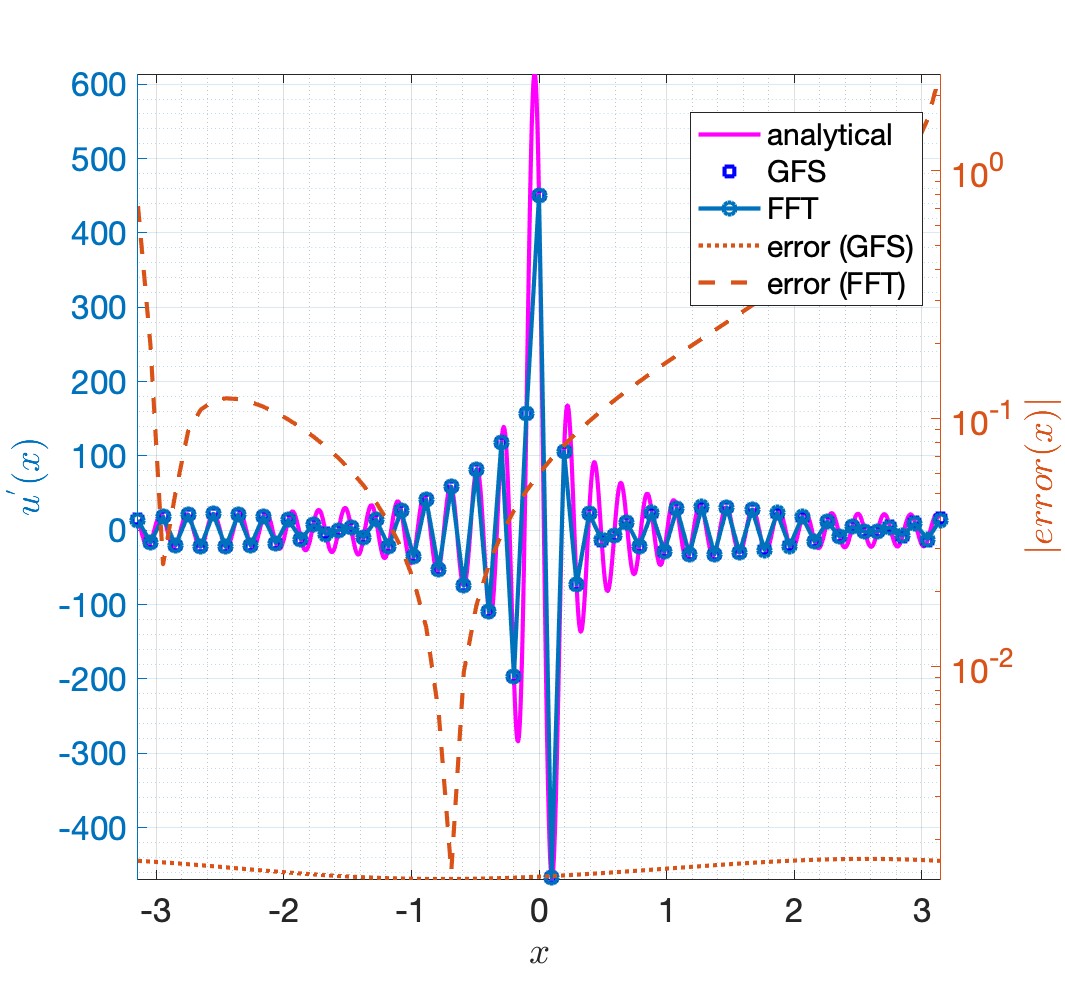} 
\\(b) 
 \end{minipage}
 \begin{minipage}{0.32\linewidth}\centering
 \vspace{.2in}
 \includegraphics[width=\linewidth,height=1.7in]{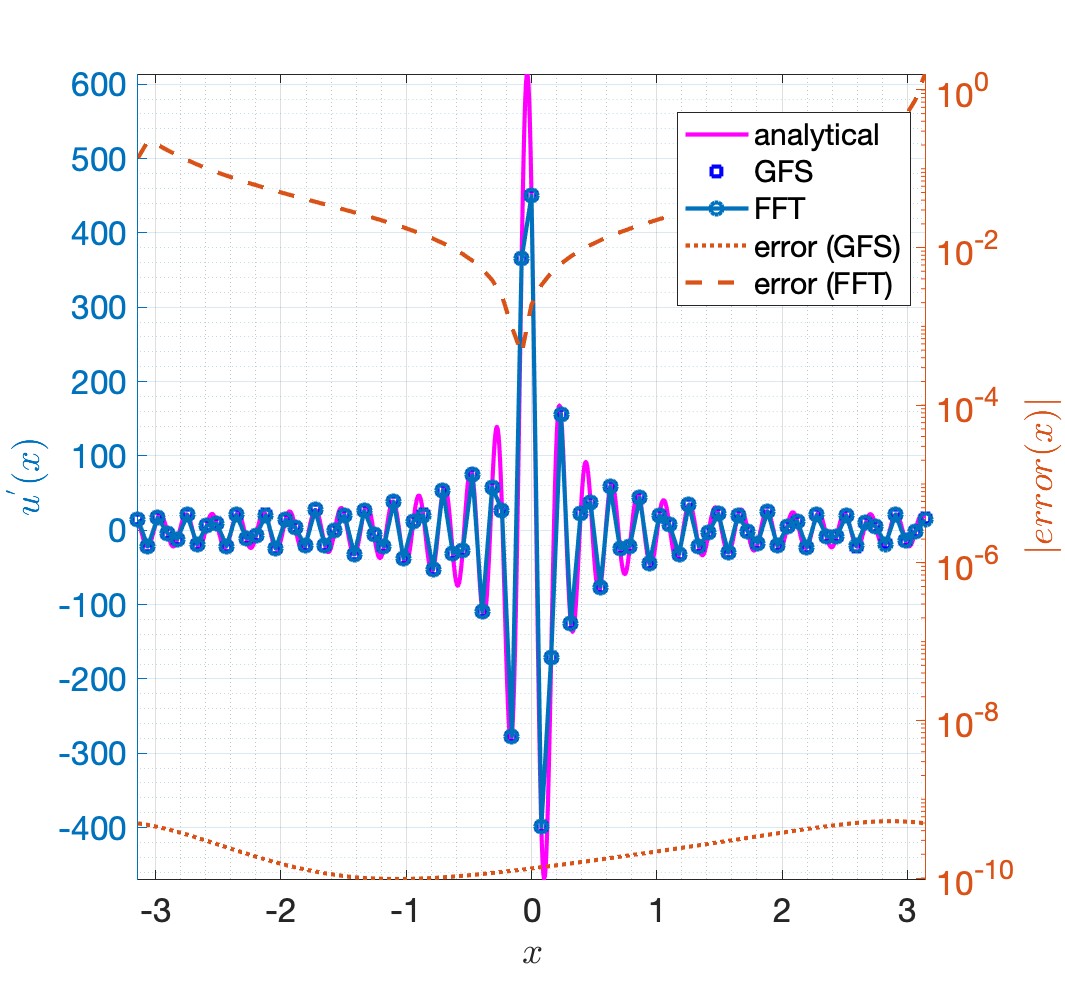} 
\\(c) 
 \end{minipage}
 \caption{A multi-mode sinusoidal function defined by Eq. \eqref{def:multimode} with $N_k=30$ is approximated using GFS with $n=6$. Panel (a) displays its periodic and aperiodic decompositions for $N=64$. Panels (b) and (c) compare its first derivative approximations and the errors from GFS and FFT against analytical results for $N=64$ and $N=80$.
 %(b) $N_k=32, N=64$, $n=6$.
 }
 \label{fig:multiple_aperiodic_modes}
\end{figure}

Fig. \ref{fig:multiple_aperiodic_modes}(a) displays $u(x)$ for $N=64, N_k=30$, demonstrating high oscillation, with the maximum wavenumber $k_{N_k-1}$ near the Nyquist limit ($N_k\sim N/2$). Fig. \ref{fig:multiple_aperiodic_modes}(b, c) compares the first derivative obtained via GFS for $n=6$ with the FFT and analytical methods for $N=64$ and $N=80$. GFS closely aligns with the analytical method, exhibiting no Gibbs oscillations, while FFT has Gibbs oscillations. 
%The $L^\infty, L^2$ norms of the error in first derivative ($||e^\prime||_p$) are quantified for various $n, N$.
For quantitative evaluation, Table \ref{tab:complex_exp} presents the $L^\infty$ and $L^2$ norms of the error in the first derivative for the current method (GFS), a sixth-order finite difference (FD) method, and the FFT. 
We set $N_k=30$ and begin with $N=64$, where the maximum wave number ($k_{N_k-1}$) is close to the Nyquist limit, and examine the impact of increasing resolution and the number of aperiodic modes $n$. Ideally, if $n=N_k$, the aperiodic part would exactly recover the total function $u'(x_j)$. However, this involves higher ($d$) derivatives, with their magnitudes scaling as $k_{N_k-1}^d$, leading to significant increases in precision loss. Moreover, as $n \longrightarrow N_k \gg 1$, this approach becomes costly, so we aim to limit $n$ to $\mathcal{O}(1)$ for practical computations to control expenses and maintain precision.

\begin{figure}[ht]
 \centering
 % \begin{minipage}{\linewidth}\centering
 % \includegraphics[width=3.5in,trim={0.in 0.in 0.in 0in},clip=true]{./Figures/convergence_multimode_2norm.jpg}
 \includegraphics[height=2.5in]{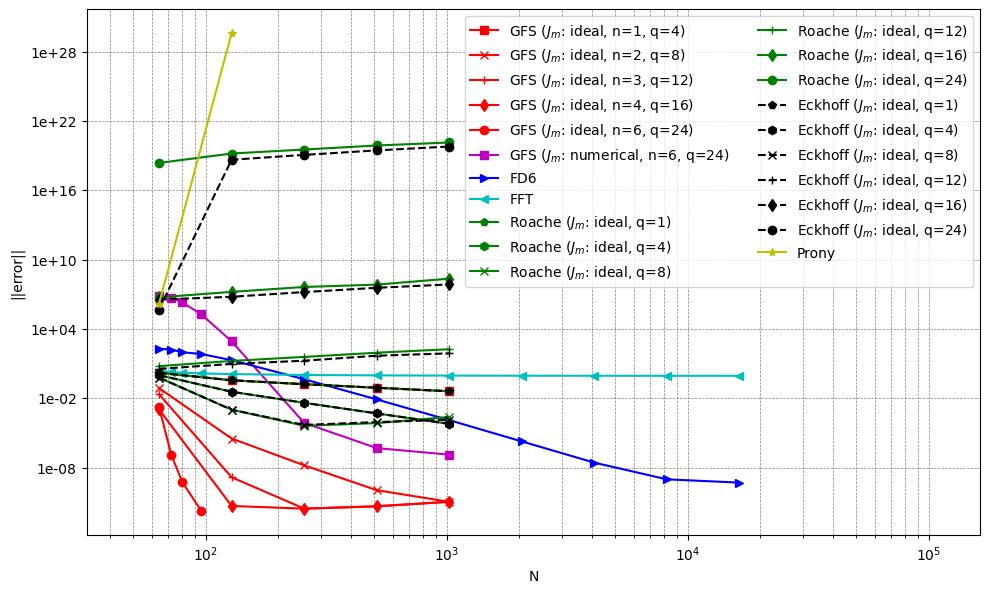}
 % \end{minipage}
  \caption{\rev{For the multi-mode function defined in Eq. \eqref{def:multimode} with $N_k=30$ non-integer modes, convergence of the numerical error ($L^\infty-$ norm) in the first derivative is compared for GFS, a sixth-order finite difference (FD6), FFT, Roache, Eckhoff, and Prony methods. For the Roache and Eckhoff methods, the jumps $J_m$ are evaluated analytically, and for the GFS method, $J_m$ are evaluated both analytically and numerically. Here, ``$J_m:$ ideal'' refers to the case where the jumps $J_m, m=0,1,\hdots,q-1, q=4n,$ are evaluated analytically. GFS shows superior convergence with mesh size $N$. As the number of jumps ($q$) increases, GFS remains robust and converges rapidly, while the other jump-based methods (Roache and Eckhoff) suffer from numerical stability issues due to the function's high oscillations and ill-conditioning at large $q$. The Prony's method ($M=N_k$) is numerically unstable.}}
  \label{fig:convergence_multimode}
\end{figure}

\rev{With $\max u'(x) \approx 600$ (Fig. \ref{fig:multiple_aperiodic_modes}(b,c)), Table \ref{tab:complex_exp} shows that machine precision is achieved with a relatively small number of aperiodic modes ($n=\mathcal{O}(1)$) using GFS. The error decreases much more rapidly with mesh refinement (increasing $N$) for GFS (with analytical $J_m$) even at $n=2$, and improves markedly as $n$ grows. 
% This is expected due to the $C^{4n-1}$ continuity of $u_{p}(x)$, which allows the periodic part to be smooth at the domain endpoints up to $4n$ derivatives, leading to  Fourier series convergence as $1/k^{4n+1}$. 
While the sixth-order FD scheme's error scales as $\mathcal{O}(\Delta x^6)$, the FFT error remains $\mathcal{O}(1)$ in non-periodic domains due to the  Gibbs phenomenon. For $N_k=30$ and a fixed target error of $\mathcal{O}(10^{-10})$, GFS requires about $N=80-90$, yielding a resolution power $R\approx 3$; in contrast, the sixth-order FD scheme needs $N=8192$ (see Fig. \ref{fig:convergence_multimode}), giving $R=300$. Thus, the FD mesh must be two orders of magnitude finer than GFS for the same accuracy, evidencing GFS's superior resolution. Even with a numerically evaluated $J_m$, GFS maintains higher accuracy and convergence than FD and FFT methods, see Table \ref{tab:complex_exp}.

Fig. \ref{fig:convergence_multimode} shows the GFS method converging faster and with lower error than the FD, FFT, Roache's, Eckhoff's, and Prony's methods. When $J_m$ is evaluated analytically, GFS rapidly improves with increasing $N$, surpassing all other methods. Even with numerical $J_m$, GFS remains markedly superior for large $N$. Polynomial-based methods like the Roache and Eckhoff approaches become ill-conditioned as the number of jumps ($q$) in the discontinuous part increases, while GFS stays robust. The Prony's method is ill-conditioned due to high oscillations, causing the computation to fail with $M=N_k$ and $M=N/2$ (not shown) in Eq. \eqref{eqn:prony}.}
  
In summary, the convergence rate of GFS improves with the number of aperiodic modes $n$, and GFS has a superior resolution power ($R\approx 3$) compared to the sixth-order finite difference scheme ($R\approx 300$). 
\iffalse
\begin{itemize}
\item When the grid is coarse but all the aperiodic modes are resolved, accuracy is only limited by the roundoff errors, i.e., machine precision is possible if the mode numbers are low; otherwise, roundoff effects due to highly ill-conditioned matrices may deteriorate the accuracy.
\item When the grid is coarse, and all the aperiodic modes are not resolved ($n\neq N_k, N_k=N/2$), the error is large. However, the error converges exponentially/algebraically ($\Delta x ^{r(n)}$) with grid refinement.
\item Low mode aperiodicity can be accurately resolved by increasing $n$.
\item High mode aperiodicity can be better resolved with grid refinement, where low modes are taken care of through $\kc_j, \ks_j,$ and high mode aperiodicity is projected onto the Fourier basis/normal modes.
\item Resolution power is a function of $n/N_k$. For $n=N_k$,the resolution power is 3 degrees of freedom per wavelength.
\end{itemize}
\fi
%%%%%%%%%%%%%%%%%%%%%%%%%
%%%%%%%%%%%%%%%%%%%%%%%%%

\begin{table}
\centering
\resizebox{\linewidth}{!}{
\begin{tabular}{|cc||c|cc|cc||cc||cc| }
\hline
 &  &	 \multicolumn{5}{c||}{GFS}	& \multicolumn{2}{c||}{FD: $\mathcal{O}(\Delta x^r)$}& \multicolumn{2}{c|}{FFT}\\
 \hline
 & & &  \multicolumn{2}{c|}{$J_m:$ analytical} & \multicolumn{2}{c||}{$J_m:$ numerical} & \multicolumn{2}{c||}{$r=6$} & \multirow{2}{*}{} & \multirow{2}{*}{} \\
\hline
$N$ & $N_k$ &	$n$ & $||e^\prime||_\infty$ & $||e^\prime||_2$ & $||e^\prime||_\infty$ & $||e^\prime||_2$ & $||e^\prime||_\infty$ & $||e^\prime||_2$ & $||e^\prime||_\infty$ & $||e^\prime||_2$ \\
\hline
64 & 30 & 2 & 7.34e-02    & 1.16e-01 & 5.52e+03 & 1.74e+03 & 1.93e+02 & 1.35e+02 & 2.45e+00 & 1.24e+00 \\
128 & 30 & 2 & 3.23e-06    & 1.80e-06 & 6.19e+01 & 1.43e+01 & 1.95e+01 & 7.23e+00 & 1.22e+00 & 3.40e-01 \\
256 & 30 & 2 & 1.55e-08    & 5.50e-09 & 3.61e-02 & 5.92e-03 & 4.40e-01 & 1.41e-01 & 1.03e+00 & 2.23e-01 \\
512 & 30 & 2 & 1.08e-10    & 2.64e-11 & 1.12e-05 & 1.27e-06 & 7.84e-03 & 2.21e-03 & 9.46e-01 & 1.55e-01 \\
% 1024 & 30 & 2 & 1.84e-11    & 9.46e-12 & 2.14e-09 & 1.72e-10 & 1.26e-04 & 3.36e-05 & 9.05e-01 & 1.09e-01 \\
% 2048 & 30 & 2 & 3.74e-11    & 1.96e-11 & 2.55e-09 & 1.58e-10 & 1.81e-06 & 5.15e-07 & 8.85e-01 & 7.72e-02 \\
\hline
64 & 30 & 4 & 1.06e-03    & 1.54e-03 & 6.97e+05 & 2.63e+05 & 1.93e+02 & 1.35e+02 & 2.45e+00 & 1.24e+00 \\
128 & 30 & 4 & 5.24e-12    & 3.56e-12 & 3.56e+02 & 9.98e+01 & 1.95e+01 & 7.23e+00 & 1.22e+00 & 3.40e-01 \\
256 & 30 & 4 & 3.79e-12    & 2.27e-12 & 1.81e-03 & 2.85e-04 & 4.40e-01 & 1.41e-01 & 1.03e+00 & 2.23e-01 \\
512 & 30 & 4 & 7.12e-12    & 4.77e-12 & 2.73e-08 & 3.81e-09 & 7.84e-03 & 2.21e-03 & 9.46e-01 & 1.55e-01 \\
% 1024 & 30 & 4 & 1.88e-11    & 8.96e-12 & 2.69e-08 & 2.11e-09 & 1.26e-04 & 3.36e-05 & 9.05e-01 & 1.09e-01 \\
% 2048 & 30 & 4 & 3.52e-11    & 1.88e-11 & -- & -- & 1.81e-06 & 5.15e-07 & 8.85e-01 & 7.72e-02 \\
\hline
64 & 30 & 6 & 1.68e-03    & 3.86e-03 & 6.57e+06 & 2.74e+06 & 1.93e+02 & 1.35e+02 & 2.45e+00 & 1.24e+00 \\
72 & 30 & 6 & 1.16e-07    & 1.97e-07 & 5.04e+06 & 1.50e+06 & 1.61e+02 & 9.46e+01 & 1.81e+00 & 6.94e-01 \\
80 & 30 & 6 & 5.27e-10    & 7.34e-10 & 1.96e+06 & 5.59e+05 & 9.50e+01 & 6.48e+01 & 1.59e+00 & 5.53e-01 \\
96 & 30 & 6 & 1.64e-12    & 1.46e-12 & 1.79e+05 & 4.68e+04 & 6.67e+01 & 3.00e+01 & 1.39e+00 & 4.36e-01 \\
128 & 30 & 6 & 2.03e-12    & 1.39e-12 & 8.34e+02 & 2.19e+02 & 1.95e+01 & 7.23e+00 & 1.22e+00 & 3.40e-01 \\
256 & 30 & 6 & 2.74e-12    &	 1.45e-12       & 7.29e-05 & 1.20e-05 & 4.40e-01 & 1.41e-01 & 1.03e+00 & 2.23e-01 \\
512 & 30 & 6 &  4.40e-12    &	 3.15e-12       & 4.80e-07 & 6.20e-08 & 7.84e-03 & 2.21e-03 & 9.46e-01 & 1.55e-01 \\
% 1024 & 30 & 6 & 1.09e-11    &	 6.00e-12       & 1.29e-07 & 1.02e-08 & 1.26e-04 & 3.36e-05 & 9.05e-01 & 1.09e-01 \\ 
% 2048 & 30 & 6 & --         & --       & --       & --       & 1.81e-06 & 5.15e-07 & 8.85e-01 & 7.72e-02 \\
% 4096 & 30 & 6 & --         & --       & --       & --       & 2.64e-08 & 7.96e-09 & 8.75e-01 & 5.46e-02 \\
% 8192 & 30 & 6 & --         & --       & --       & --       & 9.47e-10 & 1.37e-10 & 8.70e-01 & 3.86e-02 \\ 
% 16384 & 30 & 6 & --        & --       & --       & --       & 4.90e-10 & 1.11e-10 & 8.68e-01 & 2.73e-02 \\ 
\hline
\end{tabular}
}
\caption{A multi-mode sinusoidal function defined by Eq. \eqref{def:multimode} with $N_k=30$ is approximated using GFS, finite difference (FD), and FFT methods. The error of the approximation for the first derivative ($e^\prime$) is shown with increasing number of aperiodic modes ($n$). The jumps $J_m$ are computed either analytically or numerically through one sided finite difference schemes of $r^{th}$ order accuracy as indicated.}
\label{tab:complex_exp}
\end{table}

\subsection{Monomials}
We consider the monomial functions $u(x) = x^m, -\pi\le x\leq \pi$ for $m=1,3$, which are non-periodic and have finite jumps in the function values at the domain endpoints. 
%as below,
%\begin{equation}
%u(x) = \begin{cases}x^m, \quad \quad -\pi\le x\leq \pi \\0, \quad \quad \text{otherwise}.
%\end{cases}
%\end{equation}
Although the GFS is not primarily designed for approximating polynomials, it can still accommodate such functions by regularizing jumps. The GFS relies on endpoint conditions related to jumps, which can quickly diminish for monomial functions, resulting in rank-deficient matrices and vanishing right-hand sides. This issue can be resolved by regularizing jumps using very small numerical values ($10^{-15}$ in the following cases) when the actual jump is zero. It's important to note that regularization is only necessary for the analytical case; for numerical computations, the inherent error of the numerical scheme provides the required regularization. 

\begin{figure}[ht]
 \centering
 \begin{minipage}{0.45\linewidth}\centering
 \includegraphics[width=2.25in,trim={0.in 0.in 0.in 0in},clip=true]{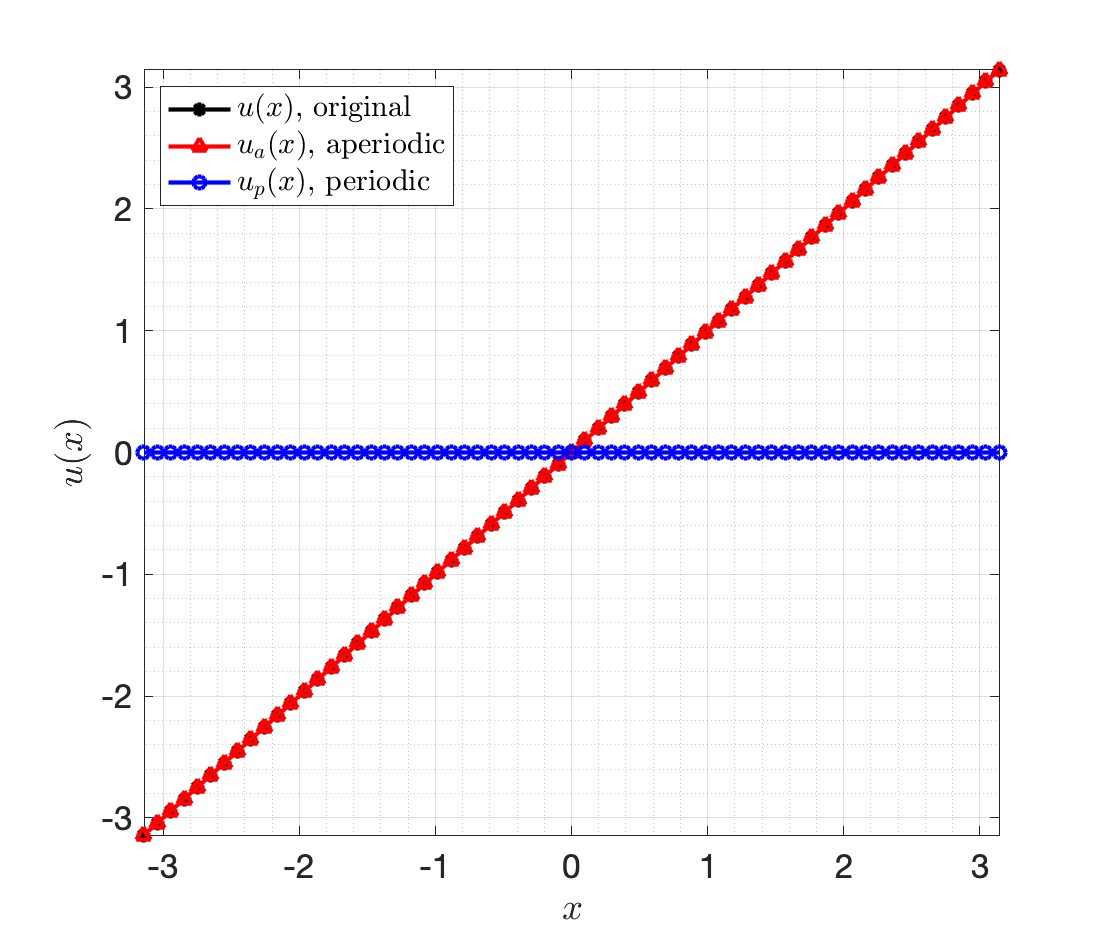}\\(a)   \end{minipage}~~
\begin{minipage}{0.5\linewidth}\centering
 \includegraphics[width=2.25in,trim={0.in 0.in 0.in 0in},clip=true]{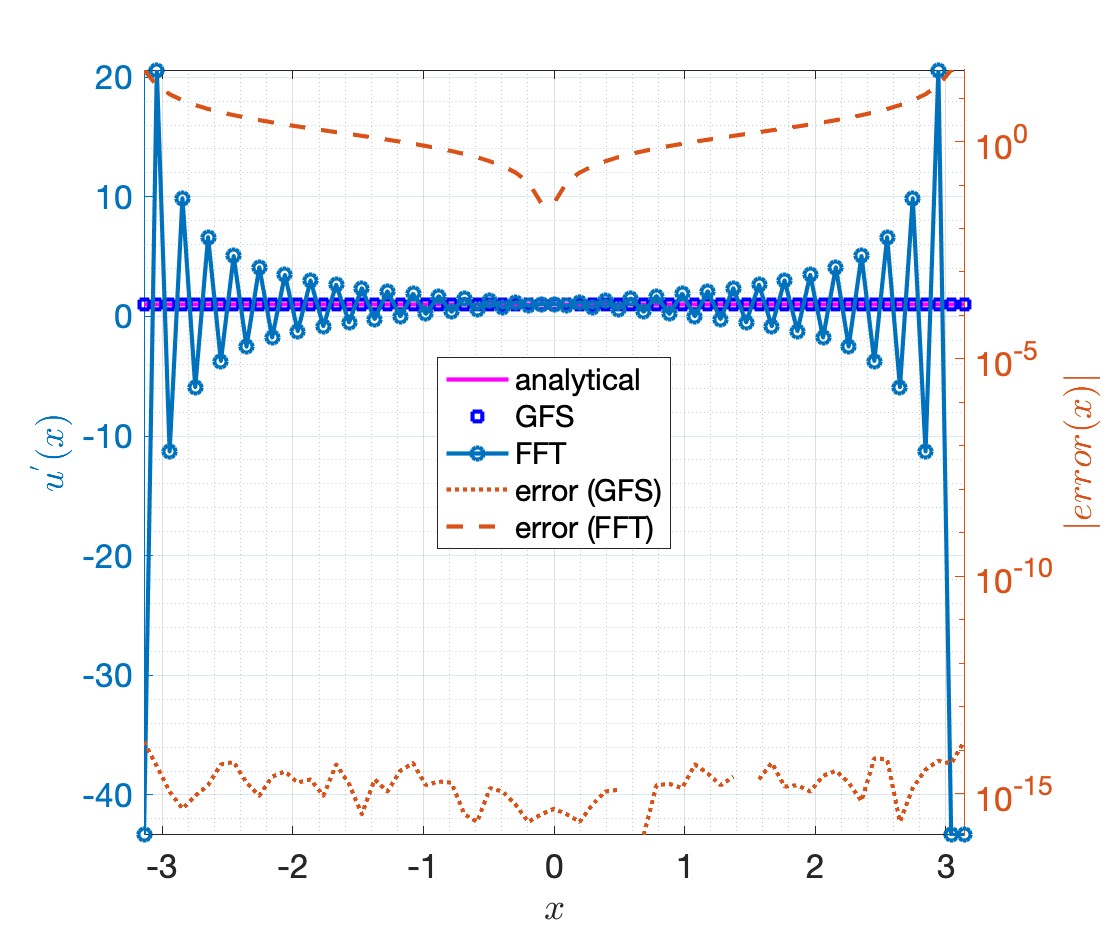}\\(b)
 \end{minipage}~~
  \caption{Ramp function $u(x)=x$ is approximated with GFS using $n=1, N=64$: Decomposition of $u(x)$ into periodic and aperiodic parts is shown in (a), approximation of the first derivative with GFS and FFT is compared with analytical values in (b); GFS does not exhibit Gibbs oscillations and aligns well with the analytic values while FFT shows Gibbs oscillations.}
  \label{fig:ramp}
\end{figure}

For $m=1$, this corresponds to the ramp function, a well-known test case that challenges Fourier spectral methods. Fig. \ref{fig:ramp}(a) shows the decomposition of $u(x)$ into periodic $u_p(x)$ and aperiodic $u_a(x)$ components obtained by GFS for $n=1$. Fig. \ref{fig:ramp}(b) compares the first derivative approximations of $u(x)$ using GFS and FFT; GFS avoids Gibbs oscillations and closely matches analytical values, whereas FFT exhibits them. Table \ref{tab:error_ramp} compares the $L^\infty, L^2$ norms of the first derivative error $||e^\prime||_p$ for GFS, FD and FFT methods at $N=64$. The numerical error for $J_m$ and the stand-alone FD method approaches machine precision, as the sixth-order FD schemes fully resolve the linear function $u(x)=x$. 

\begin{table}
\centering
\resizebox{\linewidth}{!}{
\begin{tabular}{|c||c|cc|cc||cc||cc| }
\hline
 &  \multicolumn{5}{c||}{GFS}	& \multicolumn{2}{c||}{FD: $\mathcal{O}(\Delta x^r)$}& \multicolumn{2}{c|}{FFT}\\
\hline
 & &  \multicolumn{2}{c|}{$J_m:$ analytical} & \multicolumn{2}{c||}{$J_m:$ numerical} & \multicolumn{2}{c||}{$r=6$} & \multirow{2}{*}{} & \multirow{2}{*}{} \\
\hline
$N$ &$n$ & $||e^\prime||_\infty$ & $||e^\prime||_2$ &$||e^\prime||_\infty$ & $||e^\prime||_2$ & $||e^\prime||_\infty$ & $||e^\prime||_2$ &$||e^\prime||_\infty$ & $||e^\prime||_2$ \\
\hline
64 & 1 & 1.55e-14    &	 9.50e-15 & 3.03e-14    &	 1.47e-14 & 9.38e-14    &	 3.00e-14    & 4.43e+01    &	 2.71e+01 \\
64 & 2 & 7.22e-15    &	 7.98e-15 & 6.57e-13    &	 2.91e-13 & 9.38e-14    &	 3.00e-14 & 4.43e+01    &	 2.71e+01 \\
64 & 3 & 1.03e-14    &	 7.97e-15 & 1.57e-12    &	 9.18e-13 & 9.38e-14    &	 3.00e-14 & 4.43e+01    &	 2.71e+01 \\
\hline
\end{tabular}
}
\caption{The ramp function $u(x)=x, -\pi\leq x\leq \pi$ is approximated with the GFS. Error of the approximation for the function ($e$) and its first derivative ($e^\prime$) are shown with increasing number of $\sin$ or $\cos$ modes ($n$). Here, $N=64$ and the jump conditions $J_m$ are computed from exact derivatives or numerical derivatives (one-sided finite difference method) as indicated.}
\label{tab:error_ramp}
\end{table}

% ------------------------------
%% ------------------------------
%\begin{figure}[ht]
% \centering
% \begin{minipage}{0.45\linewidth}\centering
% \includegraphics[width=2.75in,trim={0.in 0.in 0.in 0in},clip=true]{./Figures/{u_x_n2_N64}.jpg}   \end{minipage}~~
%\begin{minipage}{0.45\linewidth}\centering
% \includegraphics[width=3in,trim={0.in 0.in 0.in 4.5in},clip=true]{./Figures/{up_x_n2_N64}.jpg}
% \end{minipage}~~
%  \caption{Ramp function $u(x)=x$ approximated with $n=2, N=64$.}
%\end{figure}
%% ------------------------------
%\begin{figure}[ht]
% \centering
% \begin{minipage}{0.45\linewidth}\centering
% \includegraphics[width=2.75in,trim={0.in 0.in 0.in 0in},clip=true]{./Figures/{u_ramp_3}.jpg}   \end{minipage}~~
%\begin{minipage}{0.45\linewidth}\centering
% \includegraphics[width=3in,trim={0.in 0.in 0.in 4.5in},clip=true]{./Figures/{up_ramp_3}.jpg}
% \end{minipage}~~
%  \caption{Ramp function $u(x)=x$ with $n=3$.}
%\end{figure}

\begin{figure}[ht]
 \centering
 \begin{minipage}{0.45\linewidth}\centering
 \includegraphics[width=2.25in,trim={0.in 0.in 0.in 0in},clip=true]{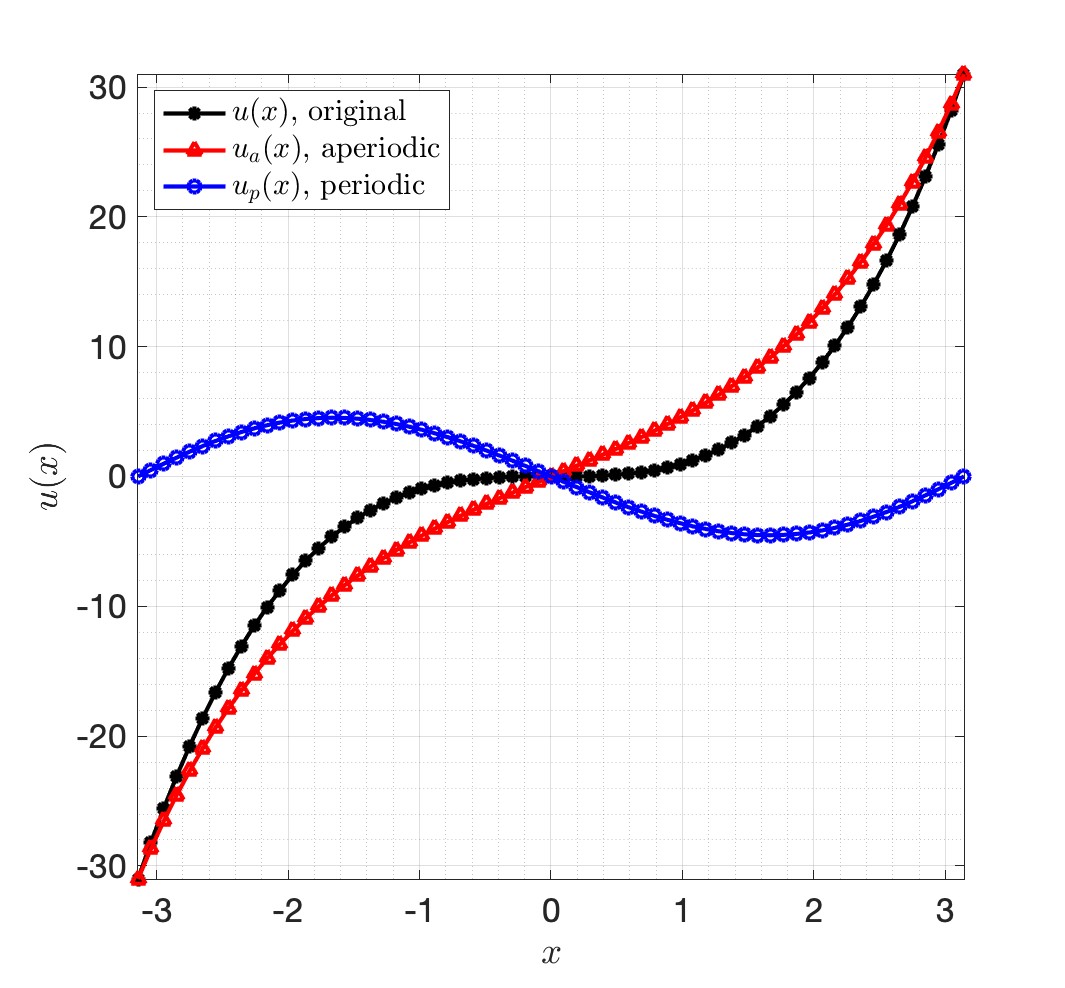} \\(a)  \end{minipage}~~
\begin{minipage}{0.5\linewidth}\centering
 \includegraphics[width=2.25in,trim={0.in 0.in 0.in 0in},clip=true]{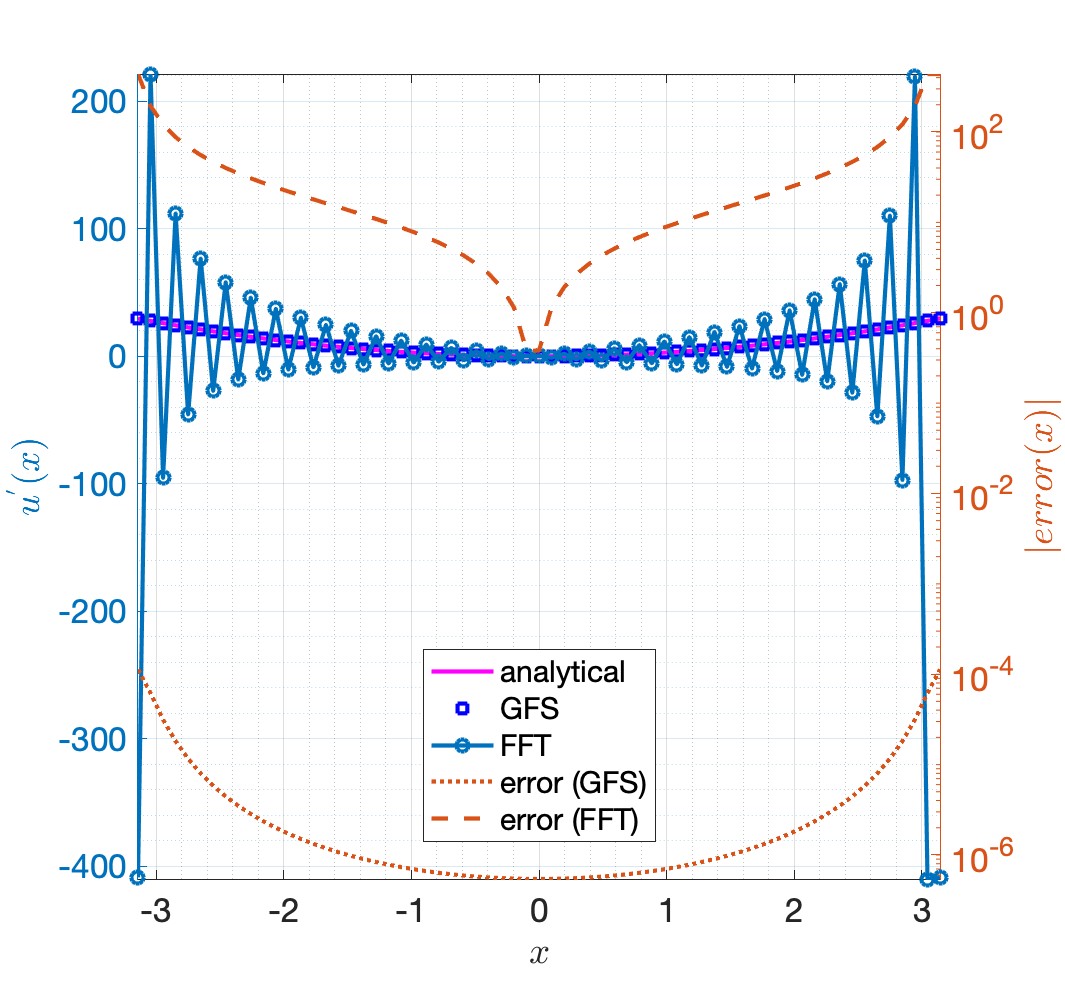}\\(b)
 \end{minipage}~~\\
 \begin{minipage}{0.45\linewidth}\centering
 \includegraphics[width=2.25in,trim={0.in 0.in 0.in 0in},clip=true]{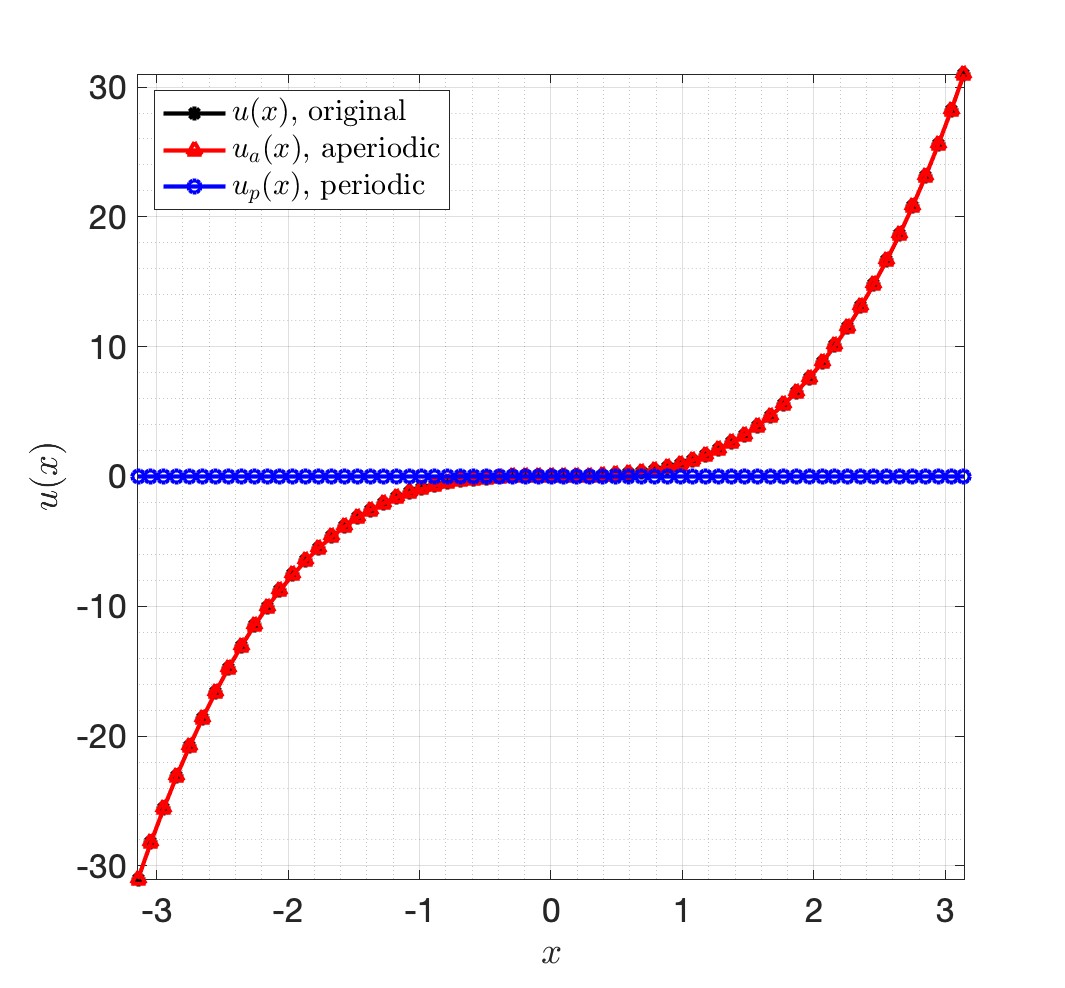} \\(c)  \end{minipage}~~
\begin{minipage}{0.5\linewidth}\centering
 \includegraphics[width=2.25in,trim={0.in 0.in 0.in 0in},clip=true]{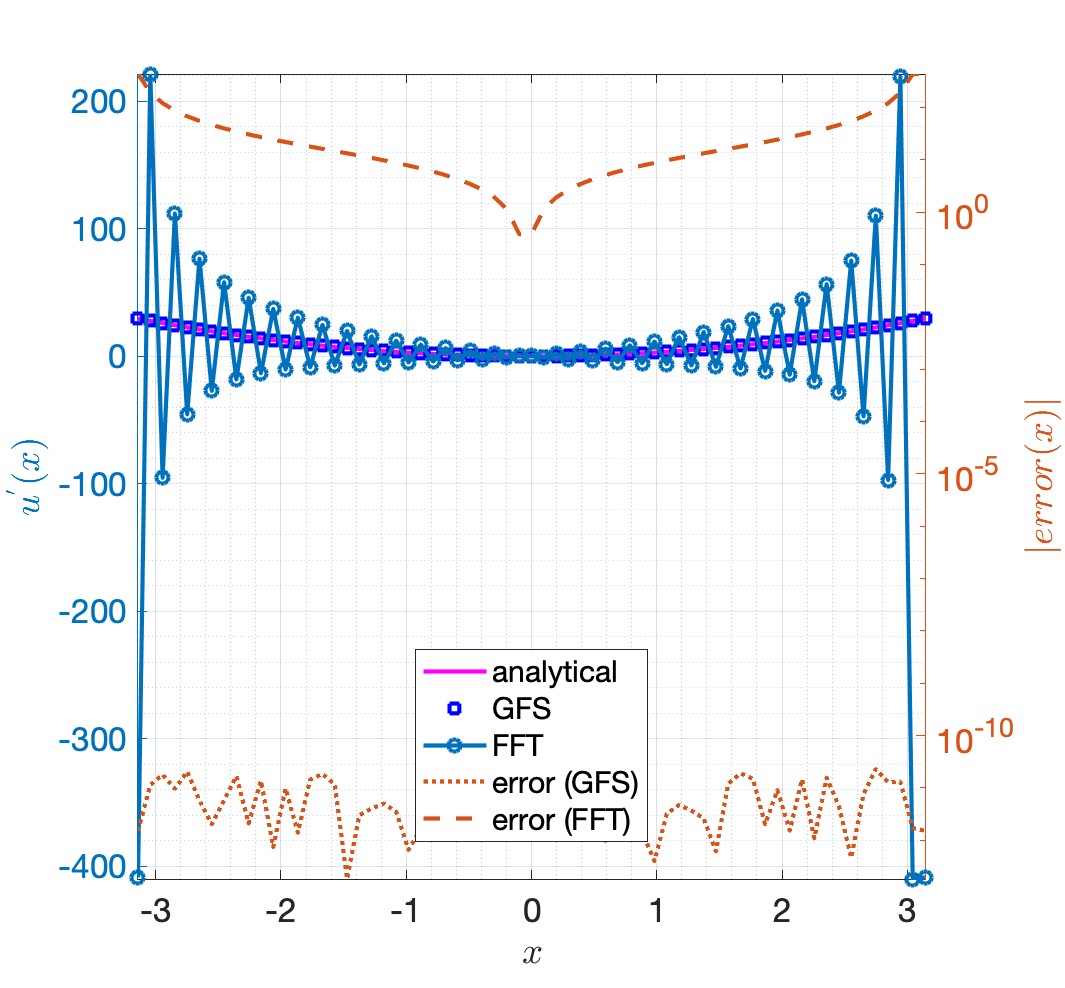} \\(d)
 \end{minipage}~~
 
  \caption{Cubic monomial $u(x)=x^3$ is approximated with GFS using $N=64$ and $n=1$ (a-b), $n=2$ (c-d): Decomposition of $u(x)$ into periodic and aperiodic parts is shown in (a,c), approximation of the first derivative with GFS and FFT is compared with analytical values in (b,d); GFS does not exhibit Gibbs oscillations and aligns well with the analytic values while FFT shows Gibbs oscillations.}
  \label{fig:cubic_n1}
\end{figure}

For $m=3$, this corresponds to a cubic monomial. Fig. \ref{fig:cubic_n1}(a,c) shows the decomposition of $u(x)$ into periodic $u_p(x)$ and aperiodic $u_a(x)$ components using GFS for $n=1$ and 2, respectively. Fig. \ref{fig:cubic_n1}(b,d) compares the first derivative approximations of $u(x)$ from GFS and FFT; GFS avoids Gibbs oscillations and aligns closely with analytical values, while FFT exhibits oscillations leading to $\mathcal{O}(1)$ error. Table \ref{tab:error_cubic} compares the $L^\infty$ and $L^2$ norms of the first derivative error $||e^\prime||_p$ for GFS, FD and FFT methods at $N=64$. 
% The analytically calculated jumps ($J_m$) are regularized with a small value of $\epsilon=10^-15$ as higher derivatives vanish. 
The numerical error for GFS converges quickly with increasing $n$, whereas the sixth-order FD scheme achieves near machine precision by exactly approximating the monomials for $m\le 5$. 
% Thus, numerical evaluation of $J_m$ with the sixth-order FD scheme does not compromise the overall accuracy of GFS.

\begin{table}
\centering
\resizebox{\linewidth}{!}{
\begin{tabular}{|c||c|cc|cc||cc||cc| }
\hline
 &  \multicolumn{5}{c||}{GFS}	& \multicolumn{2}{c||}{FD: $\mathcal{O}(\Delta x^r)$}& \multicolumn{2}{c|}{FFT}\\
\hline
 & &  \multicolumn{2}{c|}{$J_m:$ analytical} & \multicolumn{2}{c||}{$J_m:$ numerical} & \multicolumn{2}{c||}{$r=6$} & \multirow{2}{*}{} & \multirow{2}{*}{} \\
\hline
$N$ &$n$ & $||e^\prime||_\infty$ & $||e^\prime||_2$ &$||e^\prime||_\infty$ & $||e^\prime||_2$ & $||e^\prime||_\infty$ & $||e^\prime||_2$ &$||e^\prime||_\infty$ & $||e^\prime||_2$ \\
\hline
64 & 1 & 1.13e-04    &	 5.98e-05 & 1.13e-04    &\ 5.98e-05 & 7.07e-13    &\ 3.33e-13 & 4.38e+02    &\ 2.68e+02 \\
64 & 2 & 6.70e-06    &	 6.06e-06 & 2.77e-11    &	 1.94e-11 & 7.07e-13    &\ 3.33e-13 & 4.38e+02    &\ 2.68e+02 \\
64 & 3 & 8.32e-09    &	 6.90e-09 & 1.22e-10    &	 5.75e-11 & 7.07e-13    &\ 3.33e-13 & 4.38e+02    &\ 2.68e+02 \\
\hline
\end{tabular}
}
\caption{A function $u(x)=x^3, -\pi\le x \le \pi$, is approximated using GFS with $N=64$. The error of the approximation for the first derivative ($e^\prime$) is compared among various methods. Here, the jump conditions $J_m$ are computed either exactly or numerically through one-sided finite difference (FD) schemes with $\mathcal{O}(\Delta x^6)$ accuracy as indicated in the column $J_m$.}
\label{tab:error_cubic}
\end{table}

% ------------------------------
% ------------------------------
% \begin{figure}[ht]
%  \centering
%  \begin{minipage}{0.45\linewidth}\centering
%  \includegraphics[width=2.5in,trim={0.in 0.in 0.in 0in},clip=true]{./Figures/u_x3_n2_N64.jpg} \\(a)  \end{minipage}~~
% \begin{minipage}{0.5\linewidth}\centering
%  \includegraphics[width=2.5in,trim={0.in 0.in 0.in 0in},clip=true]{./Figures/up_x3_n2_N64.jpg} \\(b)
%  \end{minipage}~~
%   \caption{Cubic polynomial function $u(x)=x^3$ approximated with $n=2, N=64$.}
%   \label{fig:cubic_n2}
% \end{figure}
%% ------------------------------
%\begin{figure}[ht]
% \centering
% \begin{minipage}{0.45\linewidth}\centering
% \includegraphics[width=2.75in,trim={0.in 0.in 0.in 0in},clip=true]{./Figures/{u_x3_3}.jpg}   \end{minipage}~~
%\begin{minipage}{0.45\linewidth}\centering
% \includegraphics[width=3in,trim={0.in 0.in 0.in 4.5in},clip=true]{./Figures/{up_x3_3}.jpg}
% \end{minipage}~~
%  \caption{Cubic polynomial function $u(x)=x^3$ with $n=3$.}
%  \label{fig:cubic_n3}
%\end{figure}
%% ------------------------------
%\begin{figure}[ht]
% \centering
% \begin{minipage}{0.45\linewidth}\centering
% \includegraphics[width=2.75in,trim={0.in 0.in 0.in 0in},clip=true]{./Figures/{u_x3_4}.jpg}   \end{minipage}~~
%\begin{minipage}{0.45\linewidth}\centering
% \includegraphics[width=3in,trim={0.in 0.in 0.in 4.5in},clip=true]{./Figures/{up_x3_4}.jpg}
% \end{minipage}~~
%  \caption{Cubic polynomial function $u(x)=x^3$ with $n=4$.}
%  \label{fig:cubic_n4}
%\end{figure}

\section{Conclusions and future directions}\label{sec:conclusions}
\rev{
In this work, we introduce the Generalized Fourier Series (GFS), a spectral method that extends Fourier techniques to non-periodic domains while retaining FFT-level efficiency. GFS decomposes a function into a periodic component—approximated with the Fast Fourier Transform—and an aperiodic component captured by a small set $(n=\mathcal{O}(1))$ of adaptive, non-harmonic sinusoids that handle non-periodicity with high-order accuracy. This construction avoids artificial domain extensions and their overhead, seamlessly accommodates derivative jumps and boundary discontinuities, and delivers high-resolution accuracy and robustness. With a resolution power (degrees of freedom per wavelength) comparable to FFT on periodic problems, and a low-rank approximation for the aperiodic part, the overall complexity remains linear, making GFS a practical and efficient alternative to existing techniques.

% In this work, we introduced the Generalized Fourier Series (GFS), a spectral method designed to extend the applicability of Fourier techniques to non-periodic domains while retaining their hallmark efficiency. By decomposing a function into periodic and aperiodic parts, GFS combines the speed of FFT with an adaptive representation of the non-periodic component, requiring only $n=\mathcal{O}(1)$ additional modes. This construction avoids artificial domain extensions and the associated computational overhead, while delivering high-resolution accuracy and robustness.

Comprehensive numerical experiments show that GFS consistently outperforms classical approaches—including Roache’s polynomial-based correction, Eckhoff’s reconstruction with Bernoulli polynomials, Prony’s exponential fitting, and standard finite-difference and FFT methods. Relative to Roache’s method, GFS avoids a fixed polynomial correction—which can be unstable for oscillatory functions—in favor of adaptive sinusoidal modes that yield greater robustness and accuracy; relative to Eckhoff’s method, it retains the ability to capture discontinuities while offering greater robustness and higher accuracy, particularly when a large number of jumps must be resolved; and relative to Prony’s method, it does not require the precise number of modes to be specified a priori and remains numerically stable at large mesh sizes ($N>64$), where Prony’s scheme typically becomes ill-conditioned. Across meshes, GFS delivered significantly better convergence rates and resolution power: it achieves FFT-like resolution, while a sixth-order finite-difference scheme demands far more resources—underscoring GFS’s orders-of-magnitude performance advantage. While GFS attains optimal accuracy when exact jump conditions are known, approximating these jumps with finite differences can limit overall accuracy; the quality of jump evaluation depends on the underlying numerical scheme and could be improved with better estimators. A detailed treatment of jump-condition estimation lies beyond the scope of this article and will be addressed in future work.

% Comprehensive numerical experiments demonstrate that GFS consistently outperforms classical approaches, including Roache’s polynomial-based method, Eckhoff’s spectral reconstruction, and Prony’s exponential fitting scheme. Specifically: Compared to Roache’s method, GFS avoids the fixed polynomial correction, which can be unstable for oscillatory functions, and instead employs adaptive sinusoidal modes that provide greater robustness and accuracy. Compared to Eckhoff’s method, GFS retains the ability to capture discontinuities but does so in a stable manner, without the severe sensitivity to noise. Compared to Prony’s method, GFS avoids the need to specify the number of modes a priori, and its numerical stability allows it to scale effectively to large mesh sizes ($N > 64$), where Prony’s method typically becomes ill-conditioned.

Across a range of test functions—including modulated sines, Gaussians, and logarithmic functions—GFS showed rapid convergence, minimal Gibbs oscillations, and superior accuracy in approximating derivatives. The results confirm that GFS provides a robust, efficient, and accurate framework for non-periodic function approximation, offering a compelling alternative to existing methods for applications in numerical PDEs, signal processing, and beyond.
}

\rev{Future work will expand the GFS framework to multi-dimensional applications, addressing PDEs such as advection-diffusion and Poisson equations, and eventually tackling computational fluid dynamics problems. These efforts highlight the transformative potential of GFS for enhancing spectral methods in nonperiodic domains.}

\appendix
\section{Solution of the non-linear systems Eqs. \eqref{eqn:ujkj_sin} and \eqref{eqn:ujkj_cos} for $n=1,2,3$}\label{sec:proof}
\subsection{Single mode representation of the aperiodic part: $n=1$}
For $n=1$, the solution of the above equations is straightforward and involves the jumps in the derivatives of $u(x)$ at the endpoints $J_0, J_1, J_2, J_3$. First, we solve for the modes $\ks_1, \kc_1$, followed by the amplitudes $\us_1, \uc_1$ as below,
\begin{eqnarray}
2\us_1 \sin(\ks_1 \pi) = J_0 = -\frac{J_2}{\ks_1^2},&&
-2\kc_1\uc_1 \sin(\kc_1 \pi) = J_1 = - \frac{J_3}{\kc_1^2},\\
\ks_1^2 = -J_2/J_0, \quad \quad \kc_1^2 = -J_3/J_1,&& 
\us_1 = \frac{J_0}{2\sin(\ks_1\pi)}, \quad \uc_1 = -\frac{J_1}{2\kc_1\sin(\kc_1\pi)}. \label{eqn:es1ec1}
\end{eqnarray}
Note that $\kc_1 = \pm \sqrt{-J_3/J_1}$ and $\ks_1 = \pm \sqrt{-J_2/J_0}$ and the $\pm$ sign is irrelevant here because $\uc_j \cos(\kc_j x)$ and $\us_j \sin(\ks_j x)$ (therefore, $u_c(x)$ and $u_s(x)$) are even functions of $\kc$ and $\ks$, respectively. We define the complex square root as $\sqrt{z}\equiv \sqrt{|z|}e^{i\theta/2}$ where $z=|z|e^{i\theta}, -\pi\le \theta\le \pi$.
\subsection{Two-mode representation of the aperiodic part: $n=2$}
For $n=2$, we have eight unknowns. So, let us consider the jumps in the derivatives of $u(x)$ at the endpoints $J_0, J_1, J_2, \hdots, J_7$ which yields the following set of equations,
\begin{gather}
 2\begin{bmatrix} 
\us_1 \sin(\ks_1\pi) \\ 
\us_2 \sin(\ks_2\pi) 
 \end{bmatrix}
 = 
 \begin{pmatrix}
  1		&	1 \\
  -\ks_1^2 & -\ks_2^2 
 \end{pmatrix}^{-1}
  \begin{bmatrix}
   J_0 \\ J_2
   \end{bmatrix}
 = 
 \begin{pmatrix}
  \ks_1^4 & \ks_2^4 \\
  -\ks_1^6 & -\ks_2^6  
 \end{pmatrix}^{-1}
  \begin{bmatrix}
   J_4 \\ J_6
   \end{bmatrix},
   \label{eqn:u2k2_sin}
\end{gather}
\begin{gather}
 -2\begin{bmatrix} 
\uc_1 \kc_1\sin(\kc_1\pi) \\ 
\uc_2 \kc_2\sin(\kc_2\pi) 
 \end{bmatrix}
 = 
 \begin{pmatrix}
  1		&	1	 \\
  -\kc_1^2 & -\kc_2^2 
 \end{pmatrix}^{-1}
  \begin{bmatrix}
   J_1 \\ J_3
   \end{bmatrix}
 = 
 \begin{pmatrix}
  \kc_1^4 & \kc_2^4 \\
  -\kc_1^6 & -\kc_2^6 
 \end{pmatrix}^{-1}
  \begin{bmatrix}
   J_5 \\ J_7
   \end{bmatrix}. 
   \label{eqn:u2k2_cos}
\end{gather}
Expanding on both sides of the latter equality in Eq. \eqref{eqn:u2k2_sin} yields,
\iftrue
\begin{gather}
 \begin{pmatrix}
  1		&	1 \\
  -\ks_1^2 & -\ks_2^2 
 \end{pmatrix}^{-1}
  \begin{bmatrix}
   J_0 \\ J_2
   \end{bmatrix}
 = 
 \left[
\begin{pmatrix}
  1		&	1 \\
  -\ks_1^2 & -\ks_2^2 
 \end{pmatrix}
 \begin{pmatrix}
  \ks_1^4 & 0 \\
  0 & \ks_2^4  
 \end{pmatrix}
 \right]^{-1}
  \begin{bmatrix}
   J_4 \\ J_6
   \end{bmatrix},
\end{gather}
\begin{gather}
\begin{pmatrix}
  1		&	1 \\
  -\ks_1^2 & -\ks_2^2 
 \end{pmatrix}
 \begin{pmatrix}
  \ks_1^4 & 0 \\
  0 & \ks_2^4  
 \end{pmatrix}
 \begin{pmatrix}
  1		&	1 \\
  -\ks_1^2 & -\ks_2^2 
 \end{pmatrix}^{-1}
  \begin{bmatrix}
   J_0 \\ J_2
   \end{bmatrix}
 = 
  \begin{bmatrix}
   J_4 \\ J_6
   \end{bmatrix},
\end{gather}
\fi
\begin{gather}
 \begin{pmatrix}
  -\ks_1^2 \ks_2^2		&	-\left(\ks_1^2+\ks_2^2 \right) \\
  \ks_1^2 \ks_2^2 \left(\ks_1^2+\ks_2^2 \right)	 	& \left(\ks_1^4+\ks_2^4+\ks_1^2 \ks_2^2 \right)
 \end{pmatrix}
  \begin{bmatrix}
   J_0 \\ J_2
   \end{bmatrix}
 = 
  \begin{bmatrix}
   J_4 \\ J_6
   \end{bmatrix}. \label{eqn:two_mode_proof}
\end{gather}
Elimination of $J_0$ from the bottom row of Eq. \eqref{eqn:two_mode_proof} yields the following linear system for the elementary symmetric polynomials in $\ks_1^2, \ks_2^2$ defined as $\es_1\equiv \ks_1^2 + \ks_2^2 , \es_2 \equiv \ks_1^2 \ks_2^2 $,
%\begin{minipage}{\linewidth}
%\begin{gather}
% \begin{pmatrix} 
% J_0 & J_2  \\ 
% J_2 & J_4  
% \end{pmatrix}
%  \begin{bmatrix}
%   \ks_1^2 \ks_2^2 
%   \\ \ks_1^2 + \ks_2^2 
%   \end{bmatrix}
% =
% -
%  \begin{bmatrix}
%   J_4 \\ J_6
%   \end{bmatrix}. 
%\end{gather}
\begin{gather}
 \begin{bmatrix} 
 J_0 ~\left(\ks_1^2 \ks_2^2 \right) + J_2 ~\left(\ks_1^2 + \ks_2^2\right)  \\ 
 J_2 ~\left(\ks_1^2 \ks_2^2 \right) + J_4 ~\left(\ks_1^2 + \ks_2^2\right)  
 \end{bmatrix}
 =
 -
  \begin{bmatrix}
   J_4 \\ J_6
   \end{bmatrix} 
\quad \Longrightarrow \quad
  \begin{bmatrix}
   \es_2
   \\ \es_1
   \end{bmatrix}
 =
 -
 \begin{pmatrix} 
 J_0 & J_2  \\ 
 J_2 & J_4  
 \end{pmatrix}^{\dagger}
  \begin{bmatrix}
   J_4 \\ J_6
   \end{bmatrix}. 
   \label{eqn:es2}
\end{gather}
%\end{minipage}
Corresponding to Eq. \eqref{eqn:u2k2_cos}, we obtain a similar system for the elementary polynomials in $\kc_1^2, \kc_2^2$ defined as $\ec_1\equiv \kc_1^2 + \kc_2^2 , \ec_2 \equiv \kc_1^2 \kc_2^2 $ and is given below,
\begin{gather}
  \begin{bmatrix}
   \ec_2
   \\ \ec_1
   \end{bmatrix}
 =
 -
 \begin{pmatrix} 
 J_1 & J_3  \\ 
 J_3 & J_5  
 \end{pmatrix}^{\dagger}
  \begin{bmatrix}
   J_5 \\ J_7
   \end{bmatrix}. 
   \label{eqn:ec2}
\end{gather}
By solving Eqs. \eqref{eqn:u2k2_sin}-\eqref{eqn:u2k2_cos} for $\es_1, \es_2, \ec_1, \ec_2$, we have
\begin{eqnarray*}
\es_1 \equiv \frac{J_0 J_6 - J_2 J_4}{J_2^2 - J_0 J_4},\quad \quad \es_2 \equiv \frac{J_4^2 - J_2 J_6}{J_2^2 - J_0 J_4}, \quad \quad \ec_1 \equiv \frac{J_1 J_7 - J_3 J_5}{J_3^2 - J_1 J_5},\quad \quad \ec_2 \equiv \frac{J_5^2 - J_3 J_7}{J_3^2 - J_1 J_5}.
\end{eqnarray*}
Now, $\ks_j^2$ and $\kc_j^2$ are the roots of the polynomials given below, and the amplitudes $\us_j, \uc_j$ are obtained from Eqs. \eqref{eqn:u2k2_sin}-\eqref{eqn:u2k2_cos},
$$
\ks^4 - \es_1 \ks^2 + \es_2 = 0,  ~~ \kc^4 - \ec_1 \kc^2 + \ec_2 = 0,$$ 
hence 
$$\ks^2 = \frac{\es_1\pm\sqrt{\es_1^2-4\es_2}}{2},  ~ \kc^2 = \frac{\ec_1\pm\sqrt{\ec_1^2-4\ec_2}}{2}.$$
And
%\end{flalign}
\begin{eqnarray*}
&\us_1 = \frac{\ks_2^2 J_0 + J_2}{2(\ks_2^2-\ks_1^2)\sin(\ks_1\pi)}, ~ &\us_2 = \frac{\ks_1^2 J_0 + J_2}{2(\ks_1^2-\ks_2^2)\sin(\ks_2\pi)}, \\
&\uc_1 = -\frac{\kc_2^2 J_1 + J_3}{2\kc_1 (\kc_2^2-\kc_1^2)\sin(\kc_1\pi)}, ~ &\uc_2 = -\frac{\kc_1^2 J_1 + J_3}{2\kc_2 (\kc_1^2-\kc_2^2)\sin(\kc_2\pi)}.
\end{eqnarray*}
Above expressions satisfy $[u^{(m)}]=[u_a^{(m)}]$ and thereby, $[u_{p}^{(m)}]=0$ for $m=0,1,2,\hdots,7$. This means that the periodic part is smooth up to the seventh order derivative, i.e., $C^7$ continuous.

\subsection{Three mode representation of the aperiodic part: $n=3$}
With $n=3$, we have 12 unknowns. So, let us consider the jumps in the derivatives of $u(x)$ at the endpoints $J_0, J_1, J_2, \hdots, J_{11}$ which yields the following set of equations,

%\begin{minipage}{0.45\linewidth}
%\begin{gather}
% 2
%  \begin{pmatrix}
%  1		&	1	&	1 \\
%  -\ks_1^2 & -\ks_2^2 & -\ks_3^2 \\
%  \ks_1^4 & \ks_2^4 & \ks_3^4 \\
%  -\ks_1^6 & -\ks_2^6 & -\ks_3^6 \\
%  \ks_1^8 & \ks_2^8 & \ks_3^8 \\
%  -\ks_1^{10} & -\ks_2^{10} & -\ks_3^{10} 
% \end{pmatrix}
%\begin{pmatrix} 
%\us_1 \sin(\ks_1\pi) \\ 
%\us_2 \sin(\ks_2\pi) \\
%\us_3 \sin(\ks_3\pi) 
% \end{pmatrix}
% = 
%  \begin{bmatrix}
%   J_0 \\ J_2\\J_4
%\\   J_6 \\ J_8\\J_{10}
%   \end{bmatrix}, \nonumber
%\end{gather}
%\end{minipage}
%\begin{minipage}{0.45\linewidth}
%\begin{gather}
% -2
% \begin{pmatrix}
%  1		&	1	&	1 \\
%  -\kc_1^2 & -\kc_2^2 & -\kc_3^2 \\
%  \kc_1^4 & \kc_2^4 & \kc_3^4 \\
%  -\kc_1^6 & -\kc_2^6 & -\kc_3^6 \\
%  \kc_1^8 & \kc_2^8 & \kc_3^8 \\
%  -\kc_1^{10} & -\kc_2^{10} & -\kc_3^{10}
% \end{pmatrix}
% \begin{pmatrix} 
%\uc_1 \kc_1\sin(\kc_1\pi) \\ 
%\uc_2 \kc_2\sin(\kc_2\pi) \\
%\uc_3 \kc_3\sin(\kc_3\pi) 
% \end{pmatrix}
% = 
%  \begin{bmatrix}
%   J_1 \\ J_3\\J_5
%\\ J_7 \\ J_9\\J_{11}
%   \end{bmatrix}. \nonumber
%\end{gather}
%\end{minipage}
%\vspace{0.1in}
%\\
%Above equations can be rearranged to separate the unknowns $\usj$, $\ksj$ and $\ucj$, $\kcj$ as below,

\begin{gather}
 2\begin{bmatrix} 
\us_1 \sin(\ks_1\pi) \\ 
\us_2 \sin(\ks_2\pi) \\
\us_3 \sin(\ks_3\pi) 
 \end{bmatrix}
 = 
 \begin{pmatrix}
  1		&	1	&	1 \\
  -\ks_1^2 & -\ks_2^2 & -\ks_3^2 \\
  \ks_1^4 & \ks_2^4 & \ks_3^4 \\
 \end{pmatrix}^{-1}
  \begin{bmatrix}
   J_0 \\ J_2\\J_4
   \end{bmatrix}
 = 
 \begin{pmatrix}
  -\ks_1^6 & -\ks_2^6 & -\ks_3^6 \\
  \ks_1^8 & \ks_2^8 & \ks_3^8 \\
  -\ks_1^{10} & -\ks_2^{10} & -\ks_3^{10} \\
 \end{pmatrix}^{-1}
  \begin{bmatrix}
   J_6 \\ J_8\\J_{10}
   \end{bmatrix}. 
   \label{eqn:u3k3_sin}
\end{gather}
\begin{gather}
 -2\begin{bmatrix} 
\uc_1 \kc_1\sin(\kc_1\pi) \\ 
\uc_2 \kc_2\sin(\kc_2\pi) \\
\uc_3 \kc_3\sin(\kc_3\pi) 
 \end{bmatrix}
 = 
 \begin{pmatrix}
  1		&	1	&	1 \\
  -\kc_1^2 & -\kc_2^2 & -\kc_3^2 \\
  \kc_1^4 & \kc_2^4 & \kc_3^4 \\
 \end{pmatrix}^{-1}
  \begin{bmatrix}
   J_1 \\ J_3\\J_5
   \end{bmatrix}
 = 
 \begin{pmatrix}
  -\kc_1^6 & -\kc_2^6 & -\kc_3^6 \\
  \kc_1^8 & \kc_2^8 & \kc_3^8 \\
  -\kc_1^{10} & -\kc_2^{10} & -\kc_3^{10} \\
 \end{pmatrix}^{-1}
  \begin{bmatrix}
   J_7 \\ J_9\\J_{11}
   \end{bmatrix}. 
   \label{eqn:u3k3_cos}
\end{gather}

Note that the second and third parts in Eqs. \eqref{eqn:u3k3_sin}-\eqref{eqn:u3k3_cos} involving the unknowns $\ksj$ and $\kcj$ are same except that the jumps $J_0, J_2,\hdots, J_{10}$ are replaced by $J_1, J_3,\hdots, J_{11}$. So, the procedures for solving $\ksj$ and $\kcj$ are similar. 

Assuming full rank, i.e., $\ks_1\neq\ks_2\neq \ks_3$, and expanding on both sides of the latter equality in Eqs. \eqref{eqn:u3k3_sin} yields,
\begin{gather}
\begin{pmatrix}
  1		&	1	&	1 \\
  -\ks_1^2 & -\ks_2^2 & -\ks_3^2 \\
  \ks_1^4 & \ks_2^4 & \ks_3^4 \\
 \end{pmatrix}^{-1}
  \begin{bmatrix}
   J_0 \\ J_2\\J_4
   \end{bmatrix}
 = 
-
\left[
\begin{pmatrix}
  1		&	1	&	1 \\
  -\ks_1^2 & -\ks_2^2 & -\ks_3^2 \\
  \ks_1^4 & \ks_2^4 & \ks_3^4 \\
 \end{pmatrix}
 \begin{pmatrix}
  \ks_1^6 & 0 & 0 \\
  0 & \ks_2^6 & 0 \\
  0 & 0 & \ks_3^{6} 
 \end{pmatrix}
 \right]^{-1}
  \begin{bmatrix}
   J_6 \\ J_8\\J_{10},
   \end{bmatrix},
\end{gather}
\begin{gather}
-
\begin{pmatrix}
  1		&	1	&	1 \\
  -\ks_1^2 & -\ks_2^2 & -\ks_3^2 \\
  \ks_1^4 & \ks_2^4 & \ks_3^4 \\
 \end{pmatrix}
 \begin{pmatrix}
  \ks_1^6 & 0 & 0 \\
  0 & \ks_2^6 & 0 \\
  0 & 0 & \ks_3^{6} 
 \end{pmatrix}
\begin{pmatrix}
  1		&	1	&	1 \\
  -\ks_1^2 & -\ks_2^2 & -\ks_3^2 \\
  \ks_1^4 & \ks_2^4 & \ks_3^4 \\
 \end{pmatrix}^{-1}
  \begin{bmatrix}
   J_0 \\ J_2\\J_4
   \end{bmatrix}
 = 
  \begin{bmatrix}
   J_6 \\ J_8\\J_{10},
   \end{bmatrix},
\end{gather}
\iffalse
\begin{gather*}
-
\frac{1}{\ks_1^4 (\ks_2^2-\ks_3^2) + \ks_2^4(\ks_3^2-\ks_1^2) + \ks_3^4(\ks_1^2-\ks_2^2)}
\begin{pmatrix}
  1		&	1	&	1 \\
  -\ks_1^2 & -\ks_2^2 & -\ks_3^2 \\
  \ks_1^4 & \ks_2^4 & \ks_3^4 \\
 \end{pmatrix}
 \begin{pmatrix}
  \ks_1^6 & 0 & 0 \\
  0 & \ks_2^6 & 0 \\
  0 & 0 & \ks_3^{6} 
 \end{pmatrix}
 \begin{pmatrix}
  \ks_2^2 \ks_3^2 \left( \ks_2^2-\ks_3^2)\right)		&	\left(\ks_2^4-\ks_3^4 \right) & \left(\ks_2^2-\ks_3^2 \right) \\
  \ks_3^2 \ks_1^2 \left( \ks_3^2-\ks_1^2)\right)		&	\left(\ks_3^4-\ks_1^4 \right) & \left(\ks_3^2-\ks_1^2 \right) \\
  \ks_1^2 \ks_2^2 \left( \ks_1^2-\ks_2^2)\right)		&	\left(\ks_1^4-\ks_2^4 \right) & \left(\ks_1^2-\ks_2^2 \right)    \end{pmatrix}
  \begin{bmatrix}
   J_0 \\ J_2\\J_4
   \end{bmatrix}
 = 
  \begin{bmatrix}
   J_6 \\ J_8 \\ J_{10}
   \end{bmatrix}. 
   \label{eqn:es3}
\end{gather*}
\fi
%\frac{1}{\ks_1^6 \ks_2^6\ks_3^6}
\begin{gather*}
-
 \begin{pmatrix}
  \ks_1^6 & 0 & 0 \\
  0 & \ks_2^6 & 0 \\
  0 & 0 & \ks_3^{6} 
 \end{pmatrix}
 \begin{pmatrix}
  \ks_2^2 \ks_3^2 \left( \ks_2^2-\ks_3^2)\right)		&	\left(\ks_2^4-\ks_3^4 \right) & \left(\ks_2^2-\ks_3^2 \right) \\
  \ks_3^2 \ks_1^2 \left( \ks_3^2-\ks_1^2)\right)		&	\left(\ks_3^4-\ks_1^4 \right) & \left(\ks_3^2-\ks_1^2 \right) \\
  \ks_1^2 \ks_2^2 \left( \ks_1^2-\ks_2^2)\right)		&	\left(\ks_1^4-\ks_2^4 \right) & \left(\ks_1^2-\ks_2^2 \right)    \end{pmatrix}
  \begin{bmatrix}
   J_0 \\ J_2\\J_4
   \end{bmatrix}
 = 
 \begin{pmatrix}
  \ks_2^2 \ks_3^2 \left( \ks_2^2-\ks_3^2)\right)		&	\left(\ks_2^4-\ks_3^4 \right) & \left(\ks_2^2-\ks_3^2 \right) \\
  \ks_3^2 \ks_1^2 \left( \ks_3^2-\ks_1^2)\right)		&	\left(\ks_3^4-\ks_1^4 \right) & \left(\ks_3^2-\ks_1^2 \right) \\
  \ks_1^2 \ks_2^2 \left( \ks_1^2-\ks_2^2)\right)		&	\left(\ks_1^4-\ks_2^4 \right) & \left(\ks_1^2-\ks_2^2 \right)    \end{pmatrix}
  \begin{bmatrix}
   J_6 \\ J_8 \\ J_{10}
   \end{bmatrix}. 
\end{gather*}
\begin{gather*}
-
 \begin{pmatrix}
  \ks_1^6 & 0 & 0 \\
  0 & \ks_2^6 & 0 \\
  0 & 0 & \ks_3^{6} 
 \end{pmatrix}
 \begin{pmatrix}
  \ks_2^2 \ks_3^2 	&	\left(\ks_2^2+\ks_3^2 \right) & 1 \\
  \ks_3^2 \ks_1^2 		&	\left(\ks_3^2+\ks_1^2 \right) & 1\\
  \ks_1^2 \ks_2^2 		&	\left(\ks_1^2+\ks_2^2 \right) & 1    \end{pmatrix}
  \begin{bmatrix}
   J_0 \\ J_2\\J_4
   \end{bmatrix}
 = 
 \begin{pmatrix}
  \ks_2^2 \ks_3^2 	&	\left(\ks_2^2+\ks_3^2 \right) & 1 \\
  \ks_3^2 \ks_1^2 		&	\left(\ks_3^2+\ks_1^2 \right) & 1\\
  \ks_1^2 \ks_2^2 		&	\left(\ks_1^2+\ks_2^2 \right) & 1    \end{pmatrix}
  \begin{bmatrix}
   J_6 \\ J_8 \\ J_{10}
   \end{bmatrix}. 
\end{gather*}

\begin{eqnarray}
J_0 \left(\ks_1^6 \ks_2^2 \ks_3^2\right) + \left(\ks_2^2+\ks_3^2\right) \left[J_2 \ks_1^6 + J_8\right] + J_4 \ks_1^6 + J_6 \left(\ks_2^2 \ks_3^2\right) + J_{10} &=& 0,\label{eqn:es3.1}\\
J_0 \left(\ks_1^2 \ks_2^6 \ks_3^2\right) + \left(\ks_1^2+\ks_3^2\right) \left[J_2 \ks_2^6 + J_8\right] + J_4 \ks_2^6 + J_6 \left(\ks_1^2 \ks_3^2\right) + J_{10} &=& 0,\label{eqn:es3.2}\\
J_0 \left(\ks_1^2 \ks_2^2 \ks_3^6\right) + \left(\ks_1^2+\ks_2^2\right) \left[J_2 \ks_3^6 + J_8\right] + J_4 \ks_3^6 + J_6 \left(\ks_1^2 \ks_2^2\right) + J_{10} &=& 0. \label{eqn:es3.3}
\end{eqnarray}
Elimination of the pairs $(J_8, J_{10})$, $(J_0, J_{10})$, $(J_0, J_2)$ from Eqs. \eqref{eqn:es3.1}-\eqref{eqn:es3.3} yields
%After some tedious manipulation of the set of equations involving the second and the third parts in Eq. \eqref{eqn:u3k3_sin}, we obtain 
the following set of equations involving the elementary symmetric polynomials ($\es_j$) in $\ks_1^2, \ks_2^2, \ks_3^2$ defined as $\es_1 \equiv \ks_1^2 + \ks_2^2 + \ks_3^2, 
\es_2 \equiv \ks_1^2 \ks_2^2 + \ks_2^2 \ks_3^2 + \ks_1^2 \ks_3^2,
\es_3 \equiv \ks_1^2 \ks_2^2 \ks_3^2$,
\begin{gather}
 \begin{pmatrix} 
 J_0 & J_2 & J_4 \\ 
 J_2 & J_4 & J_6 \\
 J_4 & J_6 & J_8 
 \end{pmatrix}
  \begin{bmatrix}
   \ks_1^2 \ks_2^2 \ks_3^2 
   \\ \ks_1^2 \ks_2^2 + \ks_2^2 \ks_3^2 + \ks_1^2 \ks_3^2
   \\ \ks_1^2 + \ks_2^2 + \ks_3^2
   \end{bmatrix}
 =
 -
  \begin{bmatrix}
   J_6 \\ J_8\\J_{10}
   \end{bmatrix}
\quad \Longrightarrow \quad
  \begin{bmatrix}
   \es_3 \\ \es_2 \\ \es_1
   \end{bmatrix}
 =
 -
 \begin{pmatrix} 
 J_0 & J_2 & J_4 \\ 
 J_2 & J_4 & J_6 \\
 J_4 & J_6 & J_8 
 \end{pmatrix}^{\dagger}
  \begin{bmatrix}
   J_6 \\ J_8\\J_{10}
  \end{bmatrix}. 
  \label{eqn:es_3}
\end{gather}
%where, $\es_j$ are defined as below,
%\begin{eqnarray*}
%\es_1 &\equiv& \ks_1^2 + \ks_2^2 + \ks_3^2,\\
%\es_2 &\equiv& \ks_1^2 \ks_2^2 + \ks_2^2 \ks_3^2 + \ks_1^2 \ks_3^2,\\
%\es_3 &\equiv& \ks_1^2 \ks_2^2 \ks_3^2.\\
%\end{eqnarray*}

Similarly, corresponding to Eq. \eqref{eqn:u3k3_cos}, we obtain the elementary symmetric polynomials $\ec_j$ in $\kc_1^2, \kc_2^2, \kc_3^2$ defined as $\ec_1 \equiv \kc_1^2 + \kc_2^2 + \kc_3^2, 
\ec_2 \equiv \kc_1^2 \kc_2^2 + \kc_2^2 \kc_3^2 + \kc_1^2 \kc_3^2,
\ec_3 \equiv \kc_1^2 \kc_2^2 \kc_3^2$.
\begin{gather}
  \begin{bmatrix}
   \ec_3 \\ \ec_2 \\ \ec_1
   \end{bmatrix}
 =
 -
 \begin{pmatrix} 
 J_1 & J_3 & J_5 \\ 
 J_3 & J_5 & J_7 \\
 J_5 & J_7 & J_9 
 \end{pmatrix}^{\dagger}
  \begin{bmatrix}
   J_7 \\ J_9\\J_{11}
   \end{bmatrix}.  
  \label{eqn:ec_3}   
\end{gather}
%\begin{eqnarray*}
%\ec_1 &\equiv& \kc_1^2 + \kc_2^2 + \kc_3^2, \\
%\ec_2 &\equiv& \kc_1^2 \kc_2^2 + \kc_2^2 \kc_3^2 + \kc_1^2 \kc_3^2 \\
%\ec_3 &\equiv& \kc_1^2 \kc_2^2 \kc_3^2.
%\end{eqnarray*}
Now, $\ks_1^2, \ks_2^2, \ks_3^2$ and $\kc_1^2, \kc_2^2, \kc_3^2$ are obtained by solving for the roots ($\lambda$) of the polynomials below:
\begin{equation*}
\lambda^3 - \es_1 \lambda^2 + \es_2 \lambda - \es_3 = 0 = \prod_{j=1}^3 (\lambda-\ksj^2), \quad \quad
\lambda^3 - \ec_1 \lambda^2 + \ec_2 \lambda - \ec_3 = 0 = \prod_{j=1}^3 (\lambda-\kcj^2).
\end{equation*}

Above expressions satisfy $[u^{(m)}]=[u_a^{(m)}]$ and thereby, $[u_{p}^{(m)}]=0$ for $m=0,1,2,\hdots,11$. This implies that the periodic part is smooth up to the eleventh order derivative, i.e., $C^{11}$ continuous.

\section*{Acknowledgments}
Narsimha Rapaka would like to dedicate this work to the Late Professor Ravi Samtaney.
\section*{Code and data availability}
All source code and datasets to reproduce our results are publicly available at \url{https://github.com/nrapaka/GFS}. The repository includes the GFS implementation, example scripts, and instructions for reproducing the experiments.

 % \biboptions{authoryear}
 % \bibliographystyle{elsarticle-num-names} 
 % \bibliographystyle{elsarticle-harv} 
  % \bibliographystyle{elsarticle-num}
\bibliographystyle{plainnat}
\bibliography{References_mcs}

\end{document}